\documentclass[12pt, a4paper]{article}
\usepackage[left=1in,right=1in,top=1.0in,bottom=1in]{geometry}

\usepackage{amsmath,amssymb,amsthm,amsfonts,afterpage,hyperref}
\usepackage{amscd,subeqnarray}

\usepackage{graphicx}
\usepackage{caption}
\usepackage{subcaption}
\usepackage{authblk}
\usepackage{tikz}

\usepackage{array}
\usepackage{wrapfig}
\usepackage{color}
\usepackage{ragged2e}
\usepackage{stmaryrd}
\usepackage{siunitx}
\usepackage{multirow}
\usepackage{xcolor,cancel}
\usepackage{boldfonts}
\usepackage{symbols}
\usepackage{footmisc}
\usepackage{float}
\usepackage{setspace}
\usepackage{bm}
\usepackage{booktabs}
\usepackage{tikz, xcolor}
\usepackage{soul}
\usepackage{pdfpages} 

\usepackage[T1]{fontenc}
\usepackage[utf8]{inputenc}

\usepackage[linesnumbered,ruled,vlined]{algorithm2e}
\SetKwInput{KwInput}{Input}                
\SetKwInput{KwOutput}{Output} 

\usepackage[sort, numbers]{natbib}
\setcitestyle{square}

\usetikzlibrary{shapes,arrows, positioning,calc}	
\usetikzlibrary{arrows,snakes,backgrounds}

\newcommand{\bfeps}{\boldsymbol{\epsilon}}		
\newcommand{\bfeta}{\boldsymbol{\eta}}

\newcommand{\Sym}{\text{Sym}}   			%
\newcommand{\ddiv}{\text{div}}     				%
\newcommand{\tr}{\text{tr}}       				%

\DeclareMathAlphabet{\mathpzc}{OT1}{pzc}{m}{it}

\newcommand{\bfn}{\boldsymbol{n}}	
	
\newcommand{\bfu}{\boldsymbol{u}}	
\newcommand{\bfb}{\boldsymbol{b}}	
\newcommand{\bfv}{\boldsymbol{v}}	

\newcommand{\bfF}{\boldsymbol{F}}	
\newcommand{\bfE}{\boldsymbol{E}}	
\newcommand{\bfC}{\boldsymbol{C}}
\newcommand{\bfB}{\boldsymbol{B}}	
\newcommand{\bfx}{\boldsymbol{x}}	
\newcommand{\bfX}{\boldsymbol{X}}	
\newcommand{\bfS}{\boldsymbol{S}}	
\newcommand{\bfT}{\boldsymbol{T}}		
\newcommand{\bfI}{\boldsymbol{I}}	 
\newcommand{\bfzero}{\boldsymbol{0}}

\newcommand{\bfA}{\boldsymbol{A}}

\newcommand{\bff}{\boldsymbol{f}}	
\newcommand{\bfg}{\boldsymbol{g}}	

\newtheorem{remark}{Remark}

\newtheorem{formulation}{Formulation}

\providecommand{\keywords}[1]
{
  \small	
  \textbf{\textit{Keywords---}} #1
}

\title{Preferential stiffness and the crack-tip fields of an elastic porous solid 
based on the density-dependent  moduli model}

\author[1]{Hyun C. Yoon\thanks{hyun.yoon@kigam.re.kr}}
\author[2]{S. M. Mallikarjunaiah\thanks{m.muddamallappa@tamucc.edu}\thanks{corresponding author}}
\author[3]{Dambaru Bhatta\thanks{dambaru.bhatta@utrgv.edu}}

\affil[1]{Marine Geology \& Energy Division,
Korea Institute of Geoscience and Mineral Resources,
124 Gwahak-ro,
Daejeon 34132, Republic of Korea}
\affil[2]{Department of Mathematics \& Statistics,
Texas A\&M University - Corpus Christi, 
Corpus Christi, Texas 78412-5825, USA}
\affil[3]{School of Mathematical \& Statistical Sciences,
The University of Texas - Rio Grande Valley, 
Edinburg, Texas 78539, USA}
\date{}

\begin{document}

\maketitle
	    
\begin{abstract}
In this paper, we study the preferential stiffness and the crack-tip fields for an elastic porous solid of which material properties are dependent upon the density. Such a description is necessary to describe the failure {that can be caused by damaged pores 
in many porous 
bodies} such as ceramics, concrete and human bones. 
{To that end,} we revisit a new class of implicit constitutive relations 
under the assumption of small deformation.   
 Although the constitutive relationship \textit{appears linear} in both the Cauchy stress and linearized strain, the governing equation bestowed from the balance of linear momentum results in a quasi-linear partial differential equation (PDE) system. 
{For the linearization and obtaining a sequence of elliptic PDEs, we propose the solution algorithm comprise a \textit{Newton's method} coupled with a bilinear continuous Galerkin-type finite elements for the discretization.} Our algorithm exhibits an optimal rate of convergence for a  manufactured solution.
 In the numerical experiments,  
 we set the 
 boundary value problems (BVPs) with edge crack {under different modes of loading (i.e., the pure mode-I, II, and the mixed-mode). From the numerical results, we find that the {density-dependent  
 moduli model}} describes diverse phenomena that are not captured within the framework of classical linearized elasticity. {In particular,} 
 {numerical solutions} clearly indicate that the  {nonlinear \textit{modeling}} parameter  {depending on its sign and magnitude}  
{can control} preferential mechanical stiffness 
along with the change of volumetric strain; 
{larger the parameter is in the positive value}, the responses are such that 
the strength of porous solid gets weaker against the tensile loading while stiffer against the in-plane shear (or compressive) loading, 
which is vice versa for the negative value of it. 
 The modeling framework of the density-dependent material moduli   proposed in this study
 can provide a mathematical and computational foundation to further 
model the quasi-static and dynamic evolution of cracks and 
many other multi-physics applications such as the  
{fluid flow or heat transfer} in porous media.  
 \end{abstract}

\noindent \keywords{Density-dependent moduli, Preferential stiffness, Implicit constitutive relation, Porous solid, Finite element method}

\section{Introduction}\label{sec:Intro}
Mechanical stiffness of a porous solid is dependent not only on the mechanical property of solid grain but also on the bulk skeleton composed of connected or disconnected pores (or microstructures). 
For example, the shear stiffness of sandy soil 
can be lost and act like liquid when soil liquefaction occurs 
\cite{KiyT2013}.  It is also well known that the mechanical properties of metal matrix composite or alloy are largely dependent on the defects, i.e., porosity \cite{HarR2013}. 
{Furthermore, the pore space can be partially or fully saturated with fluid (e.g., gas or liquid) inside, which can induce the poromechanical effects \cite{CouO2004,YooH2021-jngse}.
In petroleum engineering and rock mechanics, 
the hydraulic fracturing \cite{WuY2018,ZobM2019} (i.e., fluid injection with high pressure gradient) is a well-known stimulation technology for the development of shale gas, 
which aims to generate the tensile fractures  
in the rock.}  {The success of hydraulic fracturing  is highly dependent on the rock mechanical properties, as 
it is {much more} efficient for  sand which is stiffer than soft clay or mud \cite{WuY2018,YewC2014} for the tensile failure. {Considering the heterogeneity 
in the composites, thereby simple averaging of rock properties might mislead, and rigorous mathematical formulation (e.g., the multiple porosity model \cite{KimJ2012-ijnme,YooH2021-jngse}) for upscaling is required.}} 

Material properties and mechanical strength
are also reliant on
certain modes of loading preferentially. 
In general, ductile materials are known to possess approximately equal strength in tension and compression, but weak in shear. Meanwhile, brittle materials are generally weak in tension, and 
its tensile strength is known to be only one tenth of compressive strength \cite{CalW2018}. 
For instance, the two-dimensional (2D) material \cite{GuoY2021} 
(e.g., the graphene) is considered to be one of the strongest materials able to withstand tensional loading due to its 2D nature of the structure. 
Concrete or ceramics is weak under tension but strong against compression, thus composite (or reinforced) material can be synthesized with them to have higher tensile strength. 
Rock is another exemplary material such that it is stiffer against the compressive loading 
but weaker against the tensile loading.  
{For the hydraulic fracturing,} some experimental results also substantiate that the shear failure may occur first before the tensile failure around the crack-tip \cite{vanDam2001}.
Therefore, it
is necessary to consider 
the realistic preferential regime related to the stiffness of a material.
Regarding the preferential stiffness,
we particularly pay attention to the density change of the porous material;  
microstructures of material including pore space are required to bond one another in order to transmit the tensile loading 
\cite{CalW2018}. 
Meanwhile, the fracture propagation is demanding under a compressive loading, as transverse cracks can emerge and tend to congregate, resulting in the increased density of material.
In essence, the preferential strength can 
be modeled such that the material moduli   
are dependent upon the density of material or
the  change in volumetric strain: the dilation versus compaction.
{In addition, this density-dependent stiffness concept is 
physically reasonable in that elastic properties of a porous material are known to have nonlinear relations along with the porosity value \cite{KovJ1999,ManA2013}, even though it is not straightforward to model these in continuum scale.} 

The conventional practice of modeling porous solid within the elastic and infinitesimal strain regime is by using the linearized elasticity, wherein the Cauchy stress is expressed as a function of deformation gradient, density, and material points, in which the material moduli are generally taken to be constants. The linearization would thus result in classical relationship, which has been successfully used in various applications. {However, the classical linearized theory of elasticity is incapable of accommodating some realistic responses as aforementioned, such as the ones found in many metallic alloys (see \cite{saito2003multifunctional,zhang2009fatigue}) and concrete (see \cite{grasley2015model}) that clearly exhibit nonlinear mechanical responses well within the range of ``small strains''.} {Recently,} Rajagopal has established a new class of constitutive relations for the description of 
non-dissipative elastic bodies  
(for more details about `` novel'' elastic constitute relationships, see \cite{rajagopal2003,rajagopal2007elasticity,rajagopal2014nonlinear,rajagopal2018note}). These implicit constitutive relations have been explored in many studies to describe the state of stress-strain in the neighborhood of crack-tip \cite{rajagopal2011modeling,Mallikarjunaiah2015,kulvait2013,kulvait2019,HCY_SMM_MMS2022,gou2015modeling,MalliPhD2015,ortiz2014numerical,ortiz2012,bulivcek2014elastic,rodriguez2021stretch}, for the behavior of thermoelastic bodies \cite{yoon2022finite,bustamante2017implicit}, the response of viscoelastic materials \cite{muliana2013new,erbay2015traveling,itou2019crack,sengul2021viscoelasticity,erbay2020thermodynamically,erbay2020local,csengul2021nonlinear}, and the quasi-static crack evolution \cite{yoon2021quasi,lee2022finite}.

The objective of this study to investigate {elastic porous solids with preferential stiffness with} an explicit 
{type} of the density-dependent moduli model, derived from the implicit constitutive relation {following} \cite{rajagopal2021a,rajagopal2021b,rajagopal2022implicit}. 
To this end, we assume an isotropic and homogeneous porous solid with crack inside, where 
the mechanical regime is under small strain and pure elasticity before any failure. Particularly for the stiffness variations, we focus on the crack (or damaged pores) and its tip under 
different modes of loadings: the mode-I, II, and mixed-mode.   
Instead of an explicit porosity, 
{i.e., one of the state variables based on the mixture theory \cite{CouO2004,SteH2019,RajK1996} that can also be saturated with fluid inside,} 
a pure solid matrix with implicit porosity is considered  
through the change of 
volumetric strain. 
Meanwhile, another goal of this contribution
is to weigh up the applicability of density-dependent material moduli model for the quasi-static or dynamic propagation of network of cracks under mechanical/thermal loading.
Thus, we employ a computational model that is physically and mathematically
consistent with  
the framework developed in \cite{yoon2021quasi,lee2022finite},  
an universal approach based on the Newton's method and standard finite element method (FEM).   
This approach yields a numerically stable and efficient algorithm against the serious nonlinear problem of our interest. Our iterative algorithm has shown to display an optimal order convergence for a BVP with manufactured solution.

In the numerical experiments 
comparing this nonlinear model with the conventional linearized elasticity model, we identify distinct variations of stress and strain distributions, especially near the crack-tip area under different modes of loading, i.e., parallel and perpendicular to the crack. We also compare the stress intensity factor, drained bulk modulus, and volumetric strain change to elucidate the preferential stiffness of the elastic porous solid. We find that the {nonlinear modeling} parameter with its magnitude and sign for the density-dependent moduli of elastic porous solid can appropriately describe the preferential stiffness related to the principal direction of the loadings. 
With larger positive value of the  
parameter,
the strength of porous solid gets weaker against the tensile loading, while it becomes stiffer against the in-plane shear (or compressive) loading, which is vice versa for the negative value of it.
Upon the efficacy of the density-dependent material moduli model and its modeling framework, 
current study 
can further be expanded to the topics such as the static fracture 
evolution, fluid transport, hydraulic fracture, or heat transfer in porous media. 

\section{Formulation {of the density-dependent material moduli model}}\label{sec:Formul}
The current investigation is prompted by a recent study (see  \cite{rajagopal2021a,rajagopal2021b,rajagopal2022implicit} for more details) on developing new nonlinear constitutive relationships between the linearized strain and the Cauchy stress that are worthwhile in characterizing the cracks  
or damaged pores in porous solids such as rocks and concrete.

\subsection{Basic notations}
In this section, we provide a brief introduction to the geometrical and mechanical description of the elastic porous solid body, highlighting its mass balance. Let $\Omega \subset \mathbb{R}^2$ be a closed and bounded Lipschitz domain in the Euclidean space that represent the material body under consideration in the reference configuration.  Let $\partial \Omega$ be the Lipschitz  continuous boundary and $\bfeta=\left(\eta_1, \, \eta_2\right)$ be the outward unit normal. Further, we consider a boundary part consists of two disjoint parts such that:
\[
 \partial \Omega = \overline{ \; \Gamma_{D}} \cup  \overline{\; \Gamma_{N}} \quad \mbox{and} \quad  \Gamma_{D} \cap   \Gamma_{N} = \emptyset,
\]
where the partition consists of a Neumann boundary $\Gamma_{N}$ and a nonempty Dirichlet boundary 
$\Gamma_{D}$. Let $\Gamma_c \subseteq \Omega$ be a geometrical boundary of an interface and we assume that the interface $\Gamma_c$ is completely contained inside of $\Omega$, not on the boundary $\partial \Omega$. 
Let $\bfX = \left(X_1, \,  X_2 \right)$  and 
$\bfx = \left(x_1, \, x_2 \right)$ denote typical points in the reference and deformed configurations of the body.  Let $\bfu \colon \Omega \to \mathbb{R}^2$ denote the displacement field, and 
 \begin{equation}
 \bfu = \bfx - \bfX.  
 \end{equation}
 In the rest of this paper, we use the usual notations of \textit{Lebesgue} and \textit{Sobolev spaces} \cite{ciarlet2002finite,evans1998partial}. We denote $L^{p}(\Omega)$ as the space of all \textit{Lebesgue integrable functions} with $p\in[1, \infty)$. In particular,  when $p=2$, $L^{2}(\Omega)$ denote the space of all square integrable (Lebesgue) functions on $\Omega$ together with its inner product  $\left( v, \; w \right):= \int_{\Omega} v \,w \; d\bfx$ and the norm $\| v \| := \| v \|_{L^2(\Omega)} = \left(v, \, v \right)^{1/2}$. Let $C^{m}(\Omega), \; m \in \mathbb{N}_0$ denote the linear space of continuous functions on $\Omega$ and let $H^{1}(\Omega)$ denote the standard Sobolev space:
 \begin{equation}
 H^{1}(\Omega) := W^{1,2} = \left\{ v \in L^{2}(\Omega) \; ; \; \partial_j \, v \in L^{2}(\Omega) \; j \in \{ 1, 2, \ldots n \} \right\}, 
 \end{equation}
with the inner (scalar) product and the norm defined, respectively, as:
\begin{subequations}
\begin{align}
\left( f, \; g \right)_{H^1} &:= \left( f, \; g \right) + \sum_{j=1}^{n} \left( \partial_j f, \; \partial_j g \right), \\
\| f \|_{H^1} &:= \left(  \| f^2 \|^2 + \sum_{j=1}^{n} \| \partial_j f^2  \|^2 \right)^{1/2}.
\end{align}
\end{subequations}
The closure of $C_0^{\infty}(\Omega)$ by the norm of $H^1(\Omega)$ is denoted as $H^{1}_0(\Omega)$, i.e.,
 \begin{equation}\label{def-H01}
 H_0^{1}(\Omega) = \overline{C_0^{\infty}(\Omega)}^{\| \cdot \|_{H^1}}.
 \end{equation}
The space of displacements satisfying the homogeneous and non-homogeneous  Dirichlet boundary conditions are defined as:
\begin{align}\label{test_fun_space}
\widehat{V}_{\bfzero} &:= \left\{ \bfu \in \left( H^{1}(\Omega)\right)^2 \colon \bfu=\bfzero \; \mbox{on} \; \Gamma_D \right\}, \\
\widehat{V}_{\bfg} &:= \left\{ \bfu \in \left( H^{1}(\Omega)\right)^2 \colon \bfu=\bfg \; \mbox{on} \; \Gamma_D \right\}.
\end{align}
Let $\Sym(\mathbb{R}^{2 \times 2})$ is the space of $2 \times 2$ symmetric tensors equipped with the inner product $\bfA \colon \bfB = \sum_{i, \, j=1}^d \, \bfA_{ij} \, \bfB_{ij}$ and for all $\bfA = (\bfA)_{ij}$ and $\bfB =(\bfB)_{ij}$ in $\Sym(\mathbb{R}^{2 \times 2})$, {the associated norm}  {$\| \bfA \| = \sqrt{\bfA \colon \bfA}$}. Further, let $\bfF \colon \Omega \to \mathbb{R}^{2 \times 2}$ denote the deformation gradient,  $\bfC \colon \Omega \to \mathbb{R}^{2 \times 2}$ denote the right Cauchy-Green stretch tensor, $\bfB \colon \Omega \to \mathbb{R}^{2 \times 2}$ denote the left Cauchy-Green stretch tensor, $\bfE \colon \Omega \to \mathbb{R}^{2 \times 2}$ denote the Lagrange strain, $\bfeps \colon \Omega \to  \Sym(\mathbb{R}^{2 \times 2})$ denote the  linearized strain tensor, respectively defined as:
\begin{subequations}
\begin{align}
\bfF &:= \nabla_r \bfx = \bfI + \nabla \bfu, \\
 \bfB &:=\bfF\bfF^{\mathrm{T}}, \; \bfC :=\bfF^{\mathrm{T}}\bfF,  \; \bfE := \dfrac{1}{2} \left( \bfC - \bfI \right),\\
\bfeps(\bfu) &:= \dfrac{1}{2} \left( \nabla \bfu + \nabla \bfu^{\mathrm{T}}\right),
\end{align}
\end{subequations}
where $\left( \cdot \right)^{\mathrm{T}}$ denotes the \textit{transpose} operator for the second-order tensors, $\bfI$ is the
two-dimensional identity tensor, $\nabla_r$ and $\nabla$ are the gradient operators in the reference and current configuration, respectively.
Under the standard assumption of linearized elasticity, we have the following 
\begin{equation}\label{small_grad}
\max_{\bfX \in \Omega} \| \nabla \bfu \|  \ll \mathcal{O}(\delta), \quad \delta \leq 1.
\end{equation}
Henceforth, we shall not distinguish the dependence of the quantities on $\bfX$ for the notational convenience.  The premise of infinitesimal strains \eqref{small_grad} implies
\begin{subequations}\label{lin_results}
\begin{align}
\bfB &\approx \bfI + 2 \bfeps, \;\;  \bfC \approx \bfI + 2 \bfeps , \;\;  \bfE \approx  \bfeps,  \\
\det \bfF &= 1 + \tr(\bfeps). \label{detF_linear}
\end{align}
\end{subequations}
Let $\bfT \colon \Omega \to  \Sym(\mathbb{R}^{2 \times 2})$ be the \textit{Cauchy stress tensor} in the current configuration and it satisfies the linear momentum balance:
\begin{equation} 
\rho \, {\ddot{\bfu}} = \ddiv \, \bfT + \rho \, \bfb,
\end{equation}
where $\rho$ is the density in the current configuration and $\bfb \colon \Omega \to \mathbb{R}^2$ is the body force in the current configuration and the notation $\dot{\left(\cdot\right)}$ denotes the time derivative. The first and second \textit{Piola-Kirchhoff stress tensor tensors},  $\bfS \colon \Omega \to \mathbb{R}^{2 \times 2}$ and $\overline{\bfS} \colon \Omega \to \mathbb{R}^{2 \times 2}$, are defined by 
\begin{equation}
\bfS =  \bfT \bfF^{-\mathrm{T}} \, \det \bfF, \quad \overline{\bfS} :=\bfF^{-1} \bfS.
\end{equation}
The principle of angular momentum balance implies that the Cauchy stress tensor is symmetric, i.e.,
\begin{equation}
\bfT = \bfT^{\mathrm{T}}.
\end{equation}
The balance of mass (or continuity equation) in the material description follows that 
\begin{equation}\label{mass_balance1}
\rho_0 = \rho \, \det \bfF,
\end{equation}
where $\rho_0$ is the reference 
density. In the view of \eqref{detF_linear}, the balance of mass \eqref{mass_balance1} reduces to 
\begin{equation}\label{mass_balance2}
\rho_0 = \rho \, \left( 1 + \tr(\bfeps)  \right).
\end{equation}
More details about the \textit{kinematics} and \textit{kinetics} can be found in \cite{truesdell2004non}. 

\subsection{Implicit constitutive relations}
In this study, our  focus is to study the behavior of \textit{elastic (nondissipative) porous solids} whose material moduli depend on the density.  The response of such materials to the mechanical loading can be best described by \textit{implicit constitutive relations} introduced by Rajagopal  (see \cite{rajagopal2021a,rajagopal2021b} and the references therein).  A generalization for the elastic body defined by Rajagopal \cite{rajagopal2003,rajagopal2007elasticity} through an implicit type constitutive relation between the Cauchy stress and the deformation gradient is of the form:
\begin{equation}
\bff(\rho, \, \bfT,\,\bfF, \, \bfX) =0,
\end{equation}
where $\bff$ is a tensor-valued function. In the case of \textit{isotropic bodies}, the above implicit response relation reduces to 
\begin{equation}\label{eq:implicit_relation}
\bff(\rho, \, \bfT,\,\bfB) =0.
\end{equation}
{Assuming that the function $\bff$ is \textit{isotropic}}, then one can write the most general form of the implicit constitutive relation:
\begin{align} 
0 &= \delta_0 \, \bfI + \delta_1 \, \bfT + \delta_2 \, \bfB +  \delta_3 \, \bfT^2 + \delta_4 \, \bfB^2 +  \delta_5 \, ( \bfT \bfB + \bfB \bfT) + {\delta_6 \, ( \bfT^2 \bfB + \bfB \bfT^2)}  \notag \\
& + \delta_7 \, ( \bfB^2 \bfT + \bfT \bfB^2) +  \delta_8 \, ( \bfT^2 \bfB^2 + \bfB^2 \bfT^2),  
\end{align}
{where the material moduli $\delta_i$ for $i=0,\, 1, \,  \ldots, 8$} are scalar functions that depend on the density and the invariants for the pair $\bfT$ and $\bfB$, i.e.,
\begin{align}
\big\{ &\rho, \, \tr (\bfT), \, \tr(\bfB), \, \tr(\bfT^2), \, \tr(\bfB^2), \,\tr(\bfT^3), \, \tr(\bfB^3), \,\tr(\bfT \, \bfB), \, \tr(\bfT^2 \, \bfB) \notag \\
&\tr(\bfT \, \bfB^2), \, \tr(\bfT^2 \, \bfB^2)   \big\}.
\end{align}
Using the linearization assumption given in \eqref{small_grad} and the subsequent asymptotic results for the classical terms \eqref{lin_results}, we obtain a special subclass of the implicit constitutive relations of the form:
\begin{equation}\label{imp_model}
\widehat{\delta}_0 \, \bfI + \widehat{\delta}_1 \, \bfeps + \widehat{\delta}_2 \, \bfT + \widehat{\delta}_3 \, \bfT^2 + \widehat{\delta}_4 \, \left( \bfeps \, \bfT + \bfT \, \bfeps \right) + \widehat{\delta}_5 \, \left( \bfeps \, \bfT^2 + \bfT^2 \,\bfeps \right)    =0,
\end{equation}
where the terms $\widehat{\delta}_i, \; i=0, \ldots, 5$ are the functions of {scalar-valued} invariants of $\bfeps$ and $\bfT$. But the terms $\widehat{\delta}_i, \; i=0, \, 2, \, 3$ are the functions linearly depend upon the invariants of $\bfeps$ and arbitrarily upon the invariants of the $\bfT$, while the terms $\widehat{\delta}_i, \; i=1, \, 4, \, 5$ depends on the invariants of $\bfT$. An important subclass of of the above general class of implicit constitutive relations \eqref{imp_model} is
\begin{equation}\label{model_1}
\bfeps = \widehat{\alpha}_0 \, \bfI + \widehat{\alpha}_1 \, \bfT +  \widehat{\alpha}_2 \, \bfT^2,
\end{equation}
where the material moduli $\widehat{\alpha}_i, \; i=0, \, 1, \, 2$ depend on $\rho, \; \tr \bfT, \; \tr \bfT^2, \; \tr \bfT^3$, and by virtue of the  
mass balance (as in \eqref{mass_balance2}).  The model  \eqref{model_1} has been studied in several investigations involving  cracks in elastic bodies exhibiting strain-limiting behavior \cite{rajagopal2011modeling,gou2015modeling,HCY_SMM_MMS2022,itou2017nonlinear,MalliPhD2015}, quasi-static crack evolution  \cite{yoon2021quasi,lee2022finite},  thermo-elastic bodies \cite{bustamante2017implicit,yoon2022finite},  quasi-linear viscoelastic bodies \cite{itou2018states,zhu2016nonlinear,csengul2021nonlinear}, and nonlinear constitutive model for rock \cite{bustamante2020novel,bustamante2021bimodular}. 

Another subclass of models of the above general class of relations \eqref{imp_model} wherein the constitutive relation is linear in both $\bfeps$ and $\bfT$ is given by (see also \cite{itou2021implicit})
\begin{equation} \label{spe_model1}
\left( 1 + \kappa_3 \, \tr \, \bfT  \right) \; \bfeps = C_1 \, \left( 1 + \kappa_1 \, \tr \, \bfeps  \right) \; \bfT + C_2 \, \left( 1 + \kappa_2 \, \tr \, \bfeps  \right) \;  \left(  \tr \bfT \right)\; \bfI .
\end{equation}
In the above relation \eqref{spe_model1}, the moduli $\kappa_1, \, \kappa_2, \, \kappa_3, \, C_1, \, C_2$ are all constants. Note that the above model is both linear in stress and strain, and it can be used to describe the mechanical response of porous solids. The classical model can be recovered from the above model by choosing $\kappa_1, \, \kappa_2, \, \kappa_3$ 
as zero, and then we can identify 
\begin{equation}\label{eq:material_coeff}
C_1 = \dfrac{1 + \nu}{E} = \dfrac{1}{2 \, \mu} >0, \quad C_2 = - \dfrac{ \nu}{E}  <0,
\end{equation}
where $E$ is the Young's modulus and $\nu$ is the Poisson's ratio, and these are related to the Lam\`e constants $\lambda$ and $\mu$:
\begin{equation}\label{eq:linear_lame}
\lambda =  \dfrac{ E \, \nu}{(1 + \nu) (1 - 2 \nu)}, \quad \mu = \dfrac{ E }{2 (1 + \nu)}.
\end{equation}
Taking $\kappa_3 =0$ in \eqref{spe_model1}, we can  obtain a special constitutive relation
 \begin{equation} \label{spe_model2}
 \bfeps = C_1 \, \left( 1 + \kappa_1 \, \tr \, \bfeps  \right) \; \bfT + C_2 \, \left( 1 + \kappa_2 \, \tr \, \bfeps  \right) \;  \left(  \tr \bfT \right)\; \bfI .
\end{equation}
By virtue of the balance of mass in the reference configuration (as in \eqref{mass_balance2}), one can express the above relation as 
 \begin{equation} \label{spe_model3}
 \bfeps = \Xi_1^{C_1}(\rho, \; \tr \, \bfeps)     \; \bfT + \Xi_2^{C_2}(\rho, \; \tr \, \bfeps)  \;  \left(  \tr \bfT \right)\; \bfI ,
\end{equation}
with 
both functions $ \Xi_1^{C_1}(\rho, \; \tr \, \bfeps), \;  \Xi_2^{C_2}(\rho, \; \tr \, \bfeps) $ being linear in the density variable $\rho$,
hence linear in $\tr \, \bfeps$. 

\begin{remark}
The model \eqref{spe_model3} and the similar density-dependent models, introduced in \cite{itou2021implicit,murru2021stress,rajagopal2021b}, share the same feature  of crack-tip strain (as well as crack-tip  stress) singularities with the classical linearized model. However, the general strain-limiting models \cite{rajagopal2011modeling,gou2015modeling,kulvait2013,kulvait2019,ortiz2012,ortiz2014numerical} lead to the  
bounded crack-tip strain behavior. The uniform upper bound on the strain in the entire body can be fixed \textit{a priori} to a value, that is as small and bounded as modeled. 
Our model is certainly useful in describing the behavior of porous elastic solids, as the functions $ \; \Xi_1^{C_1}(\rho, \; \tr \, \bfeps) \;$ and $ \;  \Xi_2^{C_2}(\rho, \; \tr \, \bfeps) \;$ can be recognized as the density-dependent  material moduli.
\end{remark}

\begin{remark}
We can also introduce a new class of implicit constitutive model to describe the response elastic material via {splitting} the stress in terms of its deviatoric and 
volumetric parts.
Within the framework of models of type \eqref{spe_model1}, such stress splitting leads to a three-field (displacement, volumetric and deviatoric stresses) mixed formulation of the nonlinear elasticity model.  All the problems of static-crack field, nucleation, and damage evolution with interesting issues, such as the non-penetrating cracks in {bodies} \cite{itou2017nonlinear}, can then be addressed appropriately. 
\end{remark}

\begin{remark}
If an initial porosity 
is assumed to exist explicitly, which can be partially or fully saturated, the porosity can be one of the state variables based on the mixture theory \cite{CouO2004,SteH2019}. Poromechanics \cite{BioM1957,CouO2004,CouO2010,SteH2019} {is the fundamental theory for modeling mechanical response of porous solid filled with fluid or fluid mixture.}  Poroelasticity is thus usually referred as the elastic behavior of the deformable porous solid which has the fluid  saturated fully/partially inside the pores, where the deformation of the solid is then not only from the external loading but also from the internal pore pressure working on the walls (or skeleton) that partition the connected pores \cite{CouO2004,CouO2010}. 
Thermodynamically, 
the free energy ($\Psi_s$) of purely elastic skeleton (or porous structure) under the isothermal condition and small strain can be set with the conjugate relations of the state variables:
\begin{equation}\label{eq:skeleton_free_energy}
\bfT=\dfrac{\partial\Psi_s}{\partial\bfeps},\:\: p_J=\dfrac{\partial\Psi_s}{\partial\phi_J},
\end{equation}
where $p$ is the fluid pressure, $\phi$ is the true porosity in the deformed configuration, and the subscript $J$ denotes the fluid phase. As the strain and stress, the explicit porosity and fluid pressure are then in pair.
\end{remark}

\begin{remark}
As poromechanics divide the different continua, i.e., the solid and fluid phases, 
the Biot coefficient \cite{BioM1957} is introduced for their coupling. From the tangent analysis of the skeleton free energy, the Biot coefficient has the tensor form as $b_{ij}=\dfrac{\partial^2 \Psi_s}{\partial \epsilon\partial p}$. 
For the isotropic condition, it reduces to the scalar,  and for a linear isotropic poroelastic skeleton, it relates the variation of porosity to the variation of volumetric strain when pressure and temperature are held constant, defined as: 
\begin{equation}\label{eq:BiotC}
b = 1 - \dfrac{K_{dr}}{K_{s}},
\end{equation}
where $K_{dr}$ denotes the drained bulk modulus and $K_{s}$ is the modulus of matrix or solid grain.  Based on the Maxwell's symmetry \cite{CouO2004}, this parameter 
 then also relates the saturated fluid hydrostatic pressure change partitioned from the porous solid to the variation of total stress.
\end{remark}
\begin{remark} 
We note that the bulk modulus relates the volumetric change (e.g., dilation) linearly to the mean stress providing the pore pressure is zero \cite{CouO2004}.  
In petroleum or reservoir engineering,  particularly for the unconsolidated soil or rock, the fully coupled effects of the fluid flow and geomechanics then occur essentially both in the deformation of porous media and in the fluid flow, such as the Mandel-Cryer's effect \cite{YooH2018,YooH2021}.
In this study, we consider only the drained porous solid without fluid or pressure. 
An explicit model out of the form \eqref{eq:implicit_relation} with the state variable of fluid pressure ``$p$" 
can be a future study.
\end{remark}

\section{Boundary value problem and numerical method}\label{sec:BVP}
In this section, we develop a boundary value problem and the numerical approach for
stress and strain fields  
in the elastic porous solid body including the crack-tip. The proposed constitutive relation where the stress and linearized strain appear linearly, the material moduli of the porous material  are dependent upon the density.

\subsection{Mathematical model}
The elastic material body under consideration is homogeneous, isotropic, initially unstrained and unstressed. Let $\Omega$ be an open, bounded, Lipschitz, and connected domain with the boundary $\partial \Omega$ consisting of two smooth disjoint parts $\Gamma_N$ and $\Gamma_D$ such that $\partial \Omega = \overline{\; \Gamma_N} \,\cup \, \overline{\; \Gamma_D}$ where $\Gamma_{D} \cap   \Gamma_{N} = \emptyset$. To describe the state stress of the material under investigation, we consider the balance of linear momentum for the static body, which, in the absence of body forces, reduces to 
\begin{equation}
 - \nabla \cdot \bfT  = \boldsymbol{0} \quad \text{in} \quad \Omega.  \label{eq_blm}
\end{equation}
To model the stress response of the material  whose material moduli are dependent upon the density, e.g., the elastic porous solids, we make use of implicit constitutive relations (see Rajagopal \cite{rajagopal2021note} for details). 
A special subclass of such elastic bodies can be obtained  if one start with \eqref{spe_model2} and use $\kappa_1 = \kappa_2 = \beta$ as:
\begin{equation}\label{const_relation}
\bfeps = \dfrac{(1+\nu)(1+ \beta \, \tr(\bfeps))}{E} \, \bfT - \dfrac{\nu(1+ \beta \, \tr(\bfeps))}{E} \, \tr(\bfT) \,  \bfI. 
\end{equation}
We note that the above constitutive relationship is invertible, and the inverted relation is written for the Cauchy stress as:
\begin{equation}\label{def-T}
 \bfT=\frac{\mathbb{E}[\bfeps]}{1 + \beta \, \tr(\bfeps)},
\end{equation}
{which is definition for the stress of the proposed density-dependent model.} {Plugging \eqref{def-T} into \eqref{const_relation}, we note that the stress and strain relation reduces to the same form of the linearized elasticity:}
\begin{equation}\label{const_relation_linear}
\bfeps = \dfrac{(1+\nu)}{E} \, \mathbb{E}[\bfeps] - \dfrac{\nu}{E} \, \tr(\mathbb{E}[\bfeps]) \,  \bfI, 
\end{equation}
where the fourth-order tensor $\mathbb{E}[\cdot]$ are defined as:
\begin{equation}
\mathbb{E}[\bfeps] := \overline{c}_1 \, \bfeps + \overline{c}_2 \, \tr(\bfeps) \, \bfI, 
\end{equation}
and both the constants $\overline{c}_1$ and $\overline{c}_2$ depend on the {linearized} material parameters, i.e., 
the  conventional Young's modulus ($E$) and Poisson's ratio ($\nu$), defined as:
\begin{equation}\label{eq:c1_c2}
\overline{c}_1 := \dfrac{E}{1 + \nu}, \quad \overline{c}_2 := \dfrac{\nu \, E}{(1 + \nu)(1-2\nu)}.
\end{equation}
Note that $\overline{c}_1$ and $\overline{c}_2$ are identical to  linearized Lam\'{e} coefficients such that $\overline{c}_1=2\mu$ and $\overline{c}_2=\lambda$. Unlike \eqref{eq:material_coeff} and \eqref{eq:linear_lame}, material moduli are variants, thus starting from linearized $E$ and $\nu$, the nonlinear Lam\'{e} coefficients for the density-dependent model are:
\begin{equation}\label{eq:c1_c2_dd_model}
\lambda := \dfrac{\nu \, E}{(1 + \nu)(1-2\nu)(1+ \beta \, \tr(\bfeps))}, \quad \mu := \dfrac{E}{2(1 + \nu)(1+ \beta \, \tr(\bfeps))},
\end{equation}
 {and the generalized (drained) bulk modulus ($K_{dr}$) 
 is then derived from \eqref{eq:c1_c2} and \eqref{eq:c1_c2_dd_model} as:}
\begin{equation}\label{eq:bulk_modulus}
K_{dr}:=\dfrac{1}{1+ \beta \, \tr(\bfeps)}\left(\overline{c}_2+\dfrac{\overline{c}_1}{3}\right)=\dfrac{1}{1+ \beta \, \tr(\bfeps)}\left(\dfrac{\nu \, E}{(1 + \nu)(1-2\nu)}+\dfrac{E}{3(1 + \nu)}\right).
\end{equation}
Note that we can also calculate the nonlinear Young's modulus  and Poisson's ratio for this model from \eqref{eq:c1_c2_dd_model}. 
With $\tr(\bfeps)$, the density-dependent material moduli are set in \eqref{eq:c1_c2_dd_model} and \eqref{eq:bulk_modulus}. For the preferential stiffness, we see the intensity with its direction can also be designed using the nonlinear parameter for this model, $\beta$.

Combining the balance of linear momentum \eqref{eq_blm} and the constitutive relationship \eqref{const_relation}, we obtain the following governing system of equations: 
\begin{subequations}\label{system1}
\begin{align}
- \nabla \cdot \bfT &=  0 \quad \text{in} \quad \Omega, \label{eq1-updated} \\
\bfeps &= \dfrac{(1+\nu)(1+ \beta \, \tr(\bfeps))}{E} \, \bfT - \dfrac{\nu(1+ \beta \, \tr(\bfeps))}{E} \, \tr(\bfT) \,  \bfI  \quad \text{in} \quad \Omega. \label{eq2-updated}
\end{align}
\end{subequations}
The above system of partial differential equations (PDEs) need to be supplemented by appropriate boundary conditions with $\Gamma_N$ and $\Gamma_D$: 
\begin{subequations}\label{bcs-v1}
\begin{align}
\bfT  \bfn &= \bfg, \;\;  \mbox{for all} \;\; \bfx \in \Gamma_N, \\
  \bfu &= \bfu^0, \;\;  \mbox{for all} \;\; \bfx \in \Gamma_D, 
\end{align}
\end{subequations}
where $\bfg \colon \Omega \to \mathbb{R}^2$ is the given traction and $\bfu^0 \colon \Omega \to \mathbb{R}^2$ is the given boundary displacement data. Finally, we have the following formulation for the equilibrium problem for a two-dimensional elastic porous solid body:
\begin{formulation}
Given the material parameters, find $\bfu = (u_1, \, u_2)$, $\bfT = \left\{ \bfT_{ij} \right\}$, and $\bfeps = \left\{ \bfeps_{ij} \right\}, \; i, \, j = 1, \, 2$ such that 
\begin{subequations}\label{f1}
\begin{align}
- \,  \bfT_{ij, \; j}(\bfu) &= \bfzero, \; \mbox{in} \;\; \Omega, \;\; \mbox{and} \;\; i=1, \, 2, \\
 \bfeps_{ij} (\bfu) &= \dfrac{(1+\nu)(1+ \beta \, \tr(\bfeps))}{E} \, \bfT_{ij} - \dfrac{\nu(1+ \beta \, \tr(\bfeps))}{E} \, \tr(\bfT) \,  \delta_{ij}  \;\; \text{in} \;\; \Omega, \;\; \mbox{and} \;\; i=1, \, 2, \\
 \bfu_i &= \bfu^0_{i}, \quad  \mbox{on} \quad \Gamma_D, \;\; \mbox{and} \;\; i=1, \, 2, \\
 \bfT_{ij}(\bfu)  \bfn_{j} &= \bfg_{i}, \;\;  \mbox{on} \;\; \Gamma_N, \;\; \mbox{and} \;\; i=1, \, 2.
\end{align}
\end{subequations}
\end{formulation}
The above boundary value problem can also be formulated in variational form as a problem of minimizing the total strain-energy density functional. This approach will be studied in a future communication. 
As \eqref{f1} is nonlinear,
it is not tractable in the current form to any well-known analytical or numerical methods. In the following {section}, we propose a Newton's method for the linearization at the differential equation level, {followed by the} bilinear finite element method as a discretization technique. 

\subsection{Newton's iterative method and variational formulation}
To construct the numerical approximation to the boundary value problem (BVP) \eqref{f1}, we first linearize the differential equation to obtain sequence of linear problems. First, consider a mapping $\mathcal{L}(\cdot)  \colon \Sym(\mathbb{R}^{2 \times 2}) \mapsto \Sym(\mathbb{R}^{2 \times 2})$, with $\mathcal{L}(\bfzero) = \bfzero$, defined as: 
\begin{equation}\label{eq-L}
\mathcal{L}(\bfeps(\bfu)) :=  \frac{\mathbb{E}[\bfeps(\bfu)]}{ (1 + \beta \, \tr(\bfeps)) }. 
\end{equation}
The {assumption of smoothness} that we made for the Sobolev space \eqref{def-H01} and the quasi-linearity of the operator in \eqref{eq-L} bestow us  to apply the \textit{Newton's method} to obtain the linearized version of the PDE model \eqref{system1}. Given the initial guess $\bfu^0 \in \left(C^2(\Omega)\right)^2$ with $\mathbb{E}[\bfeps(\bfu^0)] \in \mathcal{C} \subseteq \Sym(\mathbb{R}^{2 \times 2})$, for each $n \in \mathbb{N}_0$, find $\bfu^{n+1} \in  \left(C^2(\Omega)\right)^2$ with $\mathbb{E}[\bfeps(\bfu^{n+1})] \in \mathcal{C}$ such that 
\begin{equation}
D \mathcal{L}(\bfeps(\bfu^n)) \left( \bfu^{n+1} - \bfu^n \right) = - \mathcal{L}(\bfeps(\bfu^n)),
\end{equation}
where $D\mathcal{L}(\bfeps(\bfu))$ denotes the \textit{Fre$\acute{e}$chet derivative}, 
defined as
\begin{equation}\label{FreDerivative}
{D \mathcal{L}(\bfeps(\bfu)):=} \lim_{\xi \to 0} \dfrac{\mathcal{L}(\bfeps(\bfu) + \xi \, \bfeps( \bfv))- \mathcal{L}(\bfeps(\bfu))}{\xi},
\end{equation}
for each $\bfv \in \left(C^2(\Omega)\right)^2$. Using \eqref{eq-L}  
and \eqref{FreDerivative}, we obtain
\begin{align}\label{frechet}
 D \mathcal{L}(\bfeps(\bfu^n)) \delta  \bfu^n  =  \Bigg[ &
\frac{ (\overline{c}_1/2) \left({ \nabla \delta\bfu^n + \left(\nabla  \delta\bfu^n \right)^T}\right)  + \overline{c}_2 \, (\nabla \cdot  \delta \bfu^n)\,  \bfI }{  (1 + \beta \, \tr(\bfeps(\bfu^n)))} \notag \\ 
&-
\frac{\beta  \;  \mathbb{E}[ \bfeps(\bfu^n)] \;  \left( \nabla \cdot \delta\bfu^n \right)}{\left( 1 + \beta \,  \tr( \bfeps(\bfu^n) \right)^{2}} \Bigg]\, ,
\end{align}
where the upper-script $``n"$ is the iteration number for the \textit{Newton's method}. Then combining the equations \eqref{frechet}, {\eqref{eq1-updated} and \eqref{eq2-updated}}, we obtain the updated {model equations}: 
\begin{subequations}
\begin{align}
&-\nabla \cdot \Bigg[ \frac{ (\overline{c}_1/2) \left({ \nabla \delta\bfu^n + \left(\nabla  \delta\bfu^n \right)^T}\right)  + \overline{c}_2 \, (\nabla \cdot  \delta \bfu^n)\,  \bfI }{  (1 + \beta \, \tr(\bfeps(\bfu^n)))} -
\frac{\beta  \;  \mathbb{E}[ \bfeps(\bfu^n)] \;  \left( \nabla \cdot \delta\bfu^n \right)}{\left( 1 + \beta \,  \tr( \bfeps(\bfu^n) \right)^{2}} \Bigg]\,  \notag \\
& =  -  \; \frac{\mathbb{E}[\bfeps(\bfu^n)]}{ (1 + \beta \, \tr(\bfeps(\bfu^n))) }.  \label{strong-form-1} \\
 & \notag \\
\bfeps(\bfu) &= \dfrac{(1+\nu)(1+ \beta \, \tr(\bfeps(\bfu)))}{E} \, \bfT(\bfu) - \dfrac{\nu(1+ \beta \, \tr(\bfeps(\bfu)))}{E} \, \tr(\bfT(\bfu)) \,  \bfI.
\end{align}
\end{subequations}
Then, the solution at the next iteration level is constructed by using
\begin{equation}
\bfu^{n+1}  = \bfu^{n} + \alpha^n \; \delta \bfu^{n},
\end{equation}
where $\alpha^n \in (0, 1 ]$ is the ``damping parameter'' that controls the convergence rate. In the \textbf{Formulation 2},
we will address  the procedure of computing the damping parameter with the line search algorithm. After the linearization, we have the formulation of the BVP at the continuous level.
\begin{formulation}\label{formul2}
Given the  linearized material parameters, Young's modulus $E$ and poisson's ratio $\nu$, initial guess $\bfu^{0} \in \left( C^2(\Omega)\right)^2$, find $\delta \bfu^n$ for $n=0, \,1, \, \ldots$, $\bfT = \left\{ \bfT_{ij} \right\}$, and $\bfeps = \left\{ \bfeps_{ij} \right\}, \; i, \, j = 1, \, 2$ such that 
\begin{subequations}\label{final_BVP}
\begin{align}
 &-\nabla \cdot \Bigg[ \frac{ (\overline{c}_1/2) \left({ \nabla \delta\bfu^n + \left(\nabla  \delta\bfu^n \right)^T}\right)  + \overline{c}_2 \, (\nabla \cdot  \delta \bfu^n)\,  \bfI }{  (1 + \beta \, \tr(\bfeps(\bfu^n)))} -
\frac{\beta  \;  \mathbb{E}[ \bfeps(\bfu^n)] \;  \left( \nabla \cdot \delta\bfu^n \right)}{\left( 1 + \beta \,  \tr( \bfeps(\bfu^n) \right)^{2}} \Bigg]\,  \notag \\
& =  - \mathcal{L}(\bfeps(\bfu^n)),  \label{SF-1} \\
  \bfu^{n}_i &= \bfu^{0}_{i} \quad  \mbox{on} \quad \Gamma_D, \;\; \mbox{and} \;\; i=1, \, 2,  \label{bc-1}\\
 \bfT_{ij}(\bfu^n)  \bfn_{j} &= \bfg_{i} \;\;  \mbox{on} \;\; \Gamma_N, \;\; \mbox{and} \;\; i=1, \, 2,   \label{bc-2}\\
  \bfu^{n+1} &= \bfu^n + \alpha^n \; \delta \bfu^n,  \label{u_np1}\\
\bfeps(\bfu^{n+1}) &= \dfrac{(1+\nu)(1+ \beta \, \tr(\bfeps(\bfu^{n+1})))}{E} \, \bfT(\bfu^{n+1}) - \dfrac{\nu(1+ \beta \, \tr(\bfeps(\bfu^{n+1})))}{E} \, \tr(\bfT(\bfu^{n+1})) \,  \bfI. 
\end{align}
\end{subequations}
\end{formulation}
From equation~\eqref{u_np1}, note that it is clear that we need to impose the zero Dirichlet boundary conditions at each iteration level for the Newton's update $\delta \bfu^n$.

\subsection{Continuous weak formulation}
In this {section},  we provide a variational formulation for the linearized version of the nonlinear BVP derived in the previous section. Also, the function spaces used to define the variational formulation has already been defined in the beginning, 
and the same setting is applied.  
To pose a weak formulation, we multiply the equations in the strong formulation \eqref{SF-1} with the test function from $\widehat{V}_{\bfzero}$ as in \eqref{test_fun_space}, then via integrating by parts using Green's formula together with the boundary conditions given in \eqref{bc-1} and \eqref{bc-2}, we arrive at the following variational formulation. \\

\begin{formulation}\label{formulation_3}
Given $\bfu^{0} \in \widehat{V}_{\bfg}$, for $n=0, 1, 2, \cdots$, find $\bfu^{n+1}:=\bfu^n + \alpha^n \,  \delta \bfu^n \in {V}$, such that 
\begin{equation}\label{eq:weak_formulation}
    A(\bfu^n; \, \delta\bfu^n,\bfv) = L(\bfu^n; \, \bfv)\quad \forall\, \bfv \in \widehat{V}_{\bfzero},
\end{equation}
where the bilinear term $A(\bfu^n; \, \delta\bfu^n, \, \bfv)$ and the linear term $L(\bfu^n; \, \bfv)$ are given by
\begin{multline}\label{A-L-Def}
A(\bfu^n; \, \delta\bfu^n,\, \bfv) = \int_{\Omega}  \Bigg[ \Bigg[ \frac{ (\overline{c}_1/2) \left({ \nabla \delta\bfu^n + \left(\nabla  \delta\bfu^n \right)^T}\right)  + \overline{c}_2 \, (\nabla \cdot  \delta \bfu^n)\,  \bfI }{  (1 + \beta \, \tr(\bfeps(\bfu^n)))} \\
-
\frac{\beta  \;  \mathbb{E}[ \bfeps(\bfu^n)] \;  \left( \nabla \cdot \delta\bfu^n \right)}{\left( 1 + \beta \,  \tr( \bfeps(\bfu^n) \right)^{2}}  \Bigg] \colon \bfeps(\bfv) \Bigg] \, d\bfx \, ,
\end{multline}
\begin{align}\label{A-L-Def-2}
    L(\bfu^n; \,  \bfv) &= -  \int_{\Omega} \left[\Bigg[\frac{ (\overline{c}_1/2)  \left( { \nabla  \bfu^n + \left(\nabla  \bfu^n \right)^T }  \right)  + \overline{c}_2 \, (\nabla \cdot  \bfu^n)\,  \bfI  }{\left( 1 + \beta \, \tr(\bfeps(\bfu^n)) \right) }\Bigg] \colon  \bfeps(\bfv)   \right]  \, d\bfx  \notag \\
    & + \int_{\Gamma_{N}} \bfT \bfn \cdot \bfv \, dx.
\end{align}
\end{formulation}

\subsection{Finite element discretization}
{In this section,} we first recall some basic notions and structure of the classical finite element method (FEM) to discretize  the weak formulation \eqref{eq:weak_formulation}. The meshes 
used for all computations in the numerical examples presented in this paper are quadrilaterals.  Let $\left\{ \mathcal{T}_h \right\}_{h >0}$ be a conforming, shape-regular (in the sense of Ciarlet \cite{ciarlet2002finite}) family of triangulation of the domain $\Omega$; $\mathcal{T}_h$  is a finite family of sets such that $K \in \mathcal{T}_h $ which 
implies $K$ is an open simplex with the mesh size $h_{K} := \text{diam}(K)$ for each $K$. Furthermore, we denote the largest diameter of the triangulation by 
\[
h:=\max_{K \in \mathcal{T}_h } \;  h_K.
\]
For any  $K_1, K_2 \in \mathcal{T}_h$, we have that $\overline{K}_1 \cap \overline{K}_2$ is either a null set or a vertex or an edge or  {the whole of $\overline{K}_1$ and $\overline{K}_2$},  and $\bigcup\limits_{K \in \mathcal{T}_h} \overline{K} = \overline{\Omega}$. We now define the classical piece-wise affine \textit{finite element space} to approximate the displacement variable ($\bfu$),
\begin{equation}\label{FEM_Q}
\widehat{V}_h = \left\{  \bfu_h \in  \left( C(\overline{\Omega})\right)^2 \colon \left. \bfu_h\right|_K \in \mathbb{Q}_k^2 \;~ \forall K \in \mathcal{T}_h \right\},
\end{equation} 
where $\mathbb{Q}_k$ is a set containing the tensor-product of {polynomials in $2$ variables up to order $k$} on the reference cell $\widehat{K}$. Then, the discrete approximation space is:
\begin{equation}\label{app-spaces}
V_h := \widehat{V}_h \, \cap  \, H^1(\Omega).
\end{equation}
The discrete counterpart of the 
continuous formulation (\textbf{Formulation 3})
is then as follows:
\begin{formulation}
Given $\bfu^{0}_h \in V_h$,  and the $n^{th}$ Newton's iterative solution, i.e., $\bfu^n_h \in {V}_h$, for $n=0, 1, 2, \cdots $, find $\bfu^{n+1}_h  \in V_h$, such that 
\begin{equation}\label{discrete_wf}
   A(\bfu^n; \; \delta \bfu^{n}_h,\, \bfv_h) = L(\bfu^n; \; \bfv_h), \forall\, \bfv_h \in \widehat{V}_h,  
\end{equation}
where the linear and bilinear term are given by:
\begin{multline}\label{disc_A}
A(\bfu_h^n; \, \delta\bfu_h^n,\, \bfv_h) = \int_{\Omega}  \Bigg[ \Bigg[ \frac{ (\overline{c}_1/2) \left({ \nabla \delta\bfu_h^n + \left(\nabla  \delta\bfu_h^n \right)^T}\right)  + \overline{c}_2 \, (\nabla \cdot  \delta \bfu_h^n)\,  \bfI }{  (1 + \beta \, \tr(\bfeps(\bfu_h^n)))} \\
-
\frac{\beta  \;  \mathbb{E}[ \bfeps(\bfu_h^n)] \;  \left( \nabla \cdot \delta\bfu_h^n \right)}{\left( 1 + \beta \,  \tr( \bfeps(\bfu_h^n) \right)^{2}}  \Bigg] \colon \bfeps(\bfv_h) \Bigg] \, d\bfx \, ,
\end{multline}
\begin{align}\label{disc_L}
    L(\bfu_h^n; \,  \bfv_h) &= -  \int_{\Omega} \left[\Bigg[\frac{ (\overline{c}_1/2)  \left( { \nabla  \bfu_h^n + \left(\nabla  \bfu_h^n \right)^T }  \right)  + \overline{c}_2 \, (\nabla \cdot  \bfu_h^n)\,  \bfI  }{\left( 1 + \beta \, \tr(\bfeps(\bfu_h^n)) \right) }\Bigg] \colon  \bfeps(\bfv_h)   \right]  \, d\bfx  \notag \\
     &+ \int_{\Gamma_N} \left( T_{12} v_1 + T_{22} v_2 \right) \, dx.
\end{align}
where
\begin{align}
\mathbb{E}[\bfeps(\bfu_h^{n})] &= \overline{c}_1  \, \bfeps(\bfu_h^{n}) + \overline{c}_2  \,  \tr(\bfeps(\bfu_h^{n})) \, \bfI  \notag \\
&= (\overline{c}_1/2)  \, \left( \nabla \bfu_h^{n} + ( \nabla \bfu_h^{n} )^{\mathrm{T}} \right) + \overline{c}_2  \,  \nabla \cdot \bfu_h^{n},
\end{align}
and 
\begin{equation}
\tr(\bfeps(\bfu_h^{n})) := \nabla \cdot \bfu_h^{n}.
\end{equation}
The solution at the next iteration level is given by $\bfu_h^{n+1} = \bfu_h^{n} + \alpha^n \, \delta \bfu_h^{n}$. 
\end{formulation}
{Initially,} we obtain the solution 
for the Newton's iteration 
via solving the linear problem (i.e., with $\beta =0$ in \eqref{disc_A}), {which is an appropriate guess for the solution of nonlinear problem.}

{Finally, we briefly outline the ``line search method'' used in the overall implementation of the FEM. Before we present the algorithm, define the residual of the overall differential equation as:}
\begin{align}\label{eq:residual}
f(\bfu^n, \; \bfv) &= \left( \mathcal{L} (\bfeps(\bfu^n)), \; \bfv \right)  \notag \\
&= \int_{\Omega} \frac{\mathbb{E}[\bfeps(\bfu^n)]   }{ (1 + \beta \, \tr(\bfeps (\bfu^n))) }  \colon \bfv  \; d\bfx.
\end{align}
{The following algorithm (\textbf{Algorithm~\ref{algo:b}}) summarizes the overall steps in the line search method \cite{wright1999numerical}.} 
\begin{algorithm}[h!]
\KwIn{ $\overline{\alpha}>0$ ($\overline{\alpha}=1$ guarantees the quadratic convergence), $ c_1 \in (0, \, 1), \; r_1 \; \& \;  r_2$ satisfying $0 < r_1 < r_2 <1$}
\KwOut{$\alpha^n$}
set $\alpha = \overline{\alpha}$; \\
while $f(\bfu^n + \alpha \, \boldsymbol{p}^n, \; \bfv) > f(\bfu^n, \; \bfv) + c_1 \; \alpha \, \nabla f(\bfu^n, \; \bfv) \cdot  \boldsymbol{p}^n$ 
\linebreak
  replace $\alpha$ by a new value in $[r_1 \, \alpha, \; r_2 \, \alpha]$\; 
return $\alpha^n = \alpha$\;
\caption{Description of backtrack line search algorithm}
\label{algo:b}
\end{algorithm}
The new value for $\alpha$ is obtained by minimizing the overall residual given in \eqref{eq:residual}. {There are two possibilities for generating the iterates by appropriately choosing the search direction. One is to select relevant values of $r_1$ and $r_2$ so that the interval $[r_1 \, \alpha, \; r_2 \, \alpha]$ always contains the damping parameter $\alpha$; Another simple strategy is to iteratively replace $\overline{\alpha}$ with $\overline{\alpha}/2$, until the value of $f(\bfu^n, \; \bfv)$ is sufficiently small.}

{Finally, the overall algorithm for the boundary value problem}  with the whole nonlinear density-dependent model is presented in \textbf{Algorithm~\ref{algo001}}.
\begin{algorithm}[h!]
\SetAlgoLined
\KwInput{Choose the parameters: $\beta, \, E, \, \nu$}
Start with a sufficiently refined mesh; \\
 \While{$[\text{Iteration Number} < \text{Max. Number of Iterations}]$.AND.$[\text{Residual} > \text{Tol.}]$}{
  Assemble the Equations~\eqref{disc_A} and \eqref{disc_L} for the displacements as primary variables using $\mathbb{Q}_1$ shape functions\;
 Use a \textit{direct solver} to solve for $\delta \bfu^{n}$\;
 Construct the solution at $n$-th iteration by a direct solver and then update the solution variable {using} $\bfu_h^{n+1} = \bfu_h^{n} +{\overline{\alpha}}\delta \bfu_h^{n}$\;
  \While{$[\text{Line Search Iteration Number} < \text{Max. Number of Line Search Iterations}]$}{
  Calculate \textit{Residual} using  Equation~\eqref{eq:residual}\;
  \If{$\text{Residual}\leq\text{Tol.}$}{
   Break\;
   }
   Do Algorithm~\ref{algo:b} to find the optimal value of $\alpha^n$\;
  }
  }  
 Write the final converged solution $\bfu_h^{n+1}$ to output files for post-processing\;
 Compute the crack-tip fields (e.g., stress and strain) for visualization\;
 \caption{Algorithm for the nonlinear density-dependent model}
 \label{algo001}
\end{algorithm}

\section{{Numerical experiments and discussion}}\label{sec:NumExp}
Primary purpose of this section lies in verifying the constitutive relation of density-dependent material moduli for an elastic porous solid via modeling and computational approaches provided in the previous sections. Using the proposed and classical models, we present several numerical tests  to compare the stress and strain distributions in the computational domains under different mechanical loadings. 
Those distributions with the preferential stiffness from density-dependent model can be confirmed particularly near the crack-tip. 
Ultimately, our goal is to suggest the rationale for 
this nonlinear model, 
from which 
physically meaningful and reliable responses of the damaged pores 
or the crack-tip can be described under the mechanical loading. 

Starting from Rajagopal's \textit{implicit constitutive theory} \cite{rajagopal2021b}, the proposed model encompasses a quasi-linear PDE system. Since there are no closed-form solutions available for such a system, a stable numerical method proposed in this study is utilized, 
where the approach 
focuses on linearizing at the differential equation level and then discretize the resultant elliptic BVP using continuous bilinear finite elements. 
All numerical implementations are done using the \texttt{deal.II} \cite{dealII90,dealiiweb}, finite element library.
Against the nonlinearity,
the Newton’s method is employed with 
the lower bound of the convergence (i.e., the tolerance) 
as $10^{-8}$ 
and the maximum number of Newton iteration step is set as $50$. {For the line search algorithm, a constant damping coefficient for the line search  is taken as $0.5$, i.e., iteratively replacing $\overline{\alpha}$ with $\overline{\alpha}/2$, and the maximum number of steps is taken as $10$ (see \textbf{Algorithm~\ref{algo:b} and \ref{algo001}}).
A direct solver is used as a linear solver to compute the numerical solution of the linearized system of equations.  

In {Section}~\ref{sec:h-conv}, we perform the $h$-convergence test  first {to confirm all the algorithms and code}, and we proceed to Section~\ref{sec:Example} where we test several BVPs 
with/without a crack. The intact domain is addressed first with mode-I and mode-II in Section~\ref{sec:ex1} (\textbf{\textit{Example} 1}); 
the crack problems are then presented 
in Section~\ref{sec:ex2}  (\textbf{\textit{Example} 2}),  Section~\ref{sec:ex3} (\textbf{\textit{Example} 3}), and Section~\ref{sec:ex4} (\textbf{\textit{Example} 4}) with mode-I,  mode-II, and finally mixed-mode (mode I and II) loadings, respectively.

\subsection{$h$-convergence study}\label{sec:h-conv}
\begin{figure}[h]
\centering
\begin{tikzpicture}
\draw (0,0) -- (3,0) -- (3,3) -- (0,3) -- (0,0);

\node at (-0.3,-0.25)   {$(0,0)$};
\node at (3.3,3.2)   {$(1,1)$};

\node at (-0.3, 1.5)   {$\Gamma_{D_{1}}$};
\node at (3.4, 1.5)   {$\Gamma_{D_{2}}$};
\node at (1.5, -0.25)   {$\Gamma_{D_{3}}$};
\node at (1.5, 3.2)   {$\Gamma_{D_{4}}$};
\end{tikzpicture}
\caption{An unit square domain and the Dirichlet boundary conditions for h-convergence study.}
\label{fig:h-conv}
\end{figure}
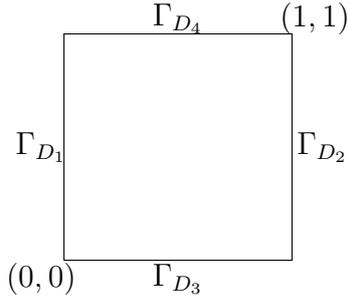
In this section, we introduce a sample problem to verify the mathematical description and algorithm of the numerical model proposed in this study. 
To this end, $h$-convergence test is performed with a manufactured solution set as  $$\bfu:=\left(\sin(\frac{\pi}{2}x), -\cos(\frac{\pi}{2}y)\right),$$ 
for an unit square domain {(1 m $\times$ 1 m)} as Figure~\ref{fig:h-conv}. The Dirichlet condition that satisfies the exact solution is applied for all the boundaries (from $\Gamma_{D_{1}}$ to $\Gamma_{D_{4}}$), resulting in $\Gamma_N=\emptyset$.  
For the nonlinear parameter, $\beta=1.0$ is taken, and for the linearized material moduli, the Young's modulus of {$E=100$ Pa and Poisson's ratio of $\nu=0.1$} are used. 
We then globally refine the whole domain with total 
6 cycles. 
In the sense of $L^{2}$-norm, 
{the optimal} convergence order of $2$  
for the bilinear polynomial is obtained. See the detailed rate values for convergence in Table~\ref{tab000}.

\begin{table}[!h]
\centering
\small
\begin{tabular}{|c|l||l|l||}
\hline
{Cycle of refinement} & \multicolumn{1}{c||}{{$h$}} &   \multicolumn{1}{c|}{$L^2$ Error} & \multicolumn{1}{c||}{Rate}   \\    \hline\hline                        
1                      &  0.5            & 0.053043115884           & -      
 \\ \hline
2                      & 0.25             & 0.013140372001           & 2.7317
 \\ \hline
3                      & 0.125             & 0.003275513882             & 2.3635
 \\ \hline
4                      & 0.0625           & 0.000818290922              & 2.1809
 \\ \hline
5                      & 0.03125            & 0.000204592608               & 2.0899
 \\ \hline
6                      & 0.015625           & 0.000051207706             & 2.0433
 \\ \hline
\end{tabular}
\caption{The results of $L^2$ error demonstrating the rate of optimal convergence.}
\label{tab000}
\end{table}%

\subsection{Boundary value problems with mechanical loadings}\label{sec:Example}

\subsubsection{Premise and setup for numerical experiments}
{In this study, {the elastic regime is of our interest, thus we do not consider any dissipation of energy.} We also do not consider any fluid saturation in the solid, thus the porous solid of our interest is in the unsaturated condition.  
Neither is introduced any explicit porosity concept into the model nor any initial 
porosity value that 
is directly measurable through experiments, or calculated via formula. However, the implicit concept of porosity 
is still considered  
with the volumetric strain, $\tr(\bfeps)$, such as in \eqref{mass_balance2}  
that relates the density of material to the mechanical moduli (\eqref{eq:c1_c2_dd_model} or \eqref{eq:bulk_modulus}). 

From \textbf{Example 1} to \textbf{Example 4} (Section~\ref{sec:ex1} to Section~\ref{sec:ex4}, respectively), {{we  
assume homogeneous isotropic material under isothermal condition.} 
For the same unit square domain in these examples {(1 m $\times$ 1 m)}, total 7 global refinements are performed, resulting in the uniformly refined mesh size of {$h=0.0078125$. For its initial linearized elastic moduli, the Young's modulus ($E$) of $100$ MPa and Poisson's ratio ($\nu$) of $0.15$ are taken.}}  
For the {Neumann boundary condition}, the same traction value of {$f_u=0.01$ MPa} 
is applied to the top boundary 
({see} Figure~\ref{Fig:ex1_setup} and Figure~\ref{Fig:ex2_3_4_setup}) with different modes of loading. 
The rest detailed boundary conditions for each example are described in each section. {The sign convention  
for stress follows such that the tensile stress is positive.}
As the study focuses on the nonlinear effects and different mechanical responses from the parameter of $\beta$-values, 
we compare 4 different cases, i.e., $\beta=-200, -50, +50,$ and $+200$, 
with the case of $\beta=0$, i.e., the classical linearized elasticity. 
Comparisons are then 
highlighted with stress and strain distributions, and strain density energy for each case, 
{focusing on} the area near the crack-tip.  Furthermore, we also compare the stress intensity factor, {volumetric strain, and bulk modulus} between the cases.

\subsubsection{\textbf{Example 1}: No crack problems}\label{sec:ex1}
\begin{figure}[!h]
\centering
\begin{minipage}{0.45\textwidth}
\centering
\begin{tikzpicture}
\coordinate (A) at (0,0);
\coordinate (B) at (3,0);
\coordinate (C) at (3,3);
\coordinate (D) at (0,3);
  \draw (A) -- (B);
\draw (B) -- (C);  
\draw (C) -- (D);
 \draw (D) -- (A);
\node at (-0.4, 0.05)   {$(0,0)$};
\node at (3.3,3.2)   {$(1,1)$};
 \node [below] at (1.5,0.5) {$\Gamma_1$};
 \node [right] at (3,1.5) {$\Gamma_2$};
 \node [left] at (0,1.5) {$\Gamma_4$};
 \node [above] at (1.5,2.4) {$\Gamma_3$};
 \node [above] at (1.52,1.45) {$\Omega_{1a}$};
 \draw[->] (0.2,3.1) -- (0.2,3.5);
\draw[->] (1.,3.1) -- (1.,3.5);
\draw[->] (2.,3.1) -- (2.,3.5);
\draw[->] (2.8,3.1) -- (2.8,3.5);

 \draw[->] (0.2,0.2) -- (0.2,0.7);
\draw[->] (0.2,0.2) -- (0.7,0.2);
\node at (0.2, 0.9)   {$y$};
 \node at (0.9,0.2) {$x$};

\draw [loosely dotted, line width=0.7mm, red]  (0.0,1.5) -- (3.0,1.5);
 \node[mark size=3.00pt] at (0.05,-0.09) {\pgfuseplotmark{triangle*}};
      \node[mark size=3.00pt] at (1.0,-0.09) {\pgfuseplotmark{triangle*}};
     \node[mark size=3.00pt] at (2.0,-0.09) {\pgfuseplotmark{triangle*}};
      \node[mark size=3.00pt] at (2.95,-0.09) {\pgfuseplotmark{triangle*}};
       \node[mark size=3.00pt] at (0.05,-0.18) {\pgfuseplotmark{*}};
      \node[mark size=3.00pt] at (1.0,-0.18) {\pgfuseplotmark{*}};
     \node[mark size=3.00pt] at (2.0,-0.18) {\pgfuseplotmark{*}};
      \node[mark size=3.00pt] at (2.95,-0.18) {\pgfuseplotmark{*}};
 \node at (1.5,3.7) {$\bfT \bfn=f_u\times(0,1)^{\mathrm{T}}$};
 \end{tikzpicture}
 \caption*{(a) \textbf{Example 1a}. }
\end{minipage}
\hspace{-2em}
\begin{minipage}{0.45\textwidth}
\centering
\begin{tikzpicture}
\coordinate (A) at (0,0);
\coordinate (B) at (3,0);
\coordinate (C) at (3,3);
\coordinate (D) at (0,3);
  \draw (A) -- (B);
\draw (B) -- (C);  
\draw (C) -- (D);
 \draw (D) -- (A);
\node at (-0.4, 0.05)   {$(0,0)$};
\node at (3.3,3.2)   {$(1,1)$};
 \node [below] at (1.5,0.5) {$\Gamma_1$};
 \node [right] at (3,1.5) {$\Gamma_2$};
 \node [left] at (0,1.5) {$\Gamma_4$};
 \node [above] at (1.5,2.4) {$\Gamma_3$};
 \node [above] at (1.52,1.45) {$\Omega_{1b}$};
\draw [loosely dotted, line width=0.7mm, red]  (0.0,1.5) -- (3.0,1.5);
 \draw[->] (0.0,3.15) -- (0.5,3.15);
\draw[->] (1.,3.15) -- (1.5,3.15);
\draw[->] (2.,3.15) -- (2.5,3.15);
 \node at (1.5,3.7) {$\bfT \bfn=f_u\times(1,0)^{\mathrm{T}}$};
  \draw[->] (0.2,0.2) -- (0.2,0.7);
\draw[->] (0.2,0.2) -- (0.7,0.2);
\node at (0.2, 0.9)   {$y$};
 \node at (0.9,0.2) {$x$};

    \node[mark size=5pt] at (0.05,-0.17) {\pgfuseplotmark{triangle*}};
     \node[mark size=5pt] at (1.0,-0.17) {\pgfuseplotmark{triangle*}};
     \node[mark size=5pt] at (2.0,-0.17) {\pgfuseplotmark{triangle*}};
      \node[mark size=5pt] at (2.95,-0.17) {\pgfuseplotmark{triangle*}};
  \end{tikzpicture}
\caption*{(b) \textbf{Example 1b}.}

\end{minipage}
\caption{{Numerical domains {(unit square, 1 m $\times$ 1 m)} for \textbf{Example 1} with traction boundary conditions in different modes: (left) \textbf{Example 1a} in mode-I and (right) \textbf{Example 1b} in mode-II. The red-dotted line is the reference line for each example.}}
\label{Fig:ex1_setup}
\end{figure} 
{In this example}, no crack or slit exists but intact porous solid is considered. Note that two different traction loads are applied to the top boundary ($\Gamma_{3}$), i.e.,  the tensile loading in mode-I (\textbf{Example 1a}) and in-plane shear loading in mode-II (\textbf{Example 1b}) in the these domains of $\Omega_{1a}$ (Figure~\ref{Fig:ex1_setup} (a)) and $\Omega_{1b}$ (Figure~\ref{Fig:ex1_setup} (b)), respectively. 
The detailed boundary conditions of each problem for the field variables are as follows: \textbf{Example 1a} has 
\begin{subequations}\label{bcs-u}
 \begin{align}
 \bfT \bfn &= \bfzero \quad\text{on} \quad \;  \Gamma_{2}, \; \Gamma_4,  \\
 \bfT \bfn \cdot \bfv & = {T_{21} \, v_1 + T_{22} \, v_2 = f_u \, v_2} \quad \text{on} \quad \Gamma_{3}, \label{trac_bc} \\
 {u_2} &= 0 \quad \text{on} \quad {\Gamma_{1}}, \label{eq:commom_boundary1}
 \end{align}
 \end{subequations}
 \noindent while \textbf{Example 1b} has 
\begin{subequations}\label{bcs-u-1}
 \begin{align}
 \bfT \bfn &= \bfzero \quad \text{on} \quad \;  \Gamma_{2}, \; \Gamma_4,  \\
 \bfT \bfn \cdot \bfv &= {T_{21} \, v_1 + T_{22} \, v_2 = f_u \, v_1} \quad\text{on} \quad \Gamma_{3}, \label{trac_bc-1} \\
 {u_1} &= 0 ~\text{and}~{u_2 = 0} \quad \text{on} \quad {\Gamma_{1}}. \label{eq:commom_boundary1-1}
 \end{align}
 \end{subequations}
 Note that we have slightly different Dirichlet boundary conditions for the bottom boundary ($\Gamma_1$) for the two problems; the roller boundary conditions are set for \textbf{Example 1a}, thus there is no displacement to its normal direction, i.e., $y$-direction, but it is free in $x$-direction. Meanwhile, hinge boundary conditions are set for \textbf{Example 1b}, thus there are no displacments in both $x$- and $y$-directions. {For these problems, we compare the displacements on the reference line (red-dotted line, {$L=1$ m}) in the $x$- 
 and $y$-directions, i.e., parallel and perpendicular to the reference line, respectively.} 
 
\begin{figure}[!h]
\centering
\includegraphics[width=1.0\textwidth]{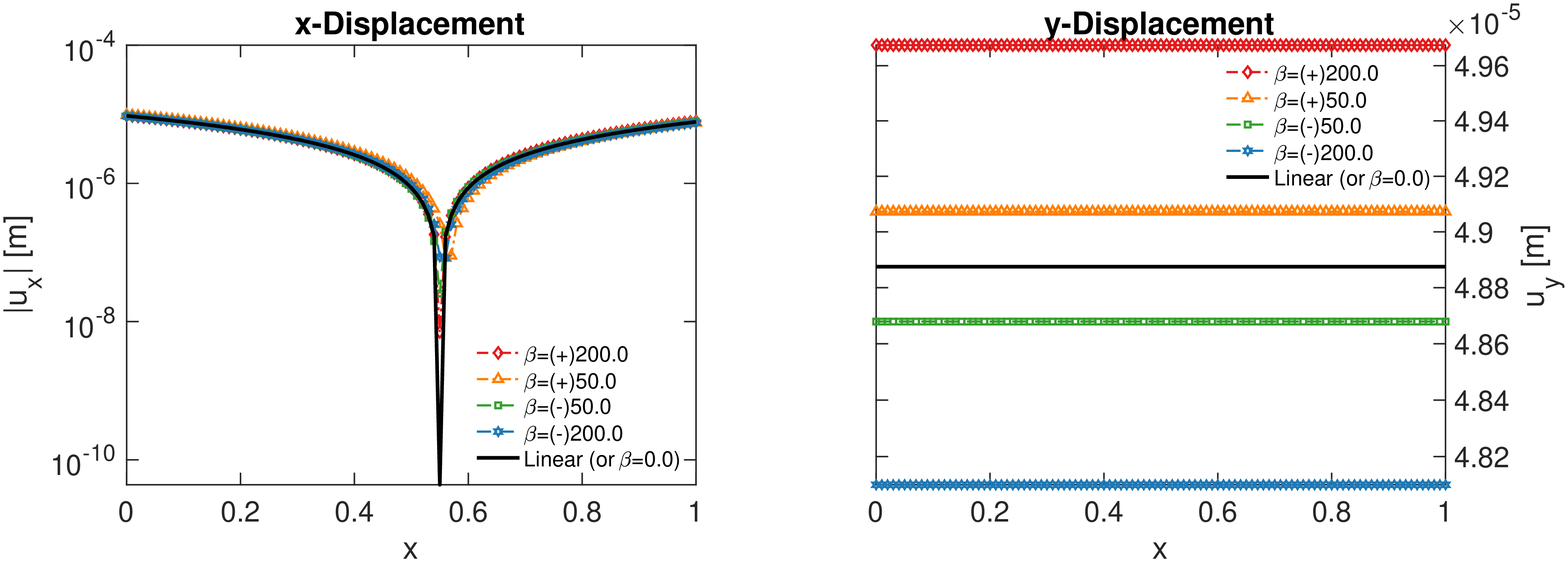}
\caption{Displacements {(unit: m)} on the reference line {($L=1$ m)} for \textbf{Example 1a}: (left) x-displacement in absolute value and (right) y-displacement. } 
\label{fig:disp_1a}
\end{figure}
\begin{figure}[!h]
\centering
\includegraphics[width=1.0\textwidth]{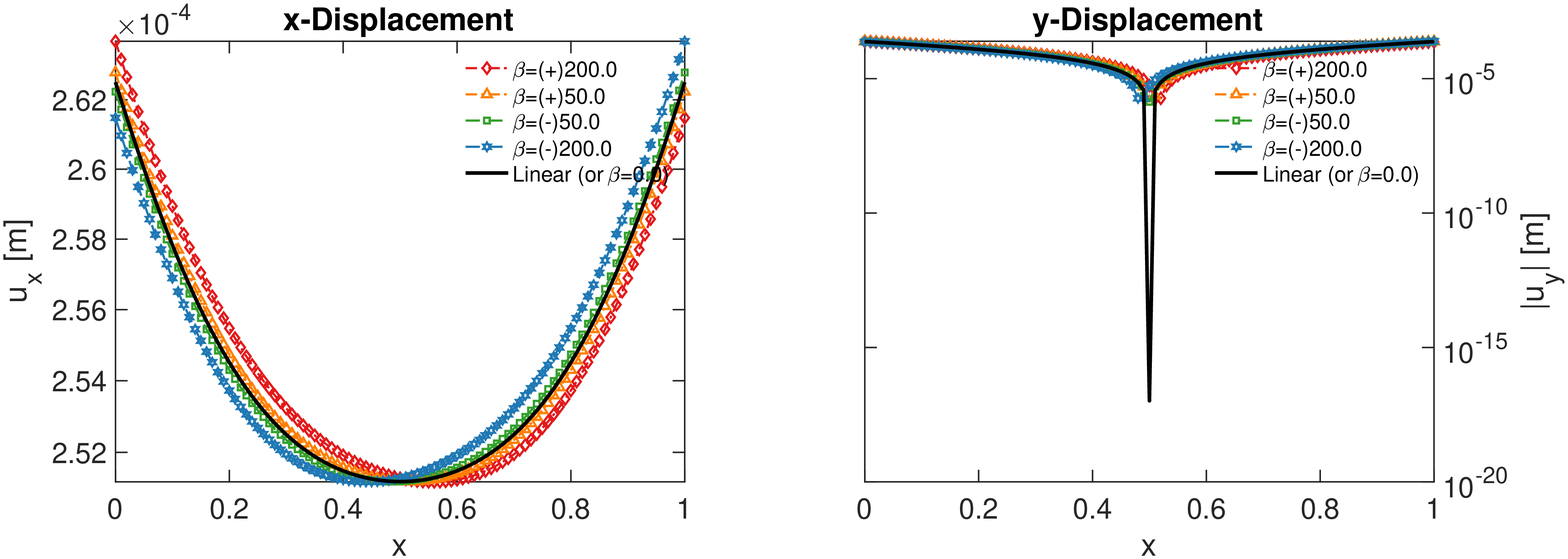}
\caption{Displacements {(unit: m)} on the reference line {($L=1$ m)} for \textbf{Example 1b}: (left) x-displacement and (right) y-displacement in absolute value.} 
\label{fig:disp_1b}
\end{figure}

{Figure}~\ref{fig:disp_1a} and Figure~\ref{fig:disp_1b} demonstrate the displacements for \textbf{Example 1a} and \textbf{Example 1b}, respectively. In each figure, 
the solutions of $x$-displacement on the left and of $y$-displacement on the right are exhibited, where we take absolute values with the semi-log scale in the direction perpendicular to the loading, i.e., $x$-direction for \textbf{Example 1a} and $y$-direction for \textbf{Example 1b}. 
In the perpendicular direction to the loading, we find that the nonlinear model with different {$\beta$-values} yields distinct responses near the center. {In the parallel direction to the loading for the mode-I and mode-II, 
we find that positive and negative $\beta$-values are reflected in the displacement in the opposite way. For example, positive $\beta$ in tension (Figure~\ref{fig:disp_1a}) yields larger 
$y$-displacement than that of the linear, while negative $\beta$ has smaller one. Thus, it implies that the strengths of a material with tensile and shear (in compressive direction) stresses are different under the nonlinear model. {
To investigate the preferential stiffness with these positive and negative $\beta$-values, 
we focus on the stress and strain distributions depending on the volumetric strain changes, particularly near an edge crack in the following examples; 
} 
{Henceforth}, the same numerical domain is addressed, having different boundary conditions depending on the mode of loading in each example but with the same right-edge crack. 
The crack is expressed with the slit ($\Gamma_{C}$, the blue lines in Figure~\ref{Fig:ex2_3_4_setup} (a), (b), (c)), where zero traction is applied. The traction value ({$f_u=0.01$ MPa}) and other mechanical properties are the same as \textbf{Example 1}. For these examples, the main comparisons 
lie in the stress and strain distributions including the stress intensity factor on the reference line (the {red-dotted} lines in Figure~\ref{Fig:ex2_3_4_setup}, {$L=0.5$ m}), where the volumetric strain and bulk modulus are also compared between the models. These post-processing works for the variables are based on the average values in each grid element of interest
using quadrature points inside.

\subsubsection{\textbf{Example 2}: {Tensile} loading with crack}\label{sec:ex2}
\begin{figure}[!h]
\centering
\begin{minipage}{0.35\textwidth}
\centering
\begin{tikzpicture}
\coordinate (A) at (0,0);
\coordinate (B) at (3,0);
\coordinate (C) at (3,3);
\coordinate (D) at (0,3);
  \draw (A) -- (B);
\draw (B) -- (C);  
\draw (C) -- (D);
 \draw (D) -- (A);
 \draw [line width=0.5mm, blue]  (1.5,1.5) -- (3,1.5);
\draw [loosely dotted, line width = 0.7mm, red] (0,1.5) -- (1.5,1.5);
\node at (-0.4, 0.05)   {$(0,0)$};
\node at (3.3,3.2)   {$(1,1)$};
\node at (2.2, 1.8)   {$\Gamma_{C}$};
 \node [below] at (1.5,0.5) {$\Gamma_1$};
 \node [right] at (3,1.5) {$\Gamma_2$};
 \node [left] at (0,1.5) {$\Gamma_4$};
 \node [above] at (1.5,2.4) {$\Gamma_3$};
 \draw[->] (0.2,3.1) -- (0.2,3.5);
\draw[->] (1.,3.1) -- (1.,3.5);
\draw[->] (2.,3.1) -- (2.,3.5);
\draw[->] (2.8,3.1) -- (2.8,3.5);
 \node at (1.5,3.7) {$\bfT \bfn=f_u\times(0,1)^{\mathrm{T}}$};
 
 \draw[->] (0.2,0.2) -- (0.2,0.7);
\draw[->] (0.2,0.2) -- (0.7,0.2);
\node at (0.2, 0.9)   {$y$};
 \node at (0.9,0.2) {$x$};
      \node[mark size=3.00pt] at (0.05,-0.09) {\pgfuseplotmark{triangle*}};
      \node[mark size=3.00pt] at (1.0,-0.09) {\pgfuseplotmark{triangle*}};
     \node[mark size=3.00pt] at (2.0,-0.09) {\pgfuseplotmark{triangle*}};
      \node[mark size=3.00pt] at (2.95,-0.09) {\pgfuseplotmark{triangle*}};
       \node[mark size=3.00pt] at (0.05,-0.18) {\pgfuseplotmark{*}};
      \node[mark size=3.00pt] at (1.0,-0.18) {\pgfuseplotmark{*}};
     \node[mark size=3.00pt] at (2.0,-0.18) {\pgfuseplotmark{*}};
      \node[mark size=3.00pt] at (2.95,-0.18) {\pgfuseplotmark{*}};
 \end{tikzpicture}
 \caption*{(a) \textbf{Example 2}.}
\end{minipage}
\hspace{-2.5em}
\begin{minipage}{0.35\textwidth}
\centering
\begin{tikzpicture}
\coordinate (A) at (0,0);
\coordinate (B) at (3,0);
\coordinate (C) at (3,3);
\coordinate (D) at (0,3);
  \draw (A) -- (B);
\draw (B) -- (C);  
\draw (C) -- (D);
 \draw (D) -- (A);
 \draw [line width=0.5mm, blue]  (1.5,1.5) -- (3,1.5);
\draw [loosely dotted, line width = 0.7mm, red] (0,1.5) -- (1.5,1.5);
\node at (-0.4, 0.05)   {$(0,0)$};
\node at (3.3,3.2)   {$(1,1)$};
\node at (2.2, 1.8)   {$\Gamma_{C}$};
 \node [below] at (1.5,0.5) {$\Gamma_1$};
 \node [right] at (3,1.5) {$\Gamma_2$};
 \node [left] at (0,1.5) {$\Gamma_4$};
 \node [above] at (1.5,2.4) {$\Gamma_3$};
 \draw[->] (0.0,3.15) -- (0.5,3.15);
\draw[->] (1.,3.15) -- (1.5,3.15);
\draw[->] (2.,3.15) -- (2.5,3.15);
 \node at (1.5,3.7) {$\bfT \bfn=f_u\times(1, 0)^{\mathrm{T}}$};
 
  \draw[->] (0.2,0.2) -- (0.2,0.7);
\draw[->] (0.2,0.2) -- (0.7,0.2);
\node at (0.2, 0.9)   {$y$};
 \node at (0.9,0.2) {$x$};
 
   \node[mark size=5pt] at (0.05,-0.17) {\pgfuseplotmark{triangle*}};
     \node[mark size=5pt] at (1.0,-0.17) {\pgfuseplotmark{triangle*}};
     \node[mark size=5pt] at (2.0,-0.17) {\pgfuseplotmark{triangle*}};
      \node[mark size=5pt] at (2.95,-0.17) {\pgfuseplotmark{triangle*}};
 \end{tikzpicture}
 \caption*{(b) \textbf{Example 3}.}
\end{minipage}
 \hspace{-2.5em}
\begin{minipage}{0.35\textwidth}
\centering
\begin{tikzpicture}
\coordinate (A) at (0,0);
\coordinate (B) at (3,0);
\coordinate (C) at (3,3);
\coordinate (D) at (0,3);
  \draw (A) -- (B);
\draw (B) -- (C);  
\draw (C) -- (D);
 \draw (D) -- (A);
 \draw [line width=0.5mm, blue]  (1.5,1.5) -- (3,1.5);
\draw [loosely dotted, line width = 0.7mm, red] (0,1.5) -- (1.5,1.5);
\node at (-0.4, 0.05)   {$(0,0)$};
\node at (3.3,3.2)   {$(1,1)$};
\node at (2.2, 1.8)   {$\Gamma_{C}$};
 \node [below] at (1.5,0.5) {$\Gamma_1$};
 \node [right] at (3,1.5) {$\Gamma_2$};
 \node [left] at (0,1.5) {$\Gamma_4$};
 \node [above] at (1.5,2.4) {$\Gamma_3$};
 \draw[->] (0.2,3.15) -- (0.7,3.15);
\draw[->] (1.2,3.15) -- (1.7,3.15);
\draw[->] (2.3,3.15) -- (2.8,3.15);
 \draw[->] (0.2,3.15) -- (0.2,3.5);
\draw[->] (1.2,3.15) -- (1.2,3.5);
\draw[->] (2.3,3.15) -- (2.3,3.5);
 \node at (1.5,3.7) {$\bfT \bfn=f_u\times(1,1)^{\mathrm{T}}$};
 
  \draw[->] (0.2,0.2) -- (0.2,0.7);
\draw[->] (0.2,0.2) -- (0.7,0.2);
\node at (0.2, 0.9)   {$y$};
 \node at (0.9,0.2) {$x$};
 
   \node[mark size=5pt] at (0.05,-0.17) {\pgfuseplotmark{triangle*}};
     \node[mark size=5pt] at (1.0,-0.17) {\pgfuseplotmark{triangle*}};
     \node[mark size=5pt] at (2.0,-0.17) {\pgfuseplotmark{triangle*}};
      \node[mark size=5pt] at (2.95,-0.17) {\pgfuseplotmark{triangle*}};
 \end{tikzpicture}
 \caption*{(c) \textbf{Example 4}.}
\end{minipage}
\caption{Numerical domains {(unit square, 1 m $\times$ 1 m)} with cracks for \textbf{Example 2} to \textbf{Example 4} with traction boundary conditions in different modes: (a) \textbf{Example 2} in mode-I, (b) \textbf{Example 3} in mode-II, and (c) \textbf{Example 4} in mixed-mode. {Crack ($\Gamma_{C}$) in blue in each domain has the same geometry with the length of 0.5 m and the red-dotted line is the reference line for each example.}}
\label{Fig:ex2_3_4_setup}
\end{figure} 
For \textbf{Example 2}, a tensile loading in pure mode-I is applied to a numerical domain with a right-edge crack as illustrated in Figure~\ref{Fig:ex2_3_4_setup} (a), which has the following boundary conditions: 
\begin{subequations}\label{bcs-u}
 \begin{align}
 \bfT \bfn &= \bfzero  \quad\text{on} \quad \;  \Gamma_{2}, \; \Gamma_4, \; \Gamma_C, \\
 \bfT \bfn \cdot \bfv & = {T_{21} \, v_1 + T_{22} \, v_2 = f_u \, v_2  \quad \text{on} \quad \Gamma_{3},} \label{trac_bc} \\
 u_2 &= 0 \quad \text{on} \quad\Gamma_{1}. \label{eq:commom_boundary1}
 \end{align}
 \end{subequations} 
The roller boundary conditions are applied to $\Gamma_{1}$, resulting in the homogeneous Dirichlet boundary in $y$-direction, i.e., $u_y = 0$, without a constraint in $x$-direction. 
\begin{figure}[!h]
\centering
\subfloat[$T_{22}~\text{with}~\beta=-200$]{\includegraphics[width= 0.3\textwidth]{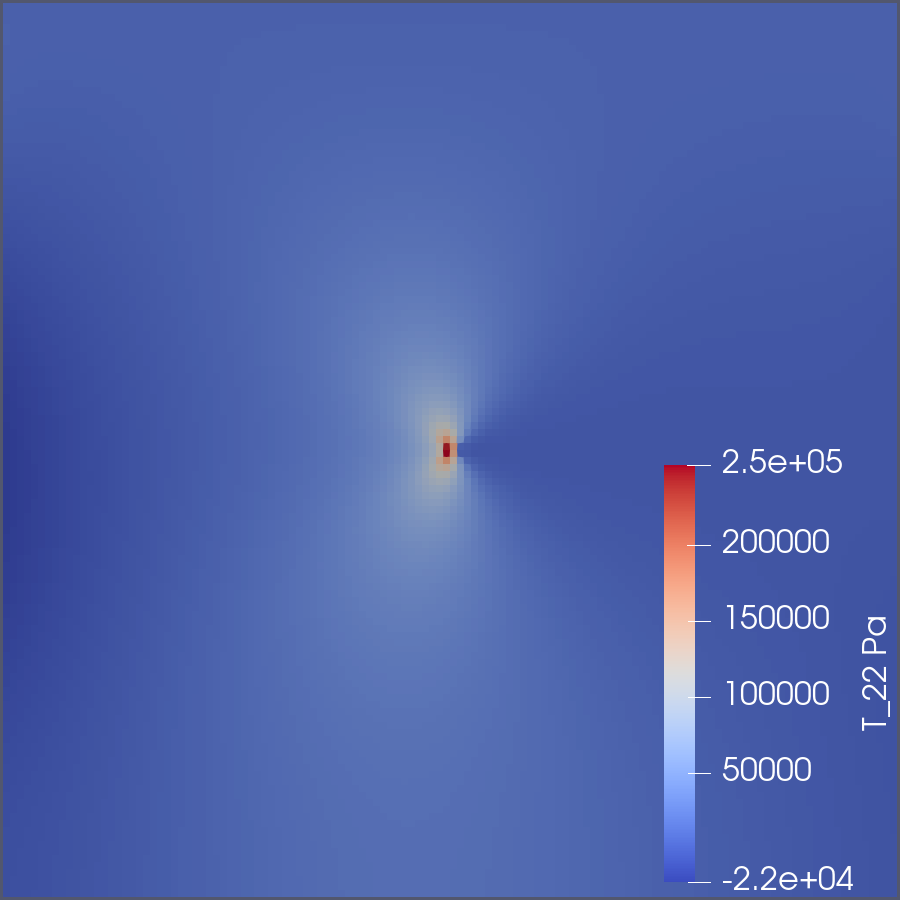}}
\hspace*{0.2in}
\subfloat[$T_{22}~\text{with}~\beta=0~({Linear})$]{\includegraphics[width= 0.3\textwidth]{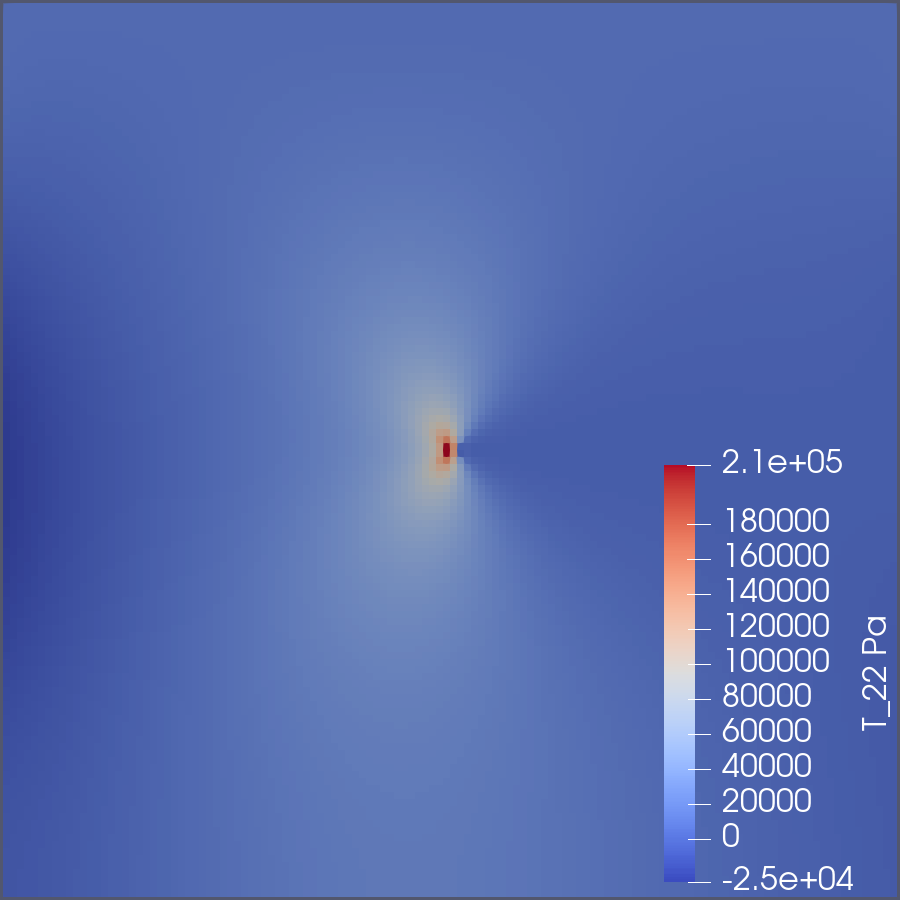}} \hspace*{0.2in}
\subfloat[$T_{22}~\text{with}~\beta=+200$]{\includegraphics[width= 0.3\textwidth]{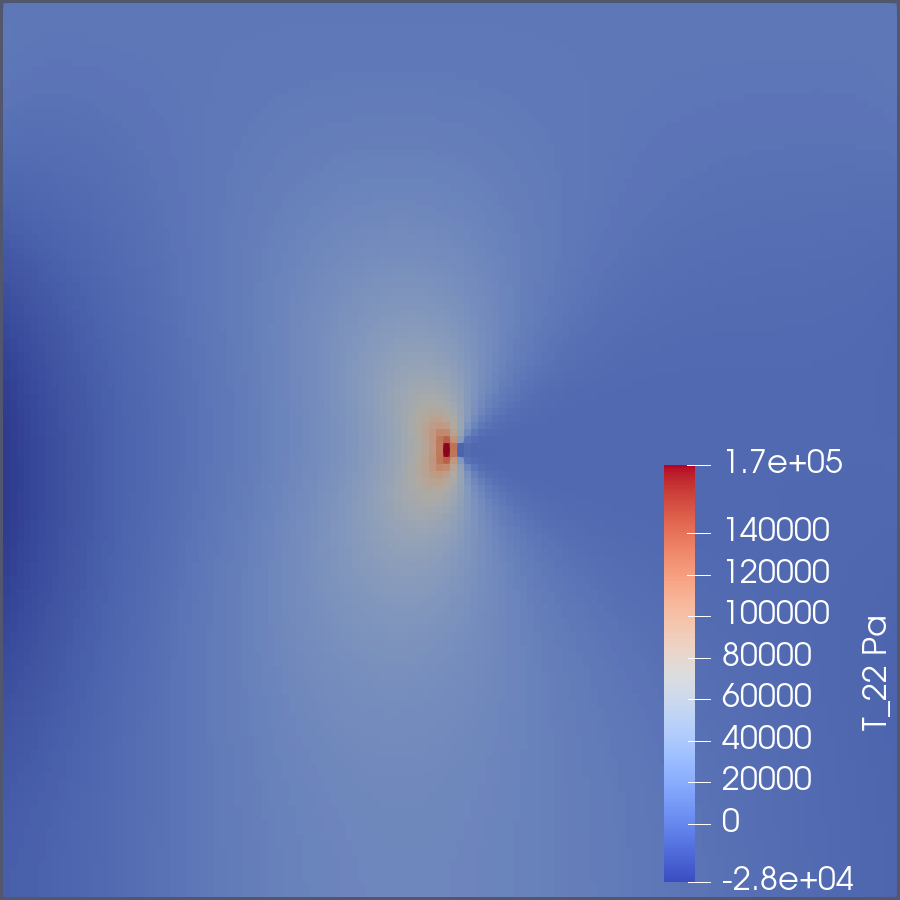}}
\\
\subfloat[$\epsilon_{22}~\text{with}~\beta=-200$]{\includegraphics[width= 0.3\textwidth]{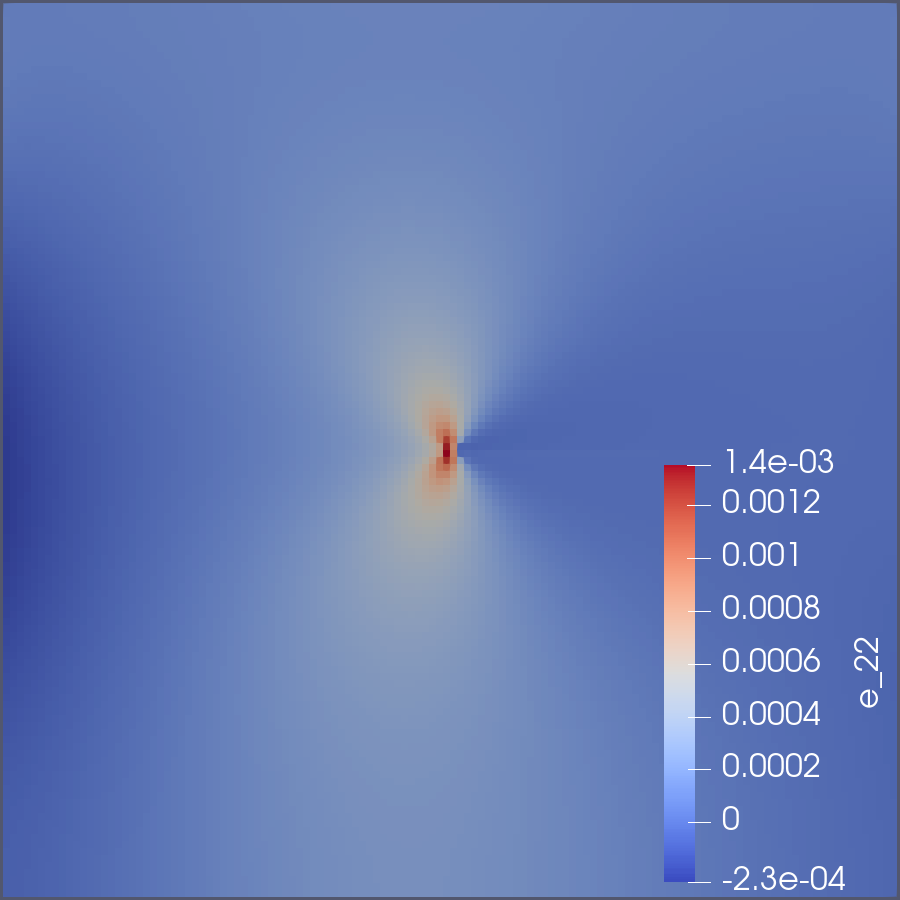}} \hspace*{0.2in}
\subfloat[$\epsilon_{22}~\text{with}~\beta=0~({Linear})$]{\includegraphics[width= 0.3\textwidth]{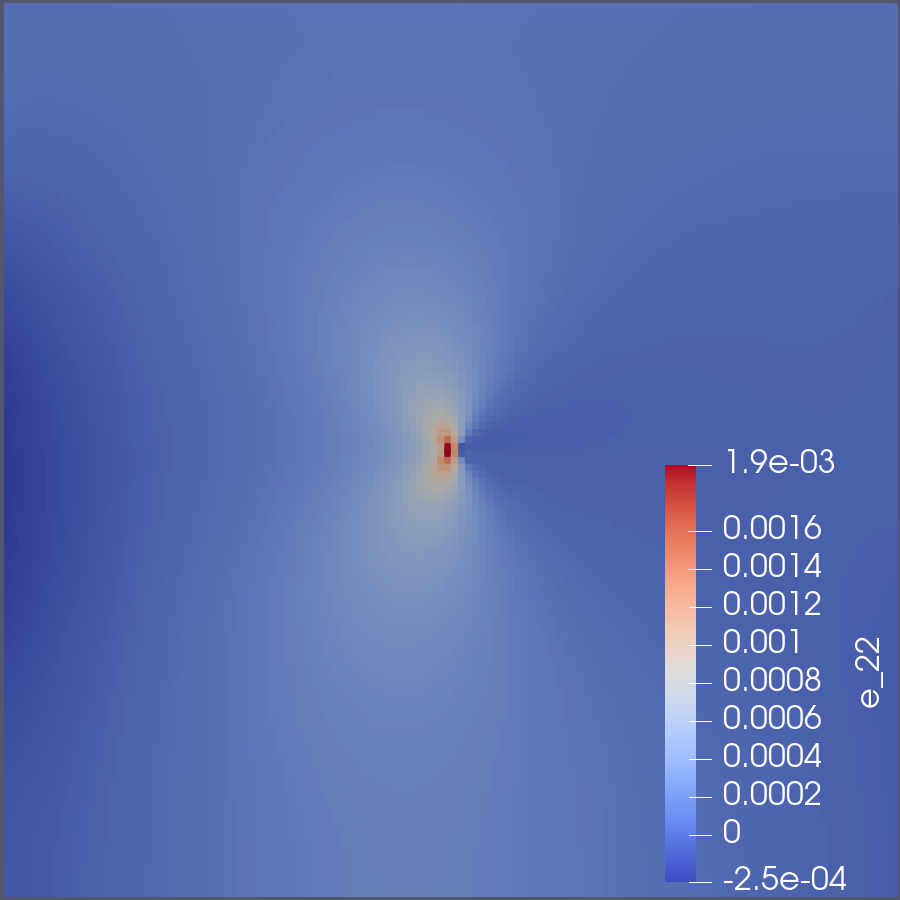}}
\hspace*{0.2in}
\subfloat[$\epsilon_{22}~\text{with}~\beta=+200$]{\includegraphics[width= 0.3\textwidth]{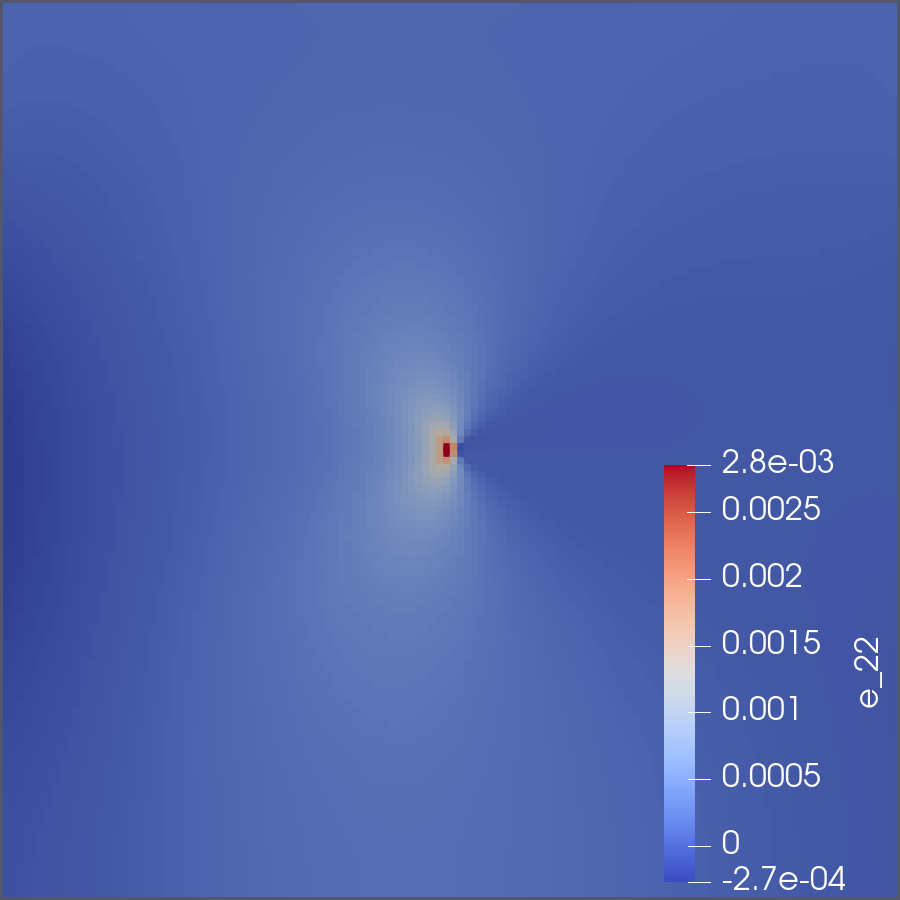}} 
\\
\subfloat[$SED~\text{with}~\beta=-200$]{\includegraphics[width= 0.3\textwidth]{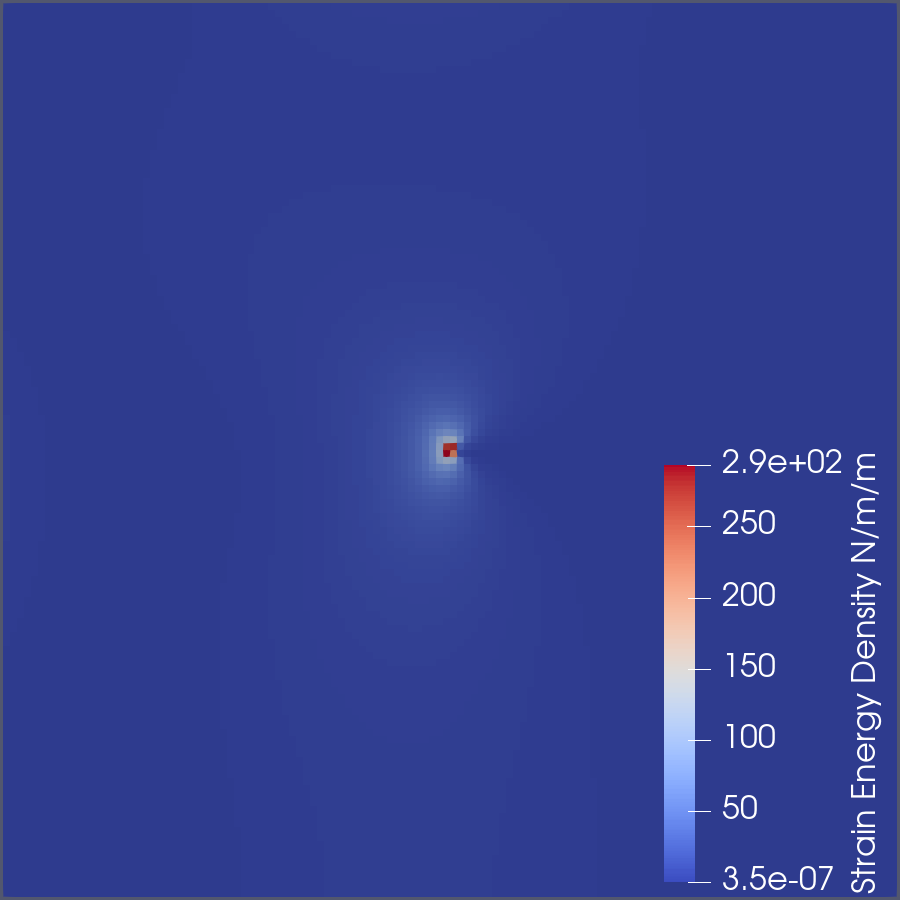}} \hspace*{0.2in}
\subfloat[$SED~\text{with}~\beta=0~({Linear})$]{\includegraphics[width= 0.3\textwidth]{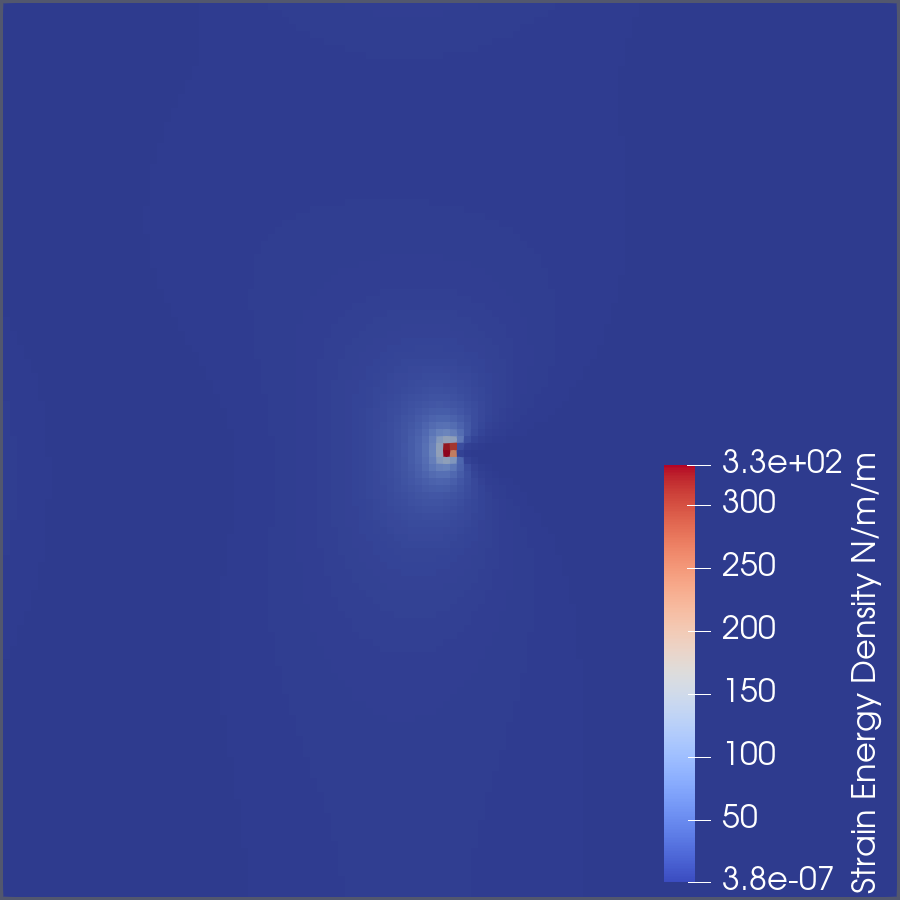}}
\hspace*{0.2in}
\subfloat[$SED~\text{with}~\beta=+200$]{\includegraphics[width= 0.3\textwidth]{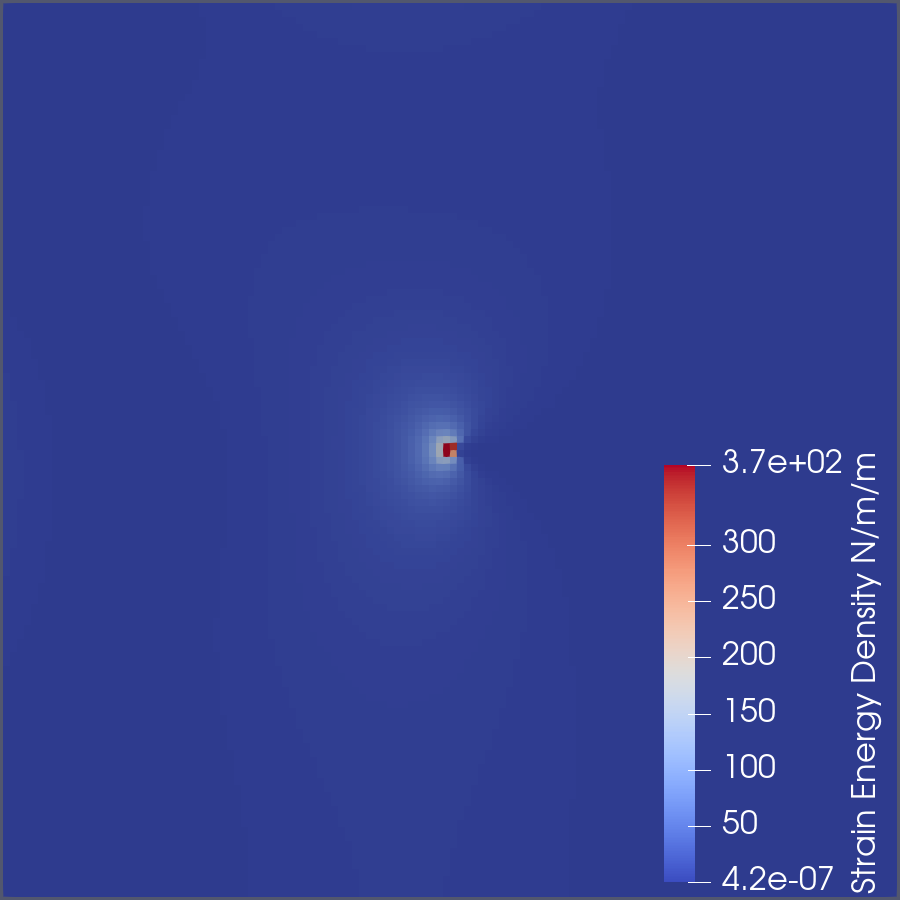}} 
\caption{{[\textbf{Example 2}] ((a), (b), (c)) stress $T_{22}$ {(unit: Pa)}, ((d), (e), (f)) strain $\epsilon_{22}$, and ((g), (h), (i)) strain energy density ($SED$, {unit: Pa}) distributions: (left) with $\beta=-200$, (middle) with $\beta=0$ (i.e., the linear model), and (right) with $\beta=+200$.}} 
\label{fig:ex2_stress_strain_density}
\end{figure}
\\
\newline
\noindent\textit{\textbf{stress and strain} }
{In Figure~\ref{fig:ex2_stress_strain_density} (a) to (c), $T_{22}$ {with $\beta=-200$, $\beta=0$, i.e., the linear model, and $\beta=+200$}, respectively, are illustrated, and the corresponding $\epsilon_{22}$ are shown in (d) to (f). For (g) to (i) in Figure~\ref{fig:ex2_stress_strain_density}, the 
strain energy density ($SED$)  
with the unit of $N/m/m$ or Pa is illustrated 
corresponding to each case. For simplicity, the {density-dependent} nonlinear model with $\beta=+50$ and $\beta=-50$ results are omitted for now.} 
{Note that $T_{22}$ and $\epsilon_{22}$ 
are parallel to the mode-I loading (normal to the horizontal line)}. 
It is found that  
larger value a variable has,  
narrower is its distribution with more localization. 
For the tensile stress ($T_{22}$) of which sign convention is positive, 
we identify that $\beta=+200$ yields the smallest with the largest of $\epsilon_{22}$, which implies the material strength becomes weaker against the tension. 
{The smallest tensile strain is obtained for $\beta=-200$. 
From the $SED$ results, we confirm that the maximum strain energy density for each case occurs right in front of the crack-tip. 
{We also identify that $SED$ for the case of $\beta=+200$ shows the largest with the maximum positive axial strain}; 
i.e., $\epsilon_{22}$ obtained from the positive $\beta$, e.g., $\beta=+200$, is larger than that from the linear and even larger than that from the negative $\beta$, i.e., $\beta=-200$, near the tip. {The detailed \textit{max} and \textit{min} values for each case are presented in Table~\ref{tab:ex2}.}  
\begin{figure}[htbp]
\centering
\includegraphics[width=1.00\textwidth]{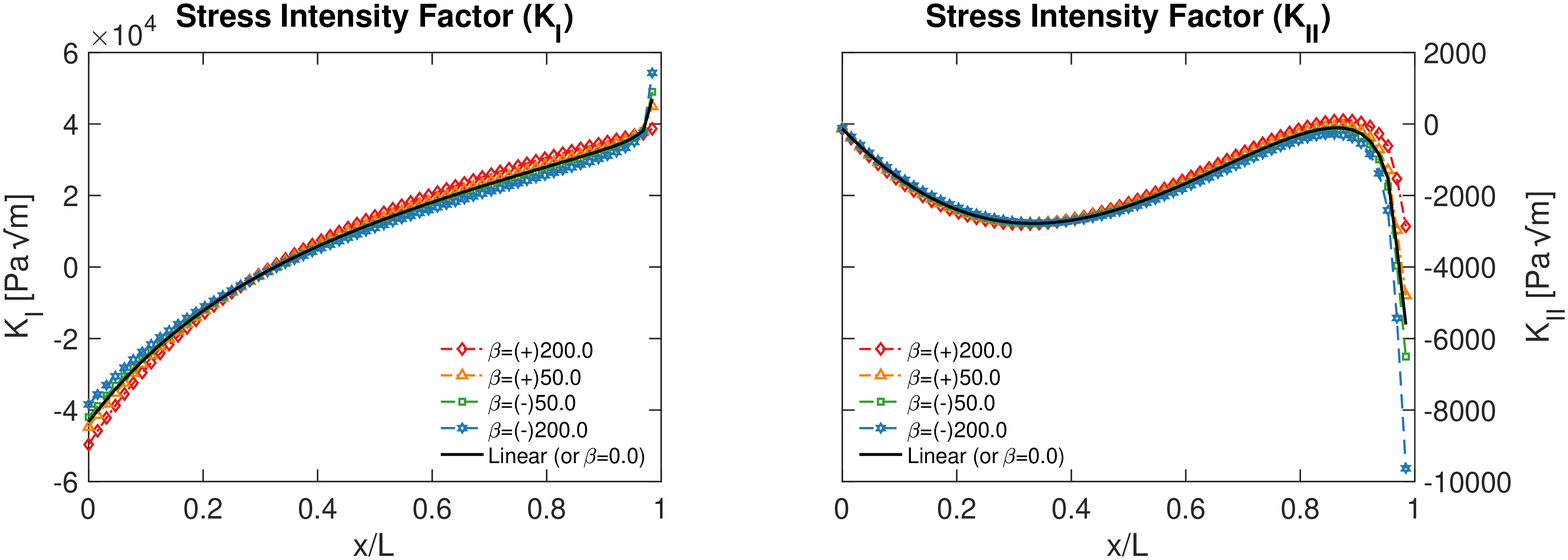}
\caption{[\textbf{Example 2}] {Stress intensity factor {(unit: Pa m$^{1/2}$)}} in mode-I (left, K$_I$) and in mode-II (right, K$_{II}$) on the reference line with {$L=0.5$ m}.}
\label{fig:EX2_SIF} 
\end{figure}
\begin{table}[h]
\centering
\begin{tabular}{|c||l|l||l|l||l|l||}
\hline
\multicolumn{1}{|c||}{\multirow{2}{*}{Variable}} & \multicolumn{2}{c||}{$\beta=-200$}                               & \multicolumn{2}{c||}{$\beta=0$ (Linear)} & \multicolumn{2}{c||}{$\beta=+200$}                            \\ \cline{2-7}                 
& \multicolumn{1}{c|}{Max} & \multicolumn{1}{c||}{Min} & \multicolumn{1}{c|}{Max} & \multicolumn{1}{c||}{Min} & \multicolumn{1}{c|}{Max} & \multicolumn{1}{c||}{Min} \\ \hline\hline                   
 $T_{22}$ [MPa]  & $0.25$ & $-0.022$    & $0.21$  & $-0.025$  & $0.17$  & $-0.028$                    \\ \cline{1-7}
 $\epsilon_{22}$ [ - ]                                  & $0.0014$                & $-0.00023$                    & $0.0019$     & $-0.00025$  & $0.0028$ & $-0.00027$                 \\ \cline{1-7}
 $SED$ {[Pa]}          & 290         & $3.5\times10^{-7}$                    & 330                & $3.8\times10^{-7}$ & 370 & $4.2\times10^{-7}$                 \\ \cline{1-7}
\hline
\end{tabular}
\caption{{[\textbf{Example 2}] The maximum and minimum values of the variables in Figure~\ref{fig:ex2_stress_strain_density} for each case using the linear ($\beta=0$) and nonlinear models with $\beta=-200$ and $\beta=+200$.}}
\label{tab:ex2}
\end{table}

{We also} investigate with the stress intensity factor (SIF) calculated on the reference line in Figure~\ref{Fig:ex2_3_4_setup} (a). The calculation of SIF for the mode-I (K$_I$) is based on the following equation: 
\begin{equation}\label{eq:SIF_mode_1}
\text{K}_{I} :=\lim_{r \to 0}\sqrt{2\pi r}T_{22}(r, \; \theta = 0),
\end{equation}
where $r=0.5$ (equivalently,  x/L~$=0$ in Figure~\ref{fig:EX2_SIF}) located on $\Gamma_4$ in Figure~\ref{Fig:ex2_3_4_setup} (a) and $r=0$ (equivalently, x/L~$=1$ in Figure~\ref{fig:EX2_SIF}) at the tip of crack contacting the reference line. Meanwhile, the stress intensity factor in mode-II loading, i.e., K$_{II}$ is based on the following:
\begin{equation}\label{eq:SIF_mode_2}
\text{K}_{II} :=\lim_{r \to 0}\sqrt{2\pi r}T_{21}(r, \; \theta = 0),
\end{equation}
with the same location for $r$. 
In each example henceforth, we plot these two SIFs together, even though the primary mode of loading is different for each example. {Note that {the negative value implies the compression due to various conditions such as the boundary condition 
and the Poisson's ratio.}}  From Figure~\ref{fig:EX2_SIF} (left), we confirm that its 
maximum K$_I$ occurs right in front of the tip for each case.}} 
Although the difference between the cases is not much discernible, {the calculated SIF of K$_I$} for $\beta=-200$ has the largest 
of at the tip in the same context as the stress 
distributions in Figure~\ref{fig:ex2_stress_strain_density}. 
Interestingly, {the case with $\beta=+200$  has the smallest 
value for the negative region of K$_I$}, {which is analogous to the hardening behavior against the compression.} 
In addition, through K$_{II}$ (Figure~\ref{fig:EX2_SIF} (right)), the shear stress with compressive direction in front of the tip is identified with its negative values for each case. It is also found that $\beta=-200$ has the smallest value in the calculated SIF of K$_{II}$.
{Thus from K$_{I}$ and K$_{II}$, the preferential stiffness is determined based on the nonlinear parameter, i.e., positive and negative $\beta$-values of the density-dependent material moduli model.}
\begin{figure}[!h]
\centering
\includegraphics[width=1.00\textwidth]{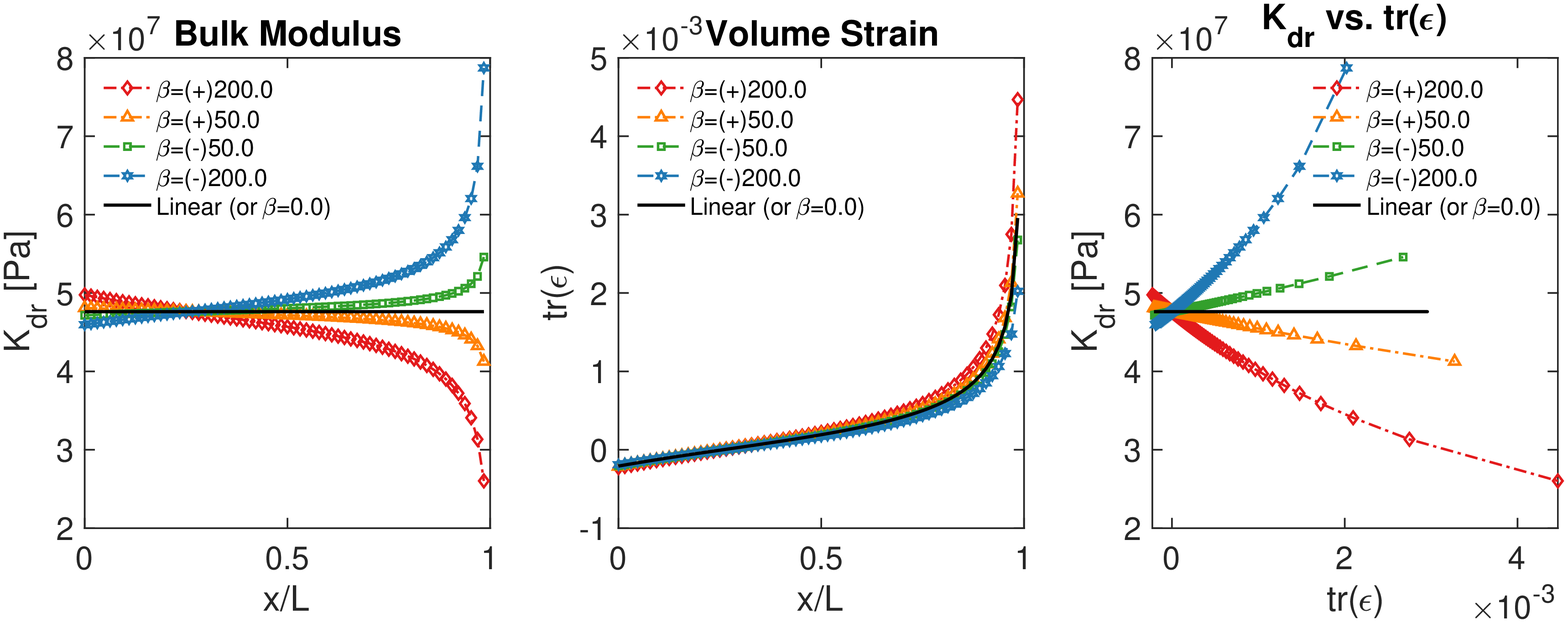}
\caption{[\textbf{Example 2}] Bulk modulus ($K_{dr}${, unit: Pa}) and volumetric strain (tr($\epsilon$)) on the reference line {($L=0.5$ m)}: (left)  $K_{dr}$, (middle) tr($\epsilon$), and (right) $K_{dr}$ vs. tr($\epsilon$).}
\label{fig:EX2_K_trace}
\end{figure} 
\\
\newline
\noindent\textit{\textbf{volumetric strain and bulk modulus} }{Here, we illustrate the change of the drained bulk modulus ($K_{dr}$) in 
\eqref{eq:bulk_modulus} and the volumetric strain ($\tr(\bfeps)$) on the reference line. {Note that the bulk modulus is the inverse of compressibility of the skeleton of the porous solid; greater it has, less compressible it is.}} 
Figure~\ref{fig:EX2_K_trace} demonstrates $K_{dr}$ on the left, $\tr(\bfeps)$ in the middle, and their relation on the right. {As previously figured in K$_I$ (Figure~\ref{fig:EX2_SIF} (left)), we see that a slight compression occurs 
{for each case} from the plot for $\tr(\bfeps)$ shown in the $x$-axis upto around {x/L~$=0.25$} (Figure~\ref{fig:EX2_K_trace} (middle)).} From their relation plot (Figure~\ref{fig:EX2_K_trace} (right)), 
{we see that} the mechanical property (i.e., $K_{dr}$ here) and the intrinsic porosity (i.e., $\tr(\bfeps)$ here) is in its reverse relation for the positive $\beta$ model with some nonlinearity. {In line with the previous finding that the material property of the negative $\beta$-value for the density-dependent model becomes stiffer against the dilation from tensile loading, 
we find from Figure~\ref{fig:EX2_K_trace} (middle) that the volumetric strains become notably differentiable approaching the vicinity of the tip based on the nonlinear parameter; {$\beta=+200$ has the largest for this pure mode-I loading, while $\beta=-200$ has the smallest.} {As compared to the linear model, it is phenomenologically similar to the {strain hardening or softening} in the elastoplasticity regime, {although we do not consider any energy dissipation and the regime of mechanical response remains under the pure elasticity}}.}

\subsubsection{\textbf{Example 3}: {In-plane shear} loading with crack}\label{sec:ex3}
{For \textbf{Example 3}}, we have the pure mode-II loading (see Figure~\ref{Fig:ex2_3_4_setup} (b)), {the volumetric strain changes under the shear force are addressed.} 
The problem has the boundary conditions 
as follows: 
\begin{subequations}\label{bcs-u}
 \begin{align}
  \bfT \bfn &= \bfzero   \quad\text{on} \quad \;  \Gamma_{2}, \; \Gamma_4, \; \Gamma_C, \\
 \bfT \bfn \cdot \bfv &= {T_{21} \, v_1 + T_{22} \, v_2 = f_u \, v_1  \quad \text{on} \quad \Gamma_{3},} \label{trac_bc} \\
 {u_1} &= 0 ~\text{and}~{u_2 = 0}  \quad \text{on} \quad \Gamma_{1}. \label{eq:commom_boundary1}
 \end{align}
 \end{subequations}
 \begin{figure}[!h]
\centering
\subfloat[$T_{21}~\text{with}~\beta=-200$]{\includegraphics[width= 0.3\textwidth]{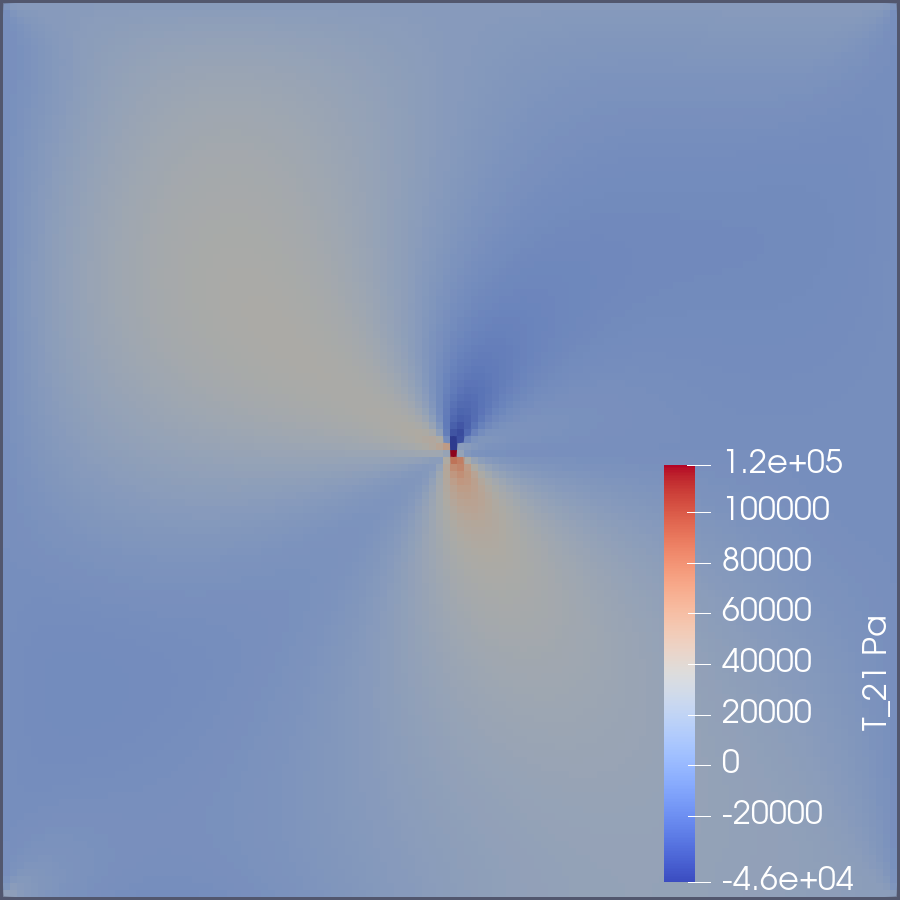}}
\hspace*{0.2in}
\subfloat[$T_{21}~\text{with}~\beta=0~({Linear})$]{\includegraphics[width= 0.3\textwidth]{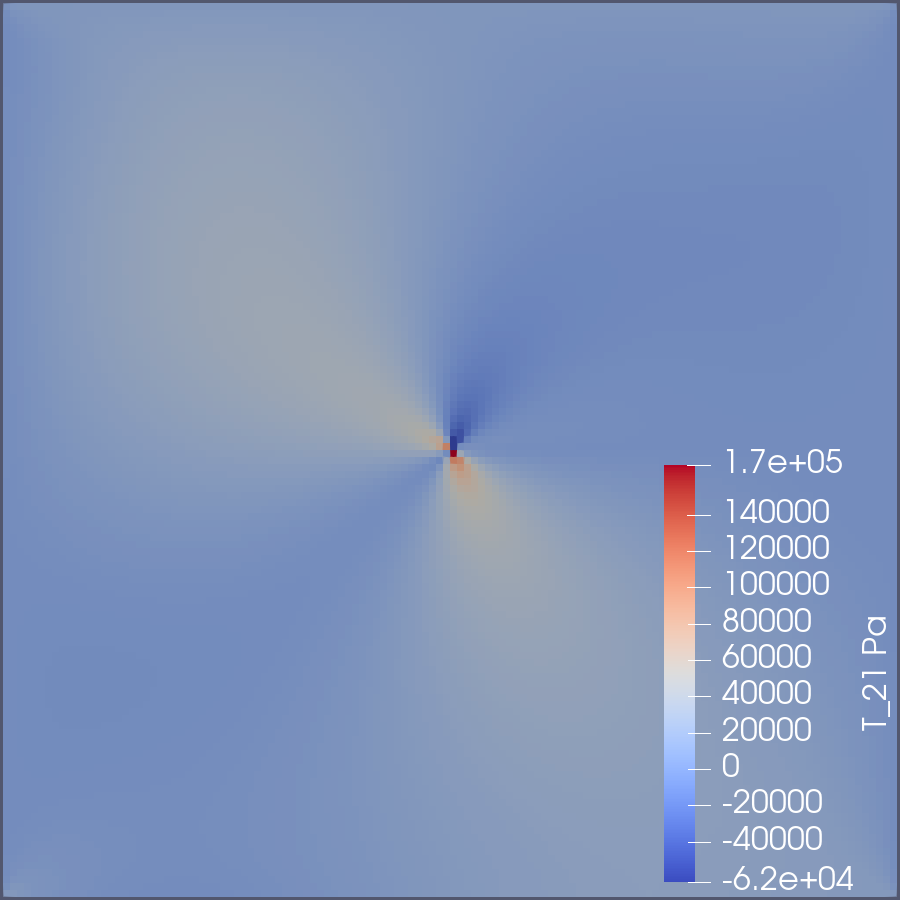}} \hspace*{0.2in}
\subfloat[$T_{21}~\text{with}~\beta=+200$]{\includegraphics[width= 0.3\textwidth]{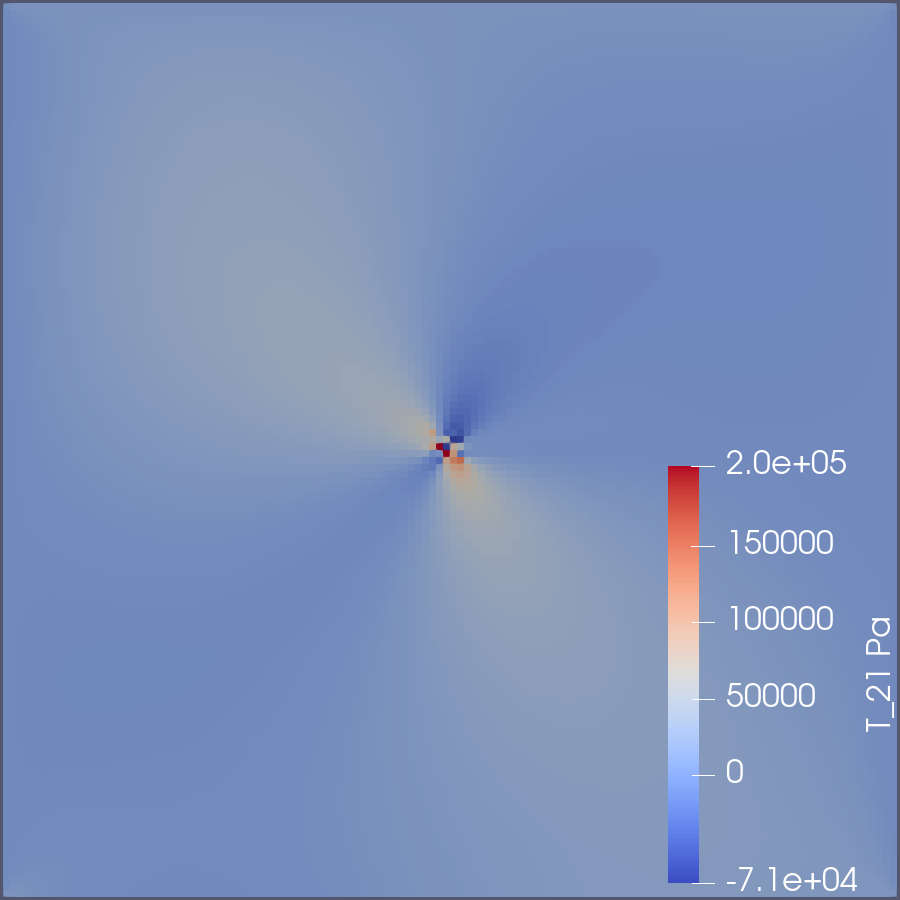}}
\\
\subfloat[$\epsilon_{21}~\text{with}~\beta=-200$]{\includegraphics[width= 0.3\textwidth]{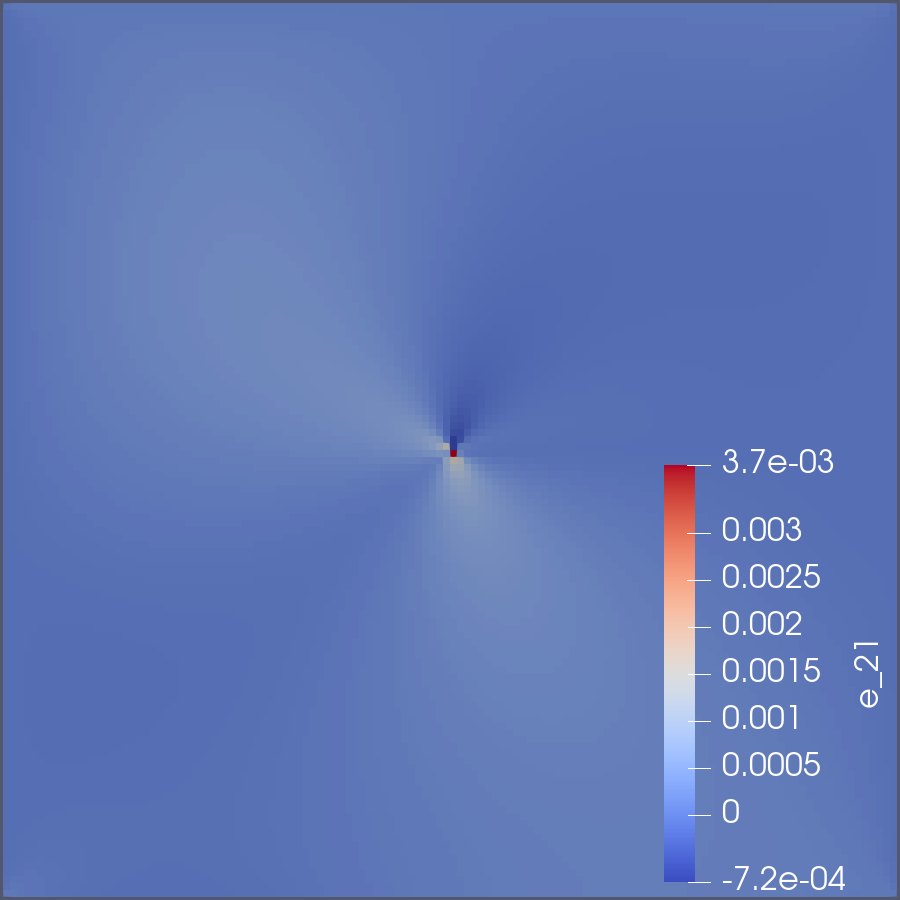}} \hspace*{0.2in}
\subfloat[$\epsilon_{21}~\text{with}~\beta=0~({Linear})$]{\includegraphics[width= 0.3\textwidth]{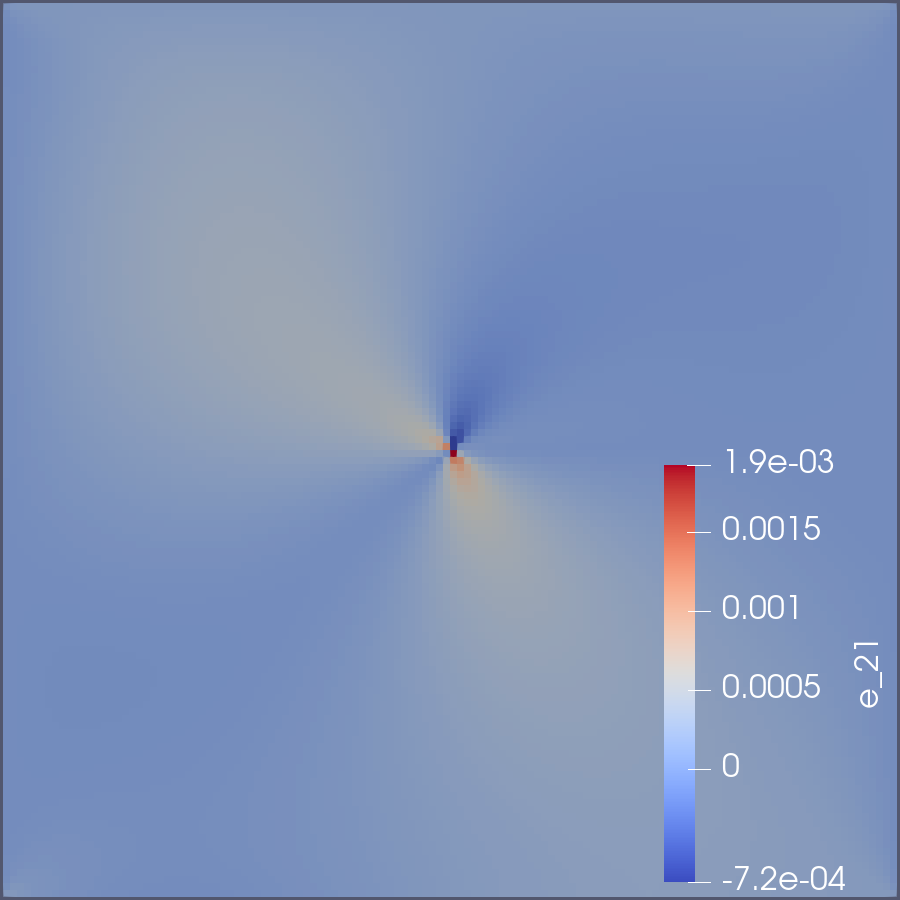}}
\hspace*{0.2in}
\subfloat[$\epsilon_{21}~\text{with}~\beta=+200$]{\includegraphics[width= 0.3\textwidth]{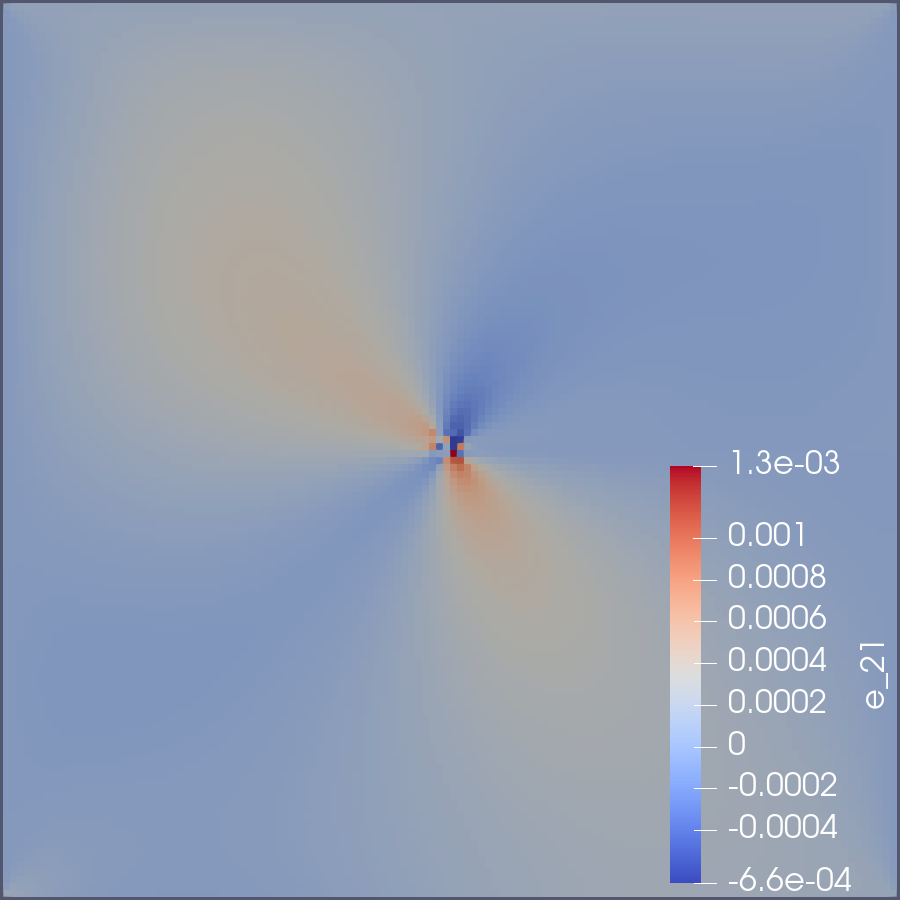}} 
\\
\subfloat[$SED~\text{with}~\beta=-200$]{\includegraphics[width= 0.3\textwidth]{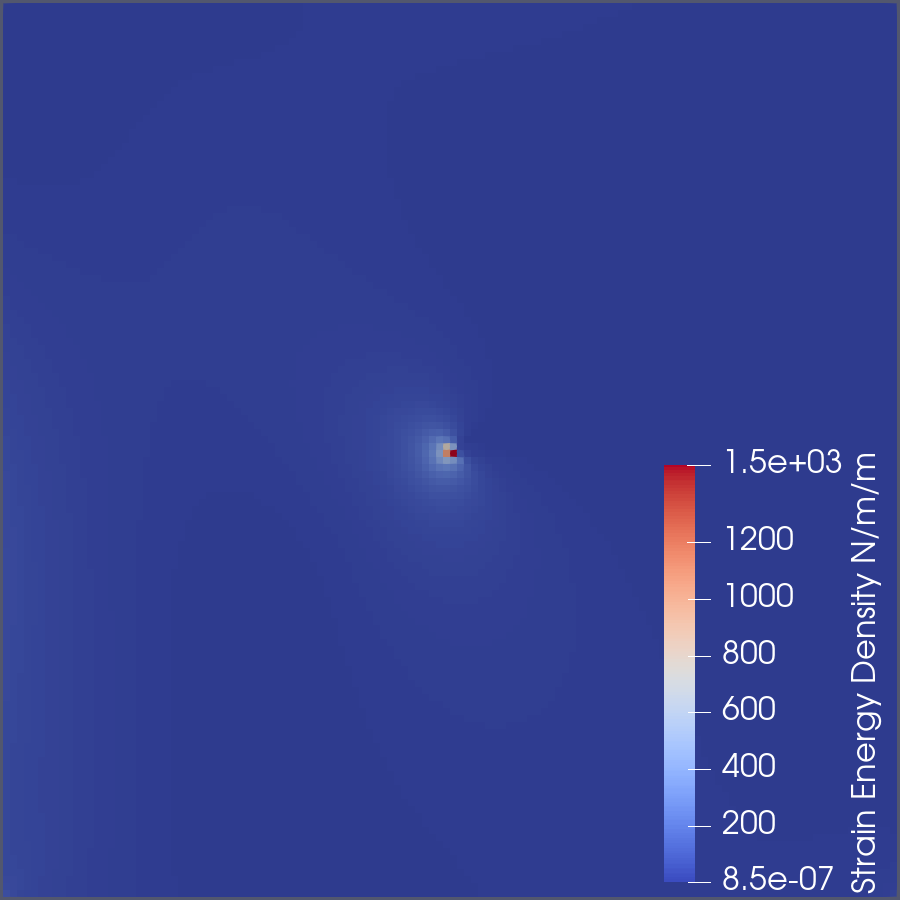}} \hspace*{0.2in}
\subfloat[$SED~\text{with}~\beta=0~({Linear})$]{\includegraphics[width= 0.3\textwidth]{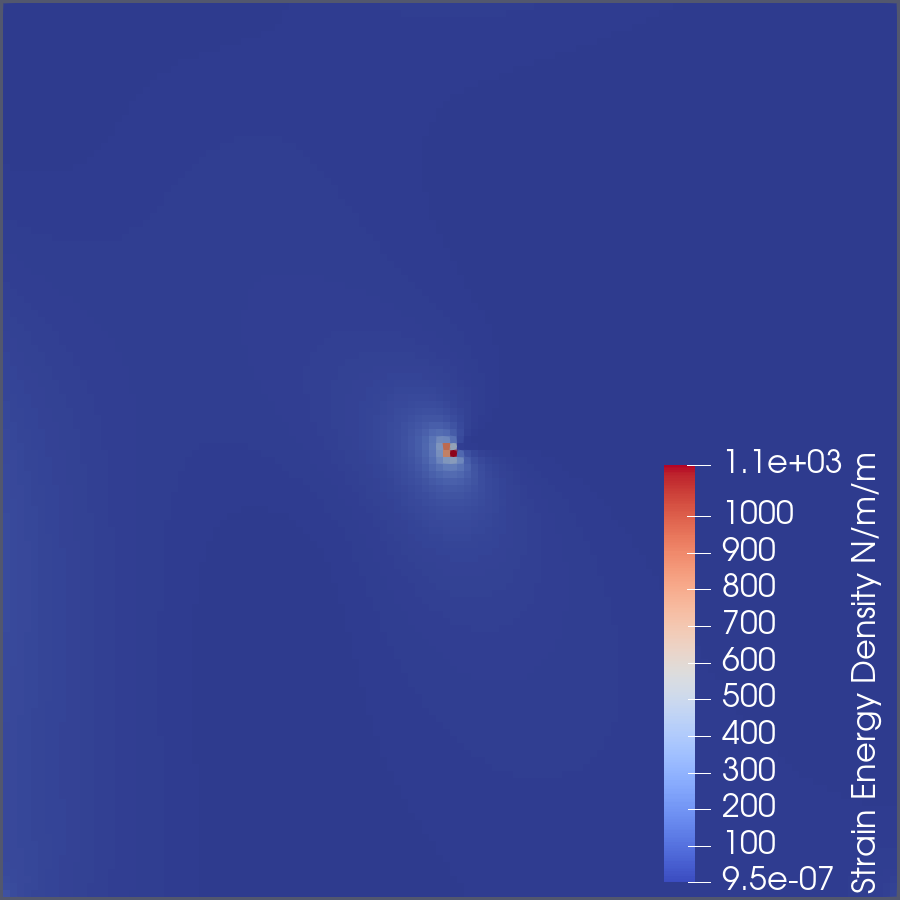}}
\hspace*{0.2in}
\subfloat[$SED~\text{with}~\beta=+200$]{\includegraphics[width= 0.3\textwidth]{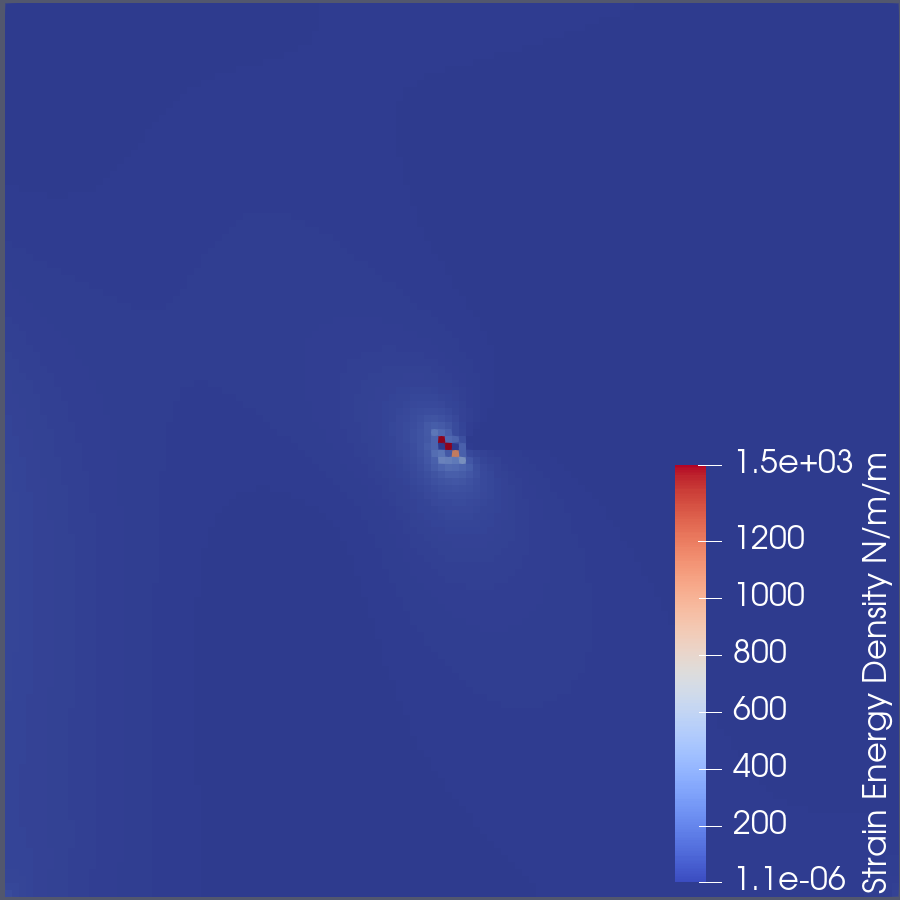}} 
\caption{[\textbf{Example 3}] ((a), (b), (c)) stress $T_{21}$ {(unit: Pa)}, ((d), (e), (f)) strain $\epsilon_{21}$, and ((g), (h), (i)) strain energy density ($SED$, {unit: Pa}) distributions: (left) with $\beta=-200$, (middle) with $\beta=0$ (i.e., the linear model), and (right) {with $\beta=+200$.}}
\label{fig:ex3_stress_strain_density}
\end{figure}

\noindent\textit{\textbf{stress and strain} } For this in-plane shear loading, we plot $T_{21}$, $\epsilon_{21}$ and $SED$ in Figure~\ref{fig:ex3_stress_strain_density} for {the same $\beta$-values}:  $\beta=-200$, $\beta=0$ (the linear), and $\beta=+200$, in the same arrangement  
for the subplots. {Due to the in-plane shear loading}, {$T_{21}$ ((a), (b), and (c) in Figure~\ref{fig:ex3_stress_strain_density}) and $\epsilon_{21}$ ((d), (e), and (f) in Figure~\ref{fig:ex3_stress_strain_density}) have their distributions of positive and negative values across, and this traverse pattern 
exists particularly around the tip. Note that the positive shear strain in the first quadrant of $x$-$y$ plane implies decreasing of the right angle.} {{Unlike \textbf{Example 2}, 
the smallest positive $T_{21}$ is for the $\beta=-200$ case, while the largest belongs to the $\beta=+200$ case}.}
We can see
more clear distinctions between the cases in the distribution of $\epsilon_{21}$.  
As positive and negative shear strains intersect around the tip, compressive stress induces the volume decrease and porosity compaction. 
{In the same context, the case of $\beta=+200$ exhibits much smaller positive values of $\epsilon_{21}$ compared to the case of $\beta=-200$, 
which implies that {the positive $\beta$ case} has more resistance against the compression (or decreasing of the right angle) along with the sliding force applied.}
From the $SED$ results, we also confirm that the maximum strain energy density for each case occurs right in front of the crack-tip. {Likewise, the detailed \textit{max \textnormal{and} min} values for each case are found in Table~\ref{tab:ex3}.}
\begin{figure}[!h]
\centering
\includegraphics[width=1.0\textwidth]{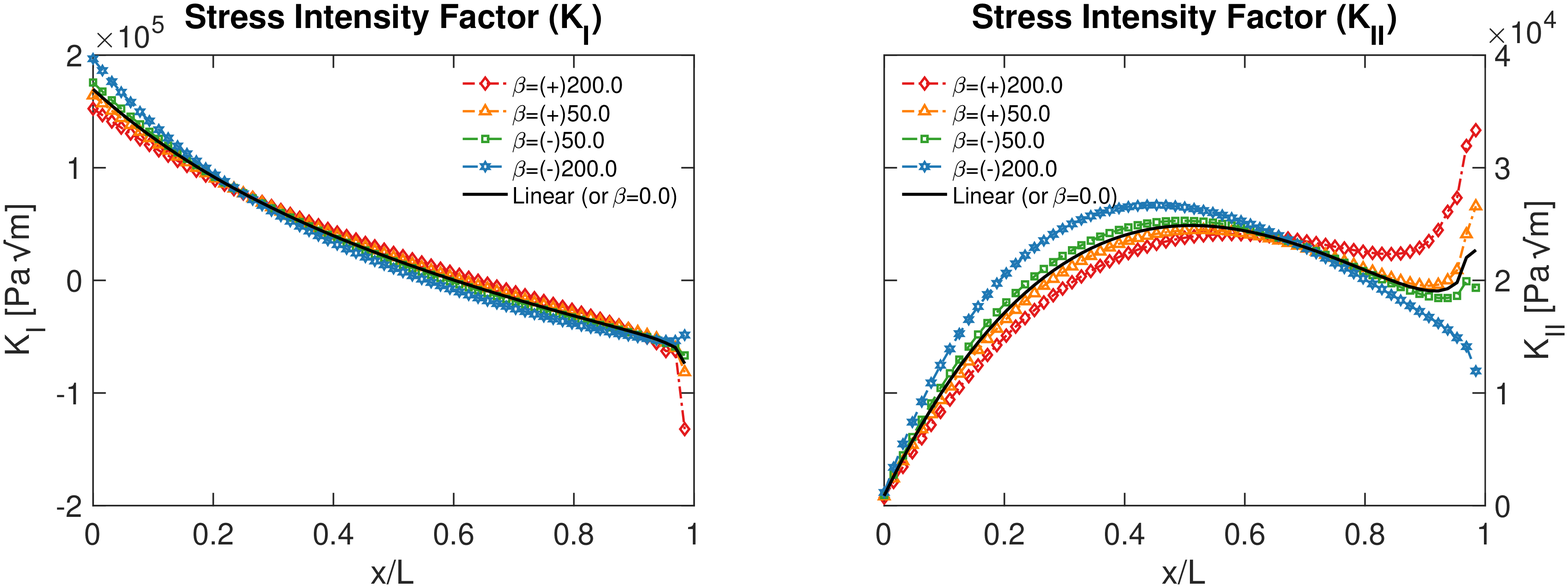}
\caption{{[\textbf{Example 3}] {Stress intensity factor} {(unit: Pa m$^{1/2}$)} in mode-I (left, K$_I$) and in mode-II (right, K$_{II}$) on the reference line with {$L=0.5$ m}.}}
\label{fig:EX3_SIF}
\end{figure}
\begin{table}[h]
\centering
\begin{tabular}{|c||l|l||l|l||l|l||}
\hline
\multicolumn{1}{|c||}{\multirow{2}{*}{Variable}} & \multicolumn{2}{c||}{$\beta=-200$}                               & \multicolumn{2}{c||}{$\beta=0$ (Linear)} & \multicolumn{2}{c||}{$\beta=+200$}                            \\ \cline{2-7}                 
& \multicolumn{1}{c|}{Max} & \multicolumn{1}{c||}{Min} & \multicolumn{1}{c|}{Max} & \multicolumn{1}{c||}{Min} & \multicolumn{1}{c|}{Max} & \multicolumn{1}{c||}{Min} \\ \hline\hline                     
 $T_{21}$ [MPa]  & $0.12$ & $-0.046$    & $0.17$  & $-0.062$  & {$0.20$}  & {$-0.071$}                    \\ \cline{1-7}
 $\epsilon_{21}$ [ - ]                                  & $0.0037$                & $-0.00072$                    & $0.0019$     & $-0.00072$  & {$0.0013$} & {$-0.00064$}                 \\ \cline{1-7}
 $SED$ {[Pa]}          & 1500         & $8.5\times10^{-7}$                    & 1100                & $9.5\times10^{-7}$ & 1500 & $1.1\times10^{-6}$                 \\ \cline{1-7}
\hline
\end{tabular}
\caption{{[\textbf{Example 3}] The maximum and minimum values of the variables in Figure~\ref{fig:ex3_stress_strain_density} for each case using the linear ($\beta=0$) and nonlinear models with $\beta=-200$ and $\beta=+200$.}}
\label{tab:ex3}
\end{table}

Two SIFs on the reference line (Figure~\ref{Fig:ex2_3_4_setup} (b)) are illustrated in Figure~\ref{fig:EX3_SIF} using \eqref{eq:SIF_mode_1} and \eqref{eq:SIF_mode_2}.
We find that at about x/L~$=0.6$ in the reference line, K$_{II}$ for the case of $\beta=+200$ surpasses the rest cases. 
In parallel, the case of $\beta=+200$ has the smallest concentration for K$_{I}$ (Figure~\ref{fig:EX3_SIF} (left)) in front of the tip  {with the negative tensile stress, i.e., under compression. Therefore for the in-plane shear loading, we find the model with negative $\beta$-values has the weaker material strength against the shear or compression.} 
\begin{figure}[!h]
\centering
\includegraphics[width=1.00\textwidth]{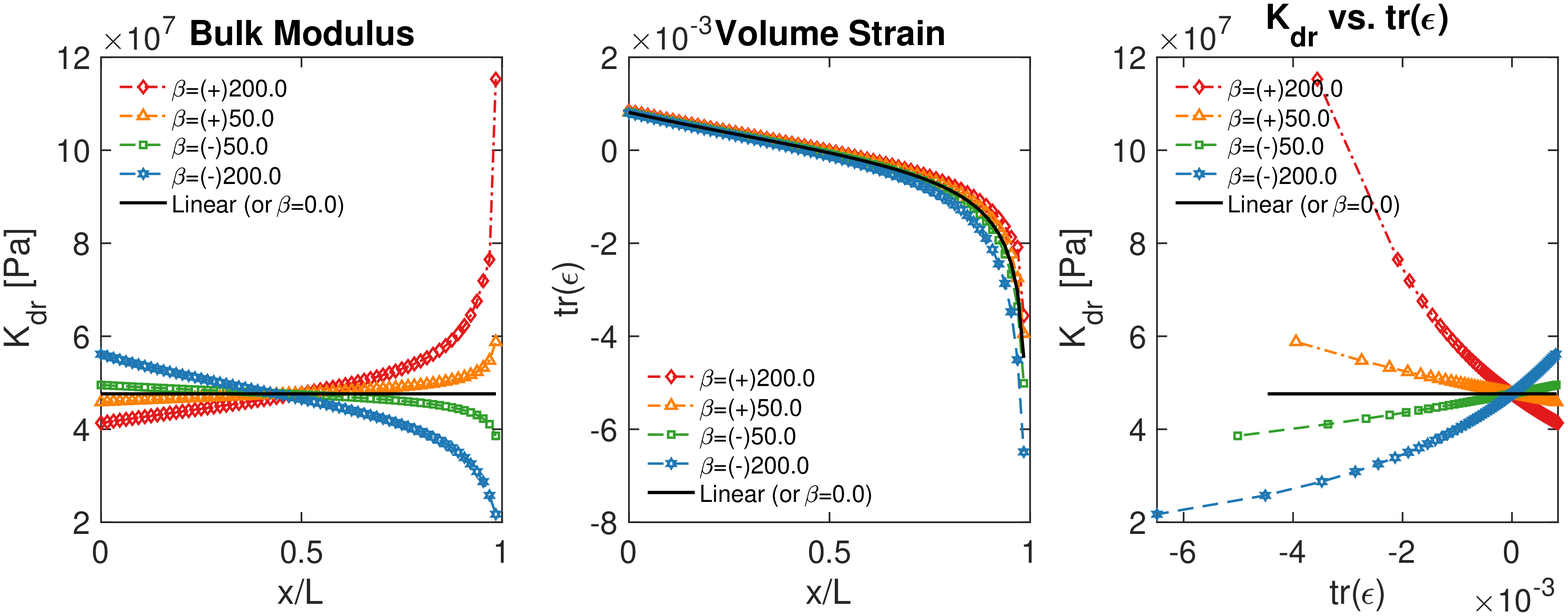}
\caption{{[\textbf{Example 3}] Bulk modulus ($K_{dr}$, unit: Pa) and volumetric strain (tr($\epsilon$)) on the reference line ({$L=0.5$ m}): (left)  {$K_{dr}$}, {(middle) tr($\epsilon$)}, and (right) $K_{dr}$ vs. tr($\epsilon$).}}
\label{fig:EX3_K_trace}
\end{figure}
\\
\newline
\noindent\textit{\textbf{volumetric strain and bulk modulus} }
In Figure~\ref{fig:EX3_K_trace}, the bulk modulus and volumetric strain for each case
are plotted. We confirm that the stiffness of each case is illustrated in the opposite compared to \textbf{Example 2} (see Figure~\ref{fig:EX2_K_trace}), as the volumetric strains are plotted in the reverse direction. Until about the half of reference line (x/L~$=0.5$) (see Figure~\ref{fig:EX3_K_trace} (middle)), positive tr($\epsilon$), i.e., the dilation of implicit porosity, is decreased {to zero} 
in each case, 
which can also be figured in Figure~\ref{fig:EX3_K_trace} (left) with $K_{dr}$. 
Note that more variation of tr($\epsilon$) with larger range is observed for the case of $\beta=-200$, particularly at the tip, {thus we confirm that {the density-dependent model with negative $\beta$} is relatively weak in the shear (or compressive) loading. 
}

\subsubsection{\textbf{Example 4}: {Mixed-mode} loading with crack}\label{sec:ex4}
In this last example, the mixed-mode (i.e., the mode-I and II) of loading  is applied. 
The boundary conditions 
are as follows:
\begin{subequations}\label{bcs-u}
 \begin{align}
 \bfT \bfn &= \bfzero  \quad\text{on} \quad \;  \Gamma_{2}, \; \Gamma_4, \; \Gamma_C, \\
 \bfT \bfn \cdot \bfv &=  {T_{21} \, v_1 + T_{22} \, v_2= f_u \, v_1+f_u \, v_2  \quad \text{on} \quad \Gamma_{3},}  \label{trac_bc} \\
 {u_1} &= 0 ~\text{and}~{u_2 = 0}  \quad \text{on} \quad \Gamma_{1}. \label{eq:commom_boundary1}
 \end{align}
 \end{subequations}
Note that we have the same bottom boundary ($\Gamma_{1}$) condition as \textbf{Example 3} (see Figure~\ref{Fig:ex2_3_4_setup} (c)) considering the mode-II loading.
\\
\newline
\noindent\textit{\textbf{stress and strain} }We plot both $T_{22}$ with $\epsilon_{22}$ and $T_{21}$ with $\epsilon_{21}$ in Figure~\ref{fig:ex4_stress_strain_density_22} and \ref{fig:ex4_stress_strain_density_21}, respectively. For $T_{22}$, unlike the pure mode-I loading, we find compressive stress near the tip, which is due to the {in-plane} shear or mode-II loading imposed simultaneously. 
As consistent with the fact that {the positive $\beta$} case has more resistance against the compression and less 
against the tension, the smallest negative $\epsilon_{22}$ is obtained for the case of $\beta=+200$ as seen in Figure~\ref{fig:ex4_stress_strain_density_22}. In addition, we figure the smallest positive $\epsilon_{22}$ is for $\beta=-200$ as the {negative $\beta$} case is relatively stiffer against the tensile loading. 
{The same pattern of distribution is also found; larger a variable has for the value, narrower is its distribution with more localization.}  {As the positive $\beta$ has more resistance against the compression, compressive stress of the case of $\beta=+200$ is more concentrated with larger value {(colored in blue in Figure~\ref{fig:ex4_stress_strain_density_22} (c))}. 
While for strain ($\epsilon_{22}$), the  
case of $\beta=-200$ has wider tensile strain region with smaller positive values {(colored in red in Figure~\ref{fig:ex4_stress_strain_density_22} (a))}, {i.e., more resistance}. 
\begin{figure}[!h]
\centering
\subfloat[$T_{22}~\text{with}~\beta=-200$]{\includegraphics[width= 0.3\textwidth]{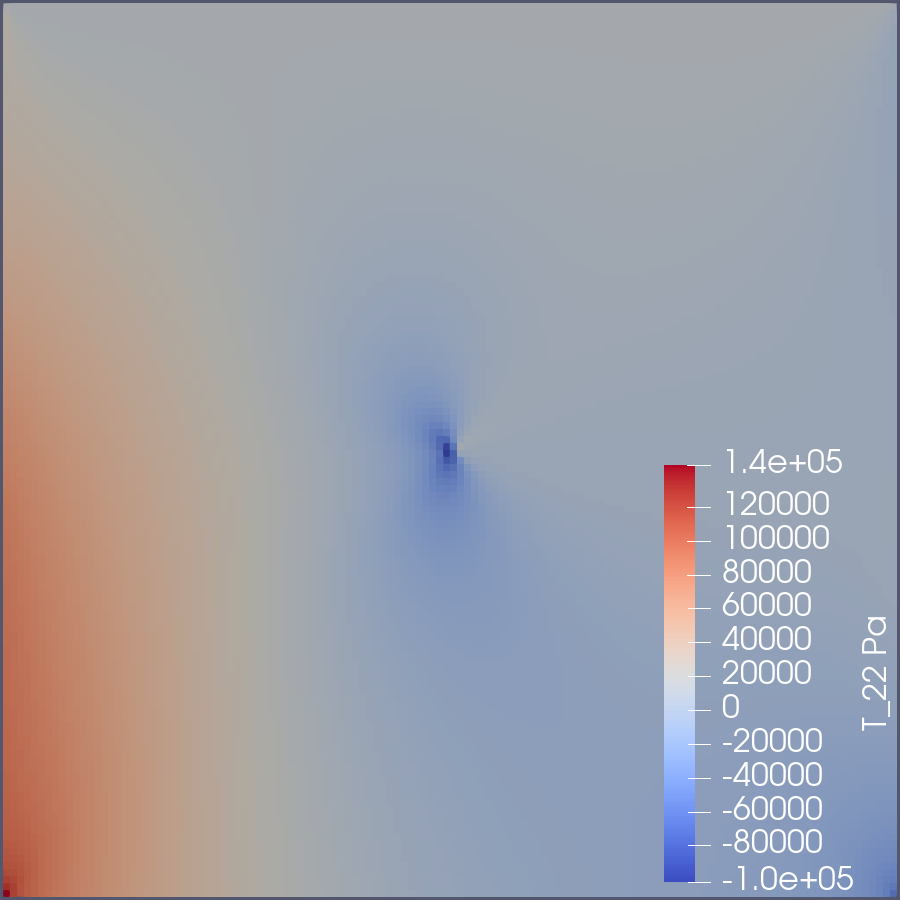}}
\hspace*{0.2in}
\subfloat[$T_{22}~\text{with}~\beta=0~({Linear})$]{\includegraphics[width= 0.3\textwidth]{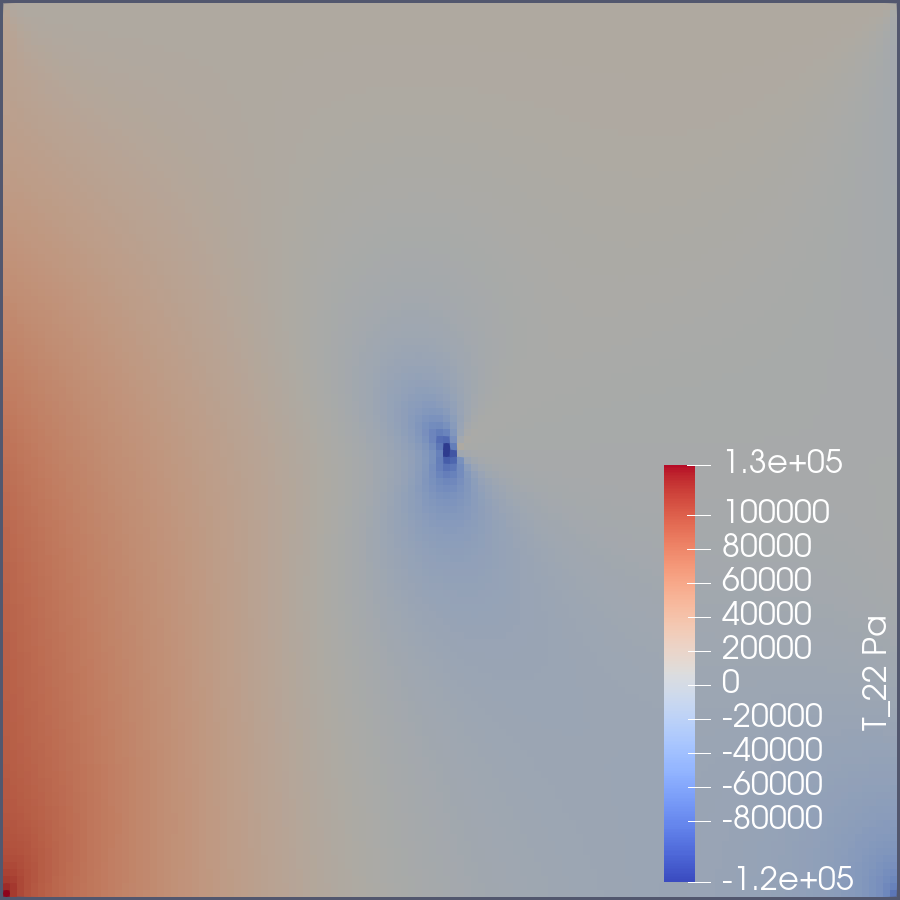}} \hspace*{0.2in}
\subfloat[$T_{22}~\text{with}~\beta=+200$]{\includegraphics[width= 0.3\textwidth]{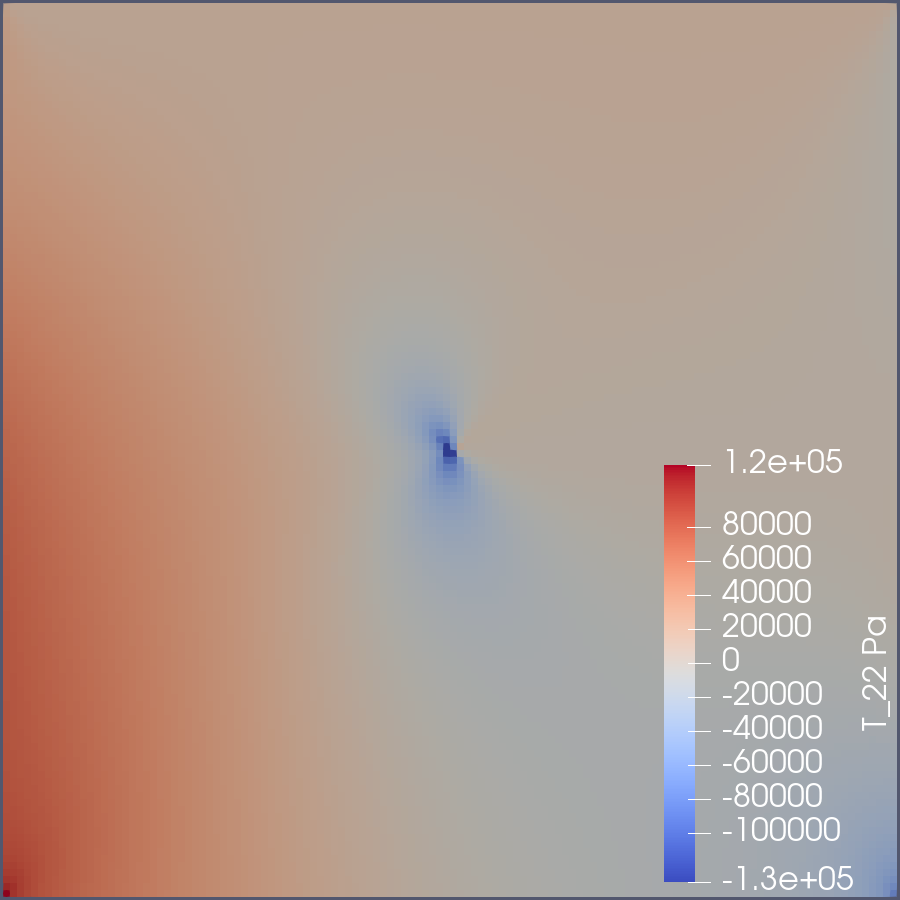}}
\\
\subfloat[$\epsilon_{22}~\text{with}~\beta=-200$]{\includegraphics[width= 0.3\textwidth]{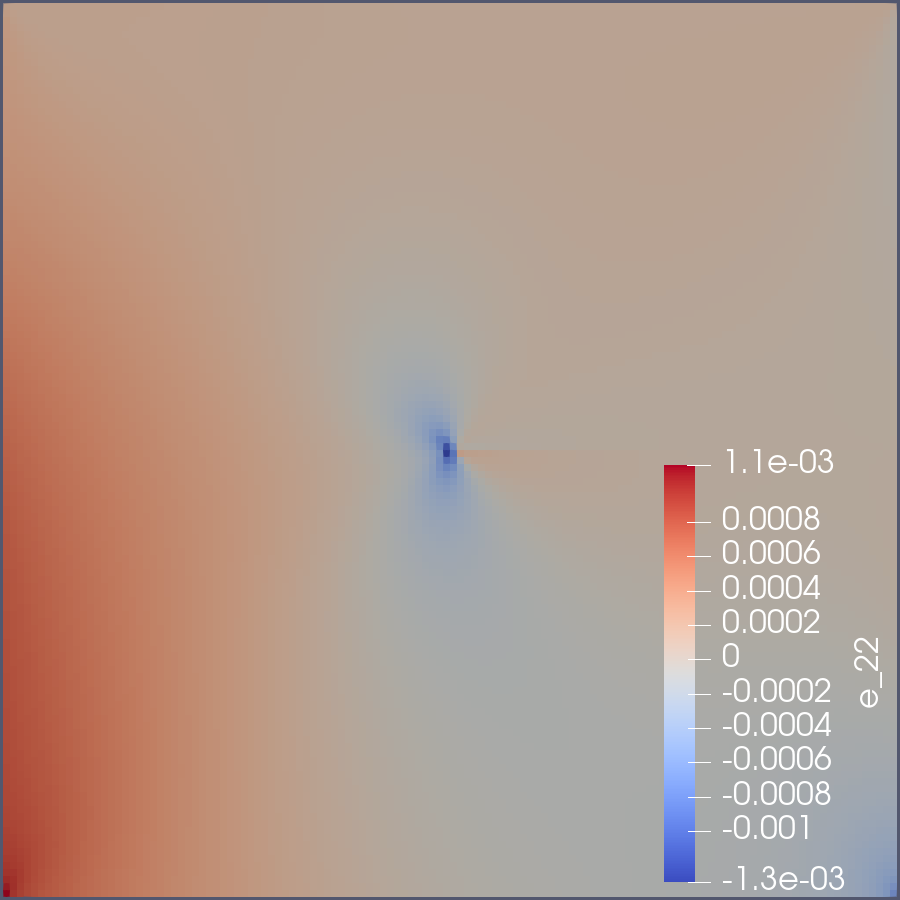}} \hspace*{0.2in}
\subfloat[$\epsilon_{22}~\text{with}~\beta=0~({Linear})$]{\includegraphics[width= 0.3\textwidth]{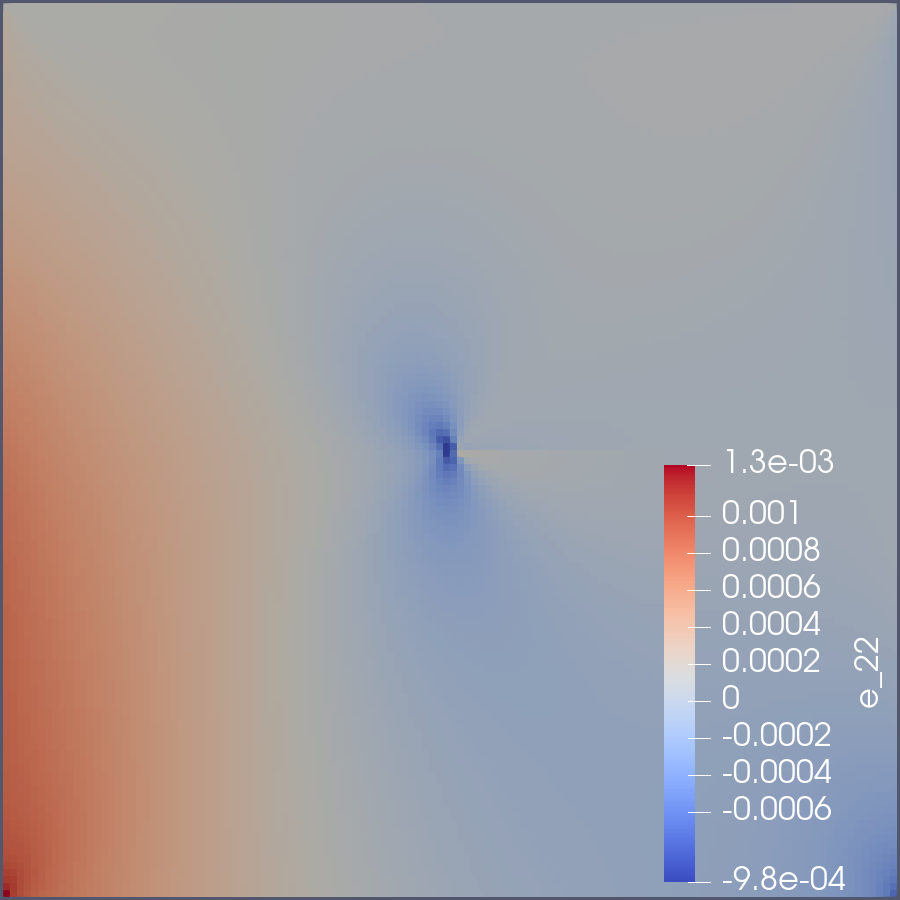}}
\hspace*{0.2in}
\subfloat[$\epsilon_{22}~\text{with}~\beta=+200$]{\includegraphics[width= 0.3\textwidth]{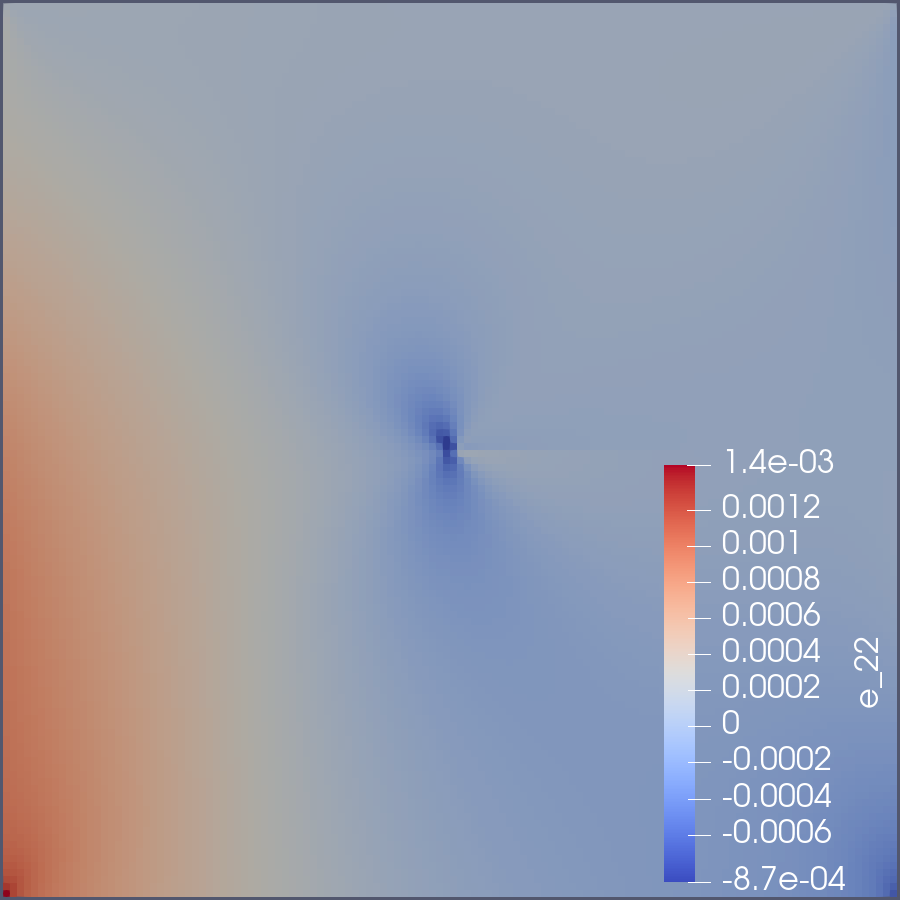}} 
\caption{[\textbf{Example 4}] ((a), (b), (c)) stress $T_{22}$ {(unit: Pa)} and ((d), (e), (f)) strain $\epsilon_{22}$ distributions: (left) with $\beta=-200$, (middle) with $\beta=0$ (i.e., the linear model), and (right) with $\beta=+200$.}
\label{fig:ex4_stress_strain_density_22}
\end{figure}
About $T_{21}$ with $\epsilon_{21}$ in Figure~\ref{fig:ex4_stress_strain_density_21}, similar 
patterns to \textbf{Example 3} with the pure mode-II loading are shown; the positive $\beta$ with the case of $\beta=+200$ is stiffer against the compression with the {smallest positive value in $\epsilon_{21}$} and the largest negative $\epsilon_{21}$. We see that the maximum strain energy density ($SED$) is 
shown in front of the tip for each case (Figure~\ref{fig:ex4_stress_strain_density_21}). {See the detailed \textit{max \textnormal{and} min} values for each case in Table~\ref{tab:ex4}.}
\begin{figure}[!h]
\centering
\subfloat[$T_{21}~\text{with}~\beta=-200$]{\includegraphics[width= 0.3\textwidth]{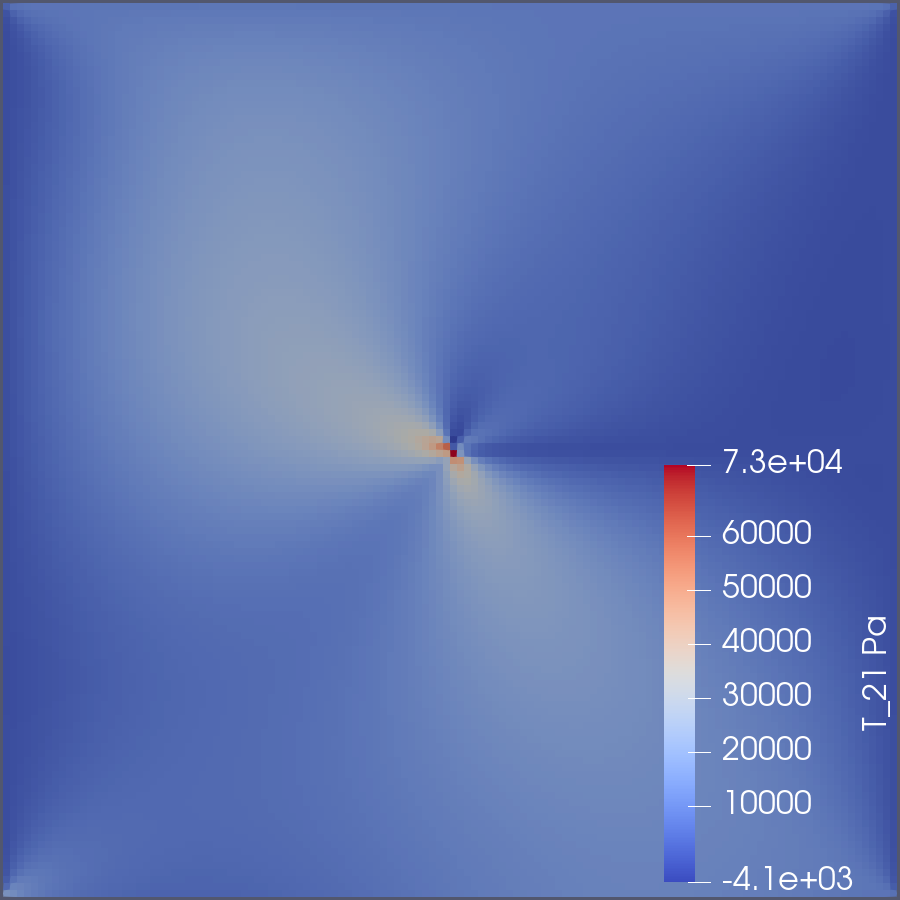}}
\hspace*{0.2in}
\subfloat[$T_{21}~\text{with}~\beta=0~({Linear})$]{\includegraphics[width= 0.3\textwidth]{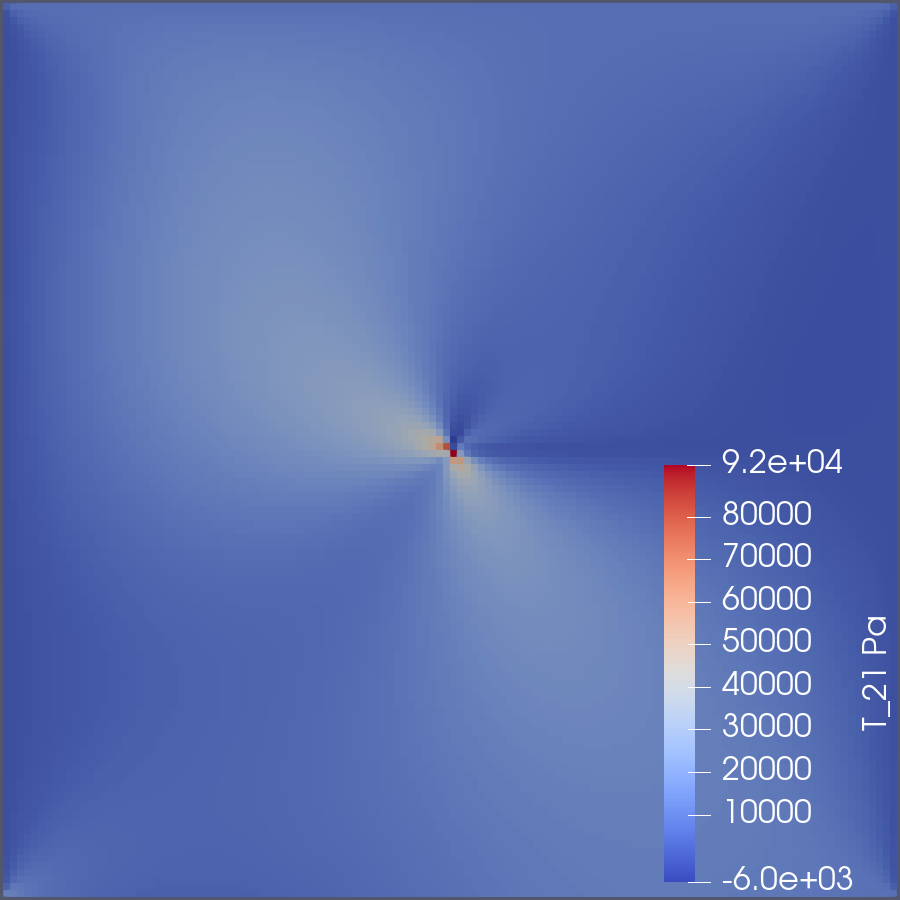}} \hspace*{0.2in}
\subfloat[$T_{21}~\text{with}~\beta=+200$]{\includegraphics[width= 0.3\textwidth]{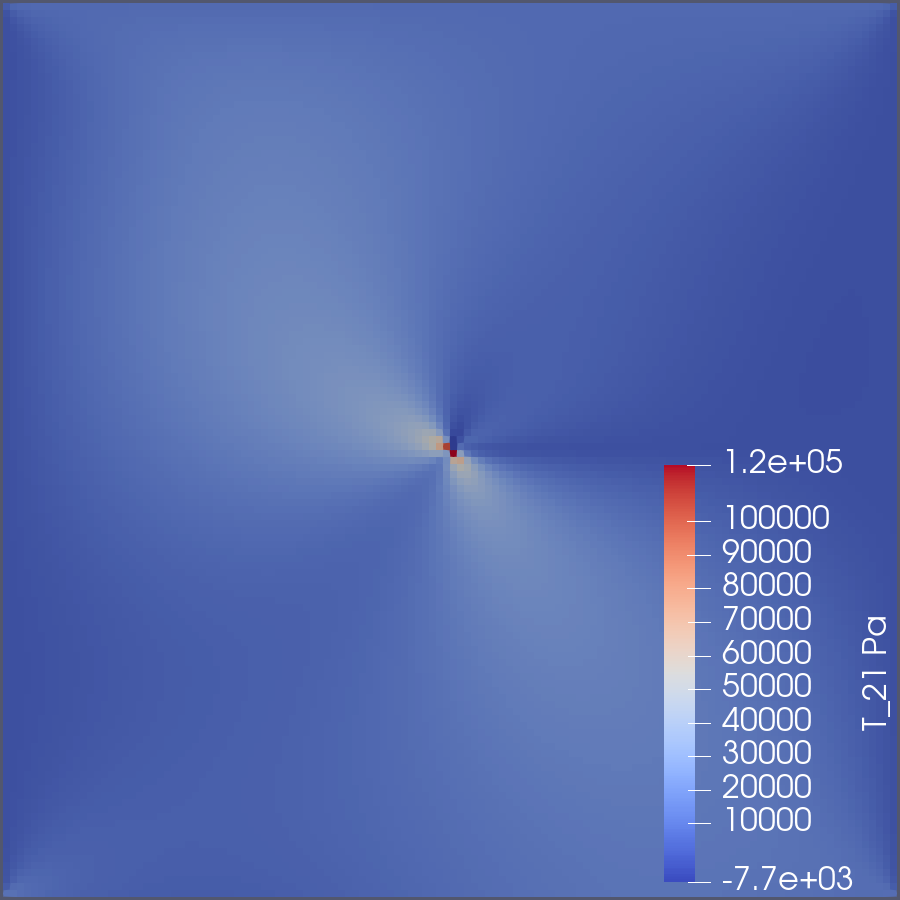}}
\\
\subfloat[$\epsilon_{21}~\text{with}~\beta=-200$]{\includegraphics[width= 0.3\textwidth]{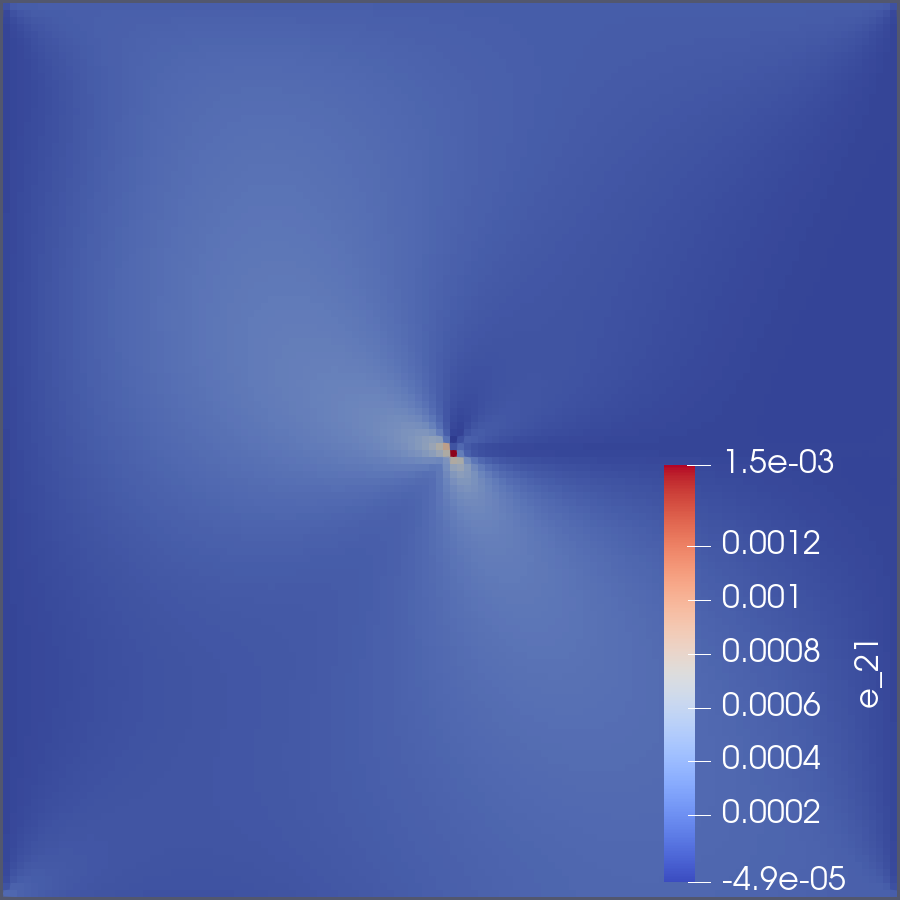}} \hspace*{0.2in}
\subfloat[$\epsilon_{21}~\text{with}~\beta=0~({Linear})$]{\includegraphics[width= 0.3\textwidth]{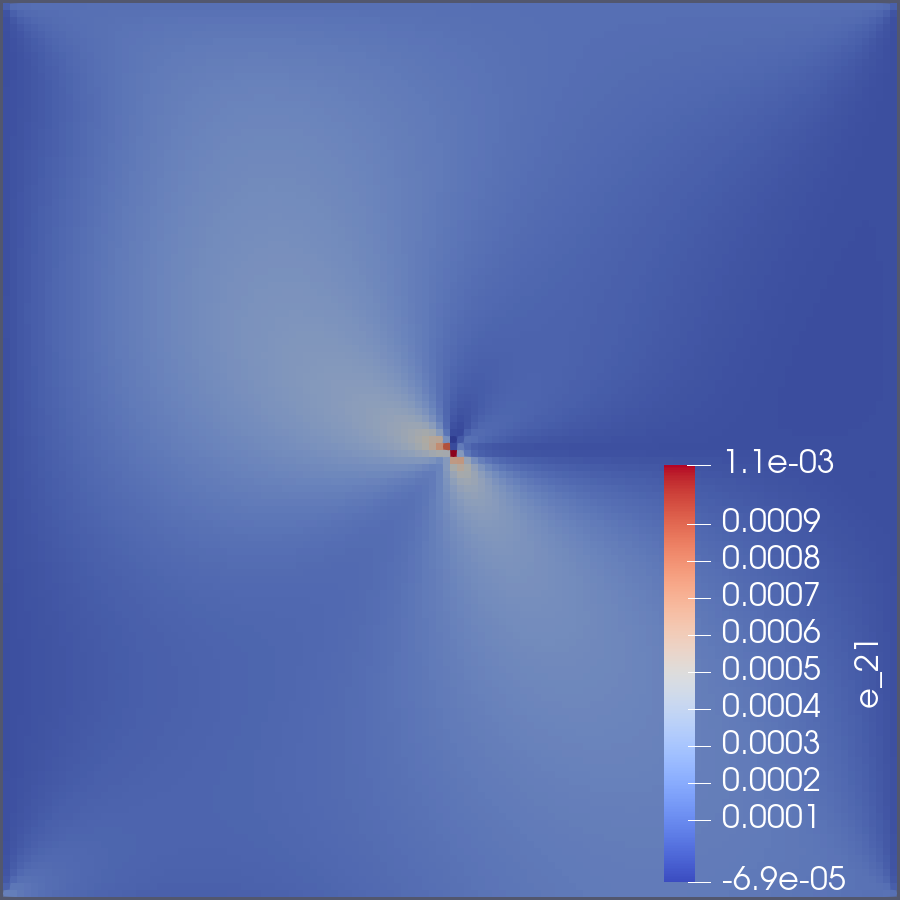}}
\hspace*{0.2in}
\subfloat[$\epsilon_{21}~\text{with}~\beta=+200$]{\includegraphics[width= 0.3\textwidth]{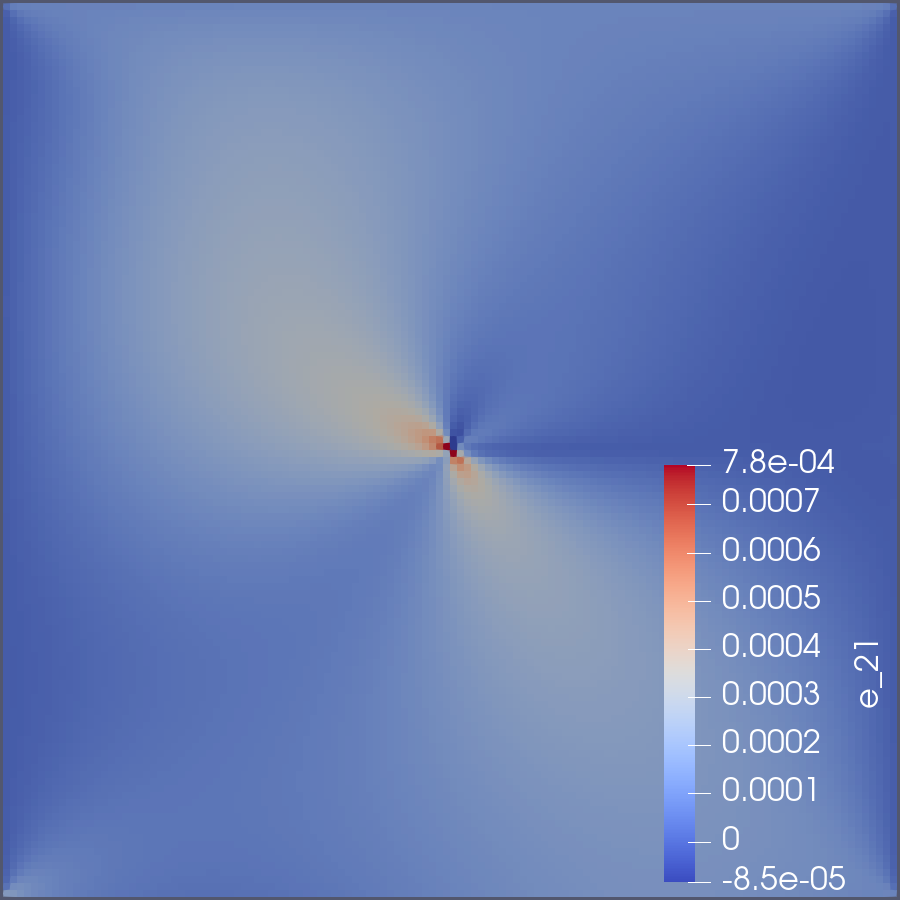}}
\\
\subfloat[$SED~\text{with}~\beta=-200$]{\includegraphics[width= 0.3\textwidth]{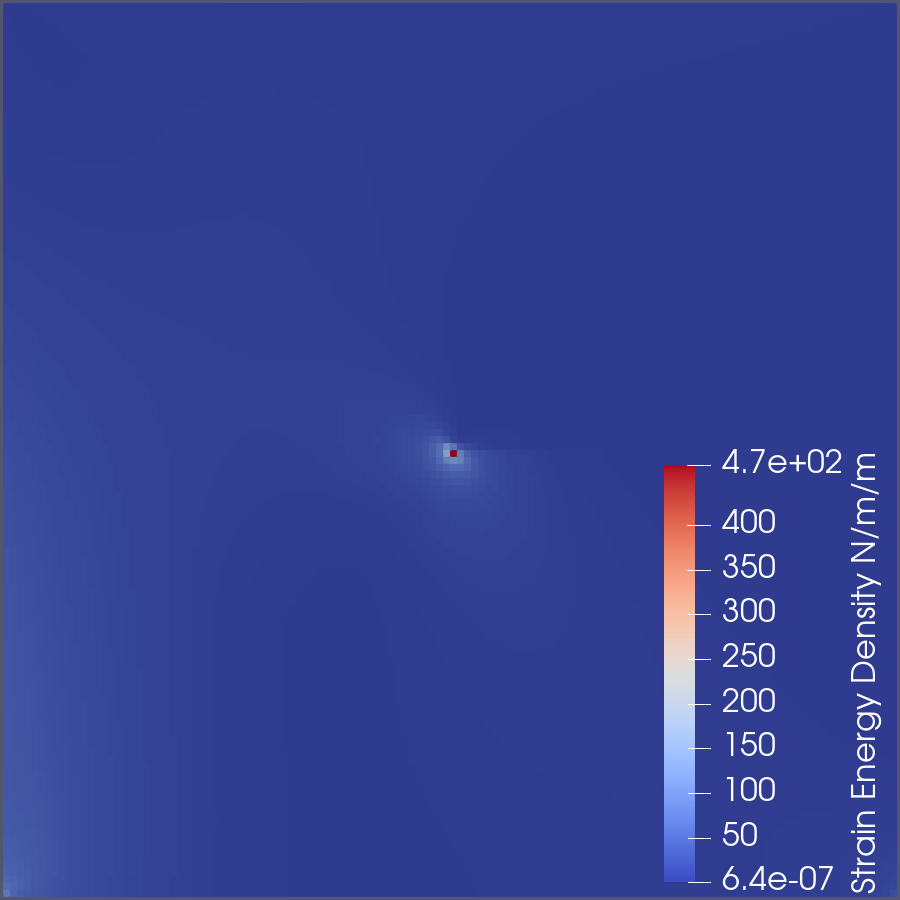}} \hspace*{0.2in}
\subfloat[$SED~\text{with}~\beta=0~({Linear})$]{\includegraphics[width= 0.3\textwidth]{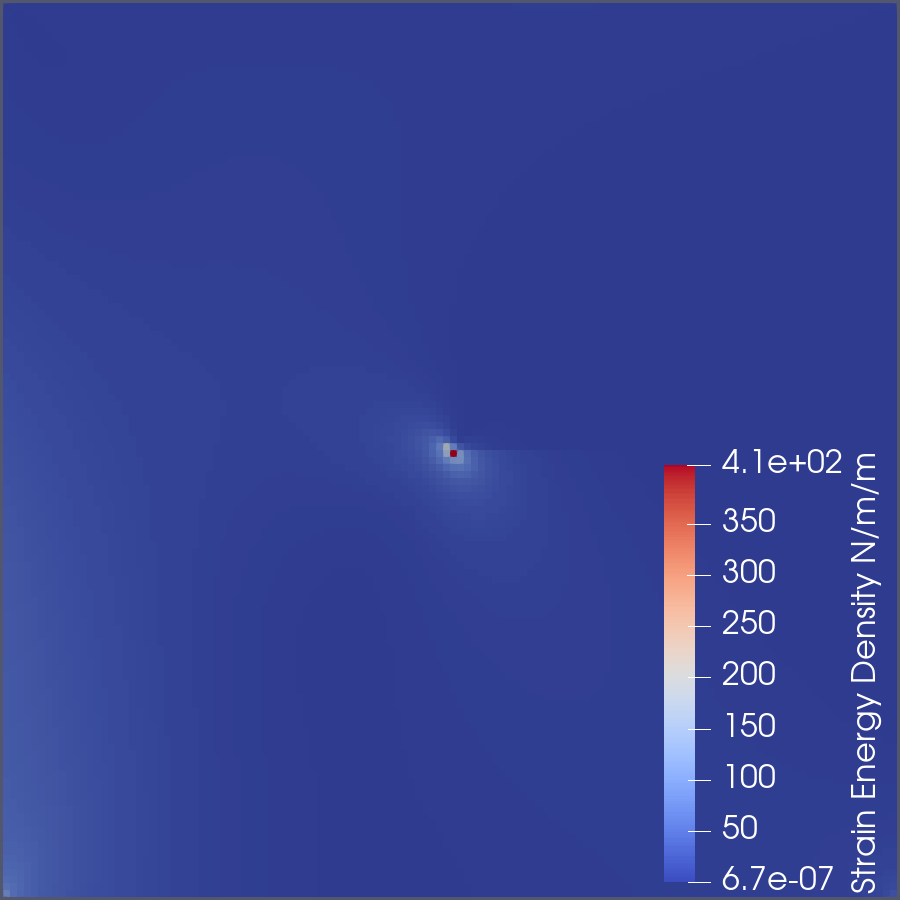}}
\hspace*{0.2in}
\subfloat[$SED~\text{with}~\beta=+200$]{\includegraphics[width= 0.3\textwidth]{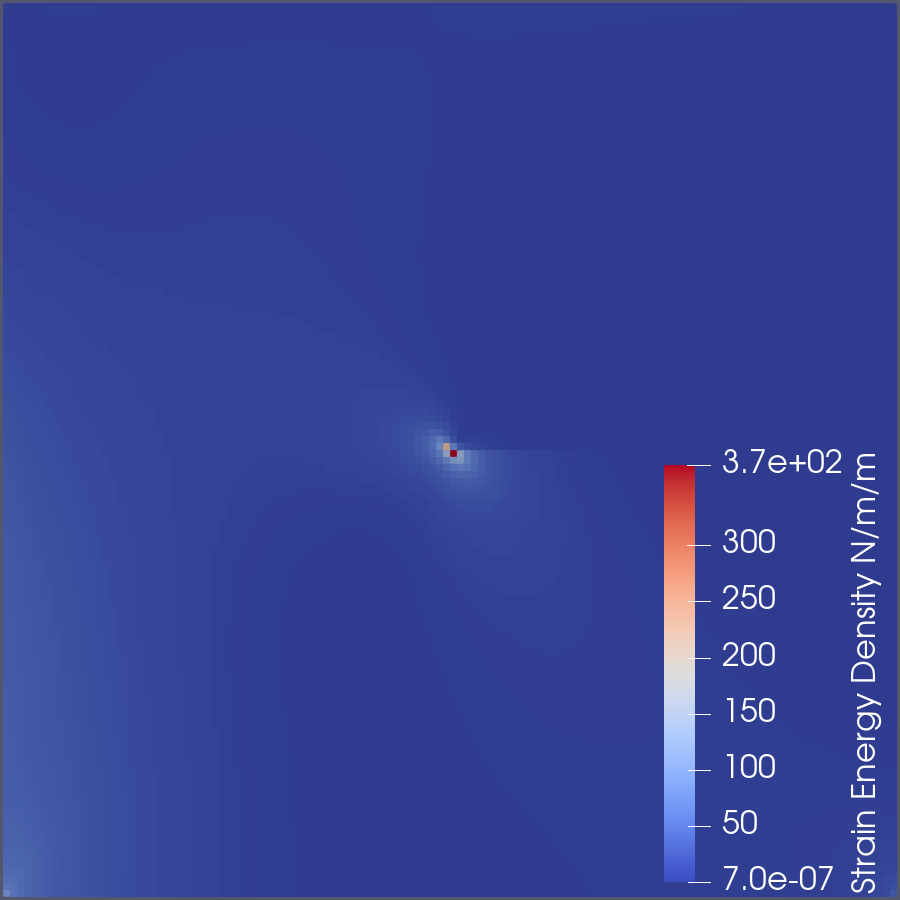}} 
\caption{[\textbf{Example 4}] ((a), (b), (c)) stress $T_{21}$ {(unit: Pa)}, ((d), (e), (f)) strain $\epsilon_{21}$, and ((g), (h), (i)) strain energy density ($SED$, {unit: Pa}) distributions: (left) with $\beta=-200$, (middle) with $\beta=0$ (i.e., the linear model), and (right) with $\beta=+200$.}
\label{fig:ex4_stress_strain_density_21}
\end{figure}
\begin{table}[h]
\centering
\begin{tabular}{|c||l|l||l|l||l|l||}
\hline
\multicolumn{1}{|c||}{\multirow{2}{*}{Variable}} & \multicolumn{2}{c||}{$\beta=-200$}                               & \multicolumn{2}{c||}{$\beta=0$ (Linear)} & \multicolumn{2}{c||}{$\beta=+200$}                            \\ \cline{2-7}                 
& \multicolumn{1}{c|}{Max} & \multicolumn{1}{c||}{Min} & \multicolumn{1}{c|}{Max} & \multicolumn{1}{c||}{Min} & \multicolumn{1}{c|}{Max} & \multicolumn{1}{c||}{Min} \\ \hline\hline   
 $T_{22}$ [MPa]  & $0.14$ & $-0.10$    & $0.13$  & $-0.12$  & $0.12$  & $-0.13$                    \\ \cline{1-7}
 $\epsilon_{22}$ [ - ]                                  & $0.0011$                & $-0.0013$                    & $0.0013$     & $-0.00098$  & $0.0014$ & $-0.00087$                 \\ \hline\hline                
 $T_{21}$ [MPa]  & $0.073$ & $-0.0041$    & $0.092$  & $-0.006$  & {$0.12$}  & $-0.0077$                    \\ \cline{1-7}
 $\epsilon_{21}$ [ - ]                                  & $0.0015$                & $-0.000049$                    & $0.0011$     & $-0.000069$  & $0.00078$ & {$-0.000085$}                 \\ \cline{1-7}
 $SED$ {[Pa]}          & 470         & $6.4\times10^{-7}$                    & 410                & $6.7\times10^{-7}$ & 370 & $7.0\times10^{-7}$                 \\ \cline{1-7}
\hline
\end{tabular}
\caption{{[\textbf{Example 4}] The maximum and minimum values of the variables in Figure~\ref{fig:ex4_stress_strain_density_22} and \ref{fig:ex4_stress_strain_density_21} for each case using the linear ($\beta=0$) and nonlinear models with $\beta=-200$ and $\beta=+200$.}}
\label{tab:ex4}
\end{table}
\begin{figure}[h]
\centering
\includegraphics[width=1.00\textwidth]{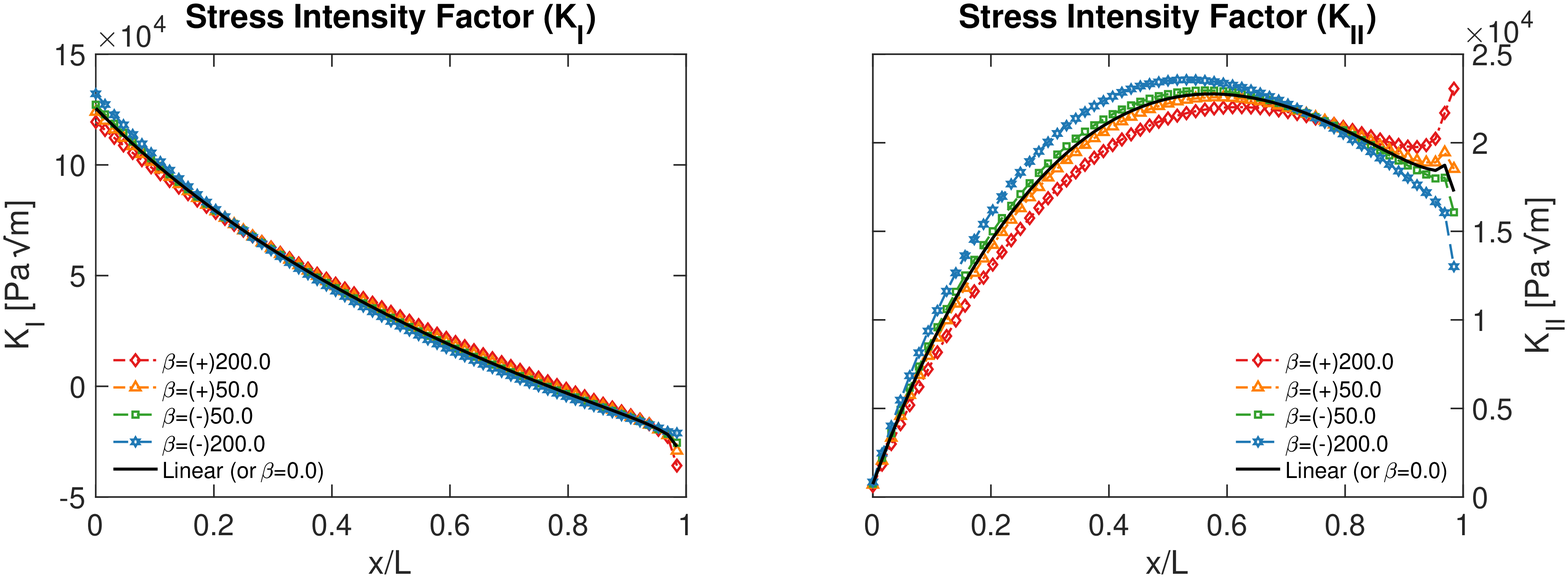}
\caption{[\textbf{Example 4}] Stress intensity factor {(unit: Pa m$^{1/2}$)} in mode-I (left, K$_{I}$) and in mode-II (right, K$_{II}$) on the reference line with {$L=0.5$ m}.}
\label{fig:EX4_SIF}
\end{figure}

{Using} \eqref{eq:SIF_mode_1} and \eqref{eq:SIF_mode_2}, the SIF 
of each mode on the reference line is presented in Figure~\ref{fig:EX4_SIF}. 
{Compared to \textbf{Example 2} in the pure mode-I loading,} {the maximum values of K$_{I}$ of all cases are on the left end boundary with much greater values, which is} 
due to the aforementioned bottom boundary condition, i.e., the hinge at the bottom as the same as the one in \textbf{Example 3}. 
{Approaching the tip, K$_{I}$ for each case is decreased 
to the negative value, and} 
we confirm the compressive stresses as seen in Figure~\ref{fig:ex4_stress_strain_density_22}. {At the tip,} the case of $\beta=+200$ has the {minimum} as it is stiffer against the compression.
Meanwhile {for K$_{II}$, the case of $\beta=+200$ has its maximum at the tip, whereas the case of $\beta=-200$ shows its maximum about the point of x/L~$=0.5$.} K$_{II}$ {of the case of} $\beta=+200$ is found to surpass the rest cases near  x/L~$=0.7$ on the reference line, 
{confirming more
compressive stress applied
approaching the tip.} 
\\
\newline
\noindent\textit{\textbf{volumetric strain and bulk modulus} }From Figure~\ref{fig:EX4_K_trace} (middle) and (right), we identify that the volumetric strains for all cases are about the same. 
Until about the point of x/L~$=0.6$ on the reference line, positive tr($\epsilon$) implying the dilation of porosity is shown in each case. With about the same change in 
volumetric strain (or implicit porosity) for each case on the reference line,
we note preferentially different and distinctive mechanical responses, {thus depending solely on the nonlinear parameter with $\beta$-values (see \eqref{eq:c1_c2_dd_model} and \eqref{eq:bulk_modulus})}. 
\begin{figure}[h]
\centering
\includegraphics[width=1.00\textwidth]{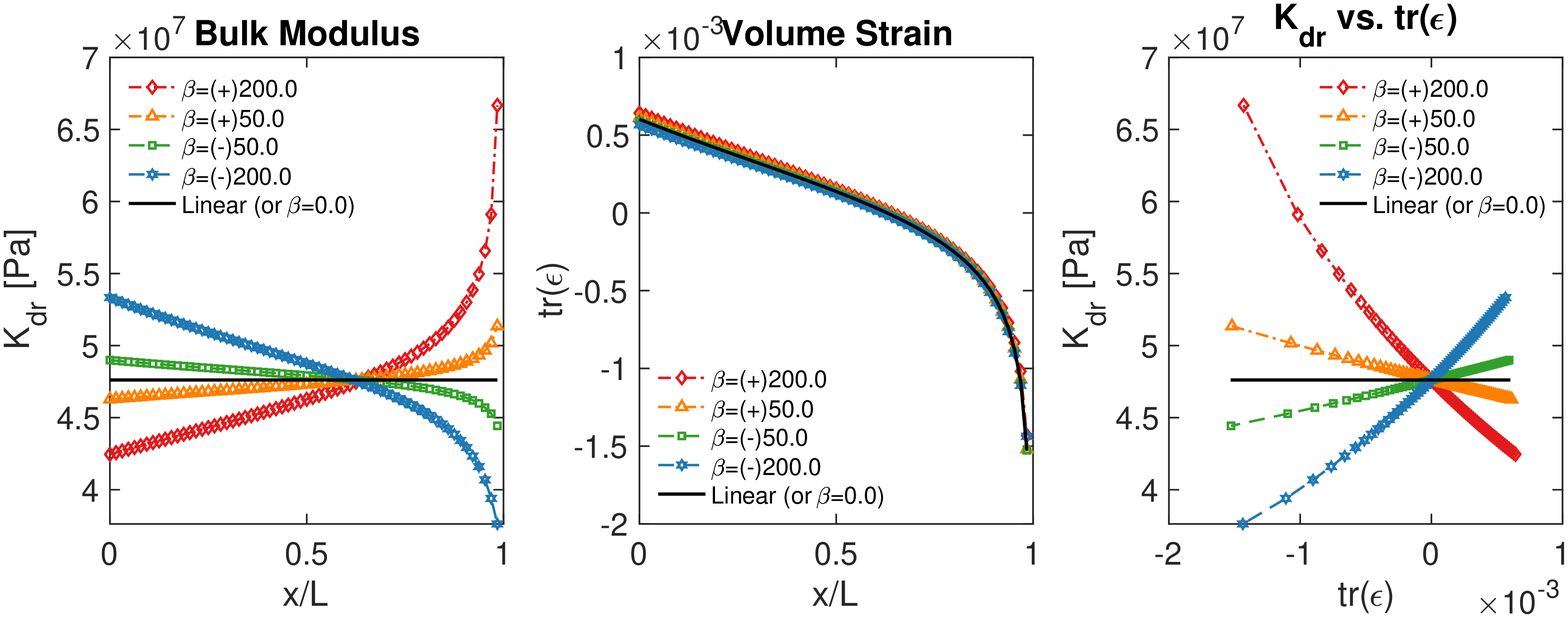}
\caption{[\textbf{Example 4}] Bulk modulus ($K_{dr}${, unit: Pa}) and volumetric strain (tr($\epsilon$)) on the reference line ({$L=0.5$ m}): (left)  $K_{dr}$, (middle) tr($\epsilon$), and (right) $K_{dr}$ vs. tr($\epsilon$). }
\label{fig:EX4_K_trace}
\end{figure}

\section{Conclusion}\label{sec:conclusion}

{The purpose} of this paper is to study 
{{mechanical responses of an elastic porous solid with preferential stiffness} 
whose material moduli are dependent upon the density, {and to provide a stable finite element solution of stress and strain fields, particularly around the crack-tip}.} 
{The proposed model based on the same linearization as the linearized elasticity, i.e., the gradient of the displacement is infinitesimal, cannot stem from the conventional framework of Cauchy elasticity.} {It is structured on a special constitutive relation \eqref{const_relation} for elastic porous solids that show significantly distinctive responses from those using classical linearized models.}
 {For our numerical method, 
we employ an universal and computationally efficient approach of the Newton’s method and FEM under the framework developed in \cite{yoon2021quasi,lee2022finite}} to overcome 
 severe nonlinearity of the model with 
 {partial differential equation}. 
 The proposed algorithm in the study 
 is verified for the optimal convergence rate {using method of manufactured solution}. Three different types of loading are considered in this work. 
{Even though the constitutive relation studied in this paper is simplified via a single nonlinear parameter describing the mechanical response of the material under scrutiny, the distinctive variations in the fields of stress and strain are found {around the crack-tip or damaged pores}.}  Some key findings of this work are:
\begin{itemize}
\item[1. ] In a domain without a crack, the nonlinear \textit{modeling} parameter ``$\beta$''
affects the deformations, 
{which is clearly found in both parallel and perpendicular displacements to the loading directions} along the {reference line}. This parameter clearly controls the strength of a porous material {with the change of volumetric strain; 
the nonlinear parameter $\beta$ with its sign and magnitude is closely related to the preferential stiffness, and when $\beta=0$, the classical linearized model can be recovered.} 
\item[2. ] For a domain with an edge crack under three different types of loading, both crack-tip stress concentration and strain energy density are found to be the maximum 
directly in front of the crack-tip, which is consistent with the observation obtained within the classical linearized elasticity model.
\item[3. ] For a domain with the edge crack under pure tensile loading (or mode-I type of loading),  the crack-tip axial strains are larger with smaller bulk modulus for higher positive $\beta$-values, which implies that the material behaves less stiff (or more compressible). But for lower negative $\beta$-values, its response is weaker against compression. {Thus, under the mode-I tensile loading, the material behavior with {positive $\beta$} can be compared to 
the strain softening, whereas one with the negative $\beta$ to the strain hardening for the elastoplasticity, even though only the pure elastic regime is considered in this study.} 
\item[4. ] In the example with in-plane shear loading (or mode-II type of loading), the density-dependent model for positive $\beta$-values distinctively shows more resistance against compression 
with the higher stress concentration {therein}. 
 \item[5. ] Finally, for the mixed-mode loading (combination of mode-I and mode-II loadings), 
 {{under about the same amount of volumetric strain changes for the cases 
 on the reference line}, {the deviations between the cases are concentrated 
 {still} around the crack-tip, and the positive $\beta$-values show larger in K$_{II}$ but less in K$_I$ in front of the tip.}  
Thus, under the similar type or mixed-mode loading, a certain failure may occur preferentially 
with different material property; for example, the shear failure may occur first for the brittle material {(e.g., rock or ceramic)}. 
{This preferential stiffness can be modeled straightforwardly with the $\beta$-values in the density-dependent material moduli.}} 
\end{itemize}
 
{As the fracture toughness or its propagation is known to be related to the newly generated porosity or damaged pores, 
the model for an elastic porous solid with the preferential material property investigated in this paper can be expanded to study a (quasi-static) crack evolution in the porous material via an appropriate numerical approach {such as the} regularized phase-field approach \cite{yoon2021quasi,lee2022finite}.}  
{Another important next topic can include {bridging the \textit{modeling} parameter to the} \textit{material} parameter and  identifying its role {for the  {characterization}} 
via comparing with some experimental data.}

\bibliographystyle{unsrtnat}
\bibliography{references}

\end{document}